\documentclass{gtart}
\input gtmonout
\volumenumber{1}
\volumeyear{1998}
\volumename{The Epstein birthday schrift}
\pagenumbers{181}{248}
\received{17 November 1997}\revised{29 November 1998}
\published{4 December 1998}
\papernumber{10}

\usepackage{amssymb,rlepsf} 
\def\Box{\sq}
\def\qua{\stdspace}

\newcommand{\maths}[1]{{\Bbb #1}}

\title{Simplicit\'e de groupes d'automorphismes \\
d'espaces \`a courbure n\'egative}
\asciititle{Simplicite de groupes d'automorphismes
d'espaces a courbure negative}
\shorttitle{Simplicit\'e de groupes d'automorphismes}
\covertitle{Simplicit\noexpand\'e de groupes d'automorphismes \\
d'espaces \noexpand\`a courbure n\noexpand\'egative}

\author{Fr\'ed\'eric Haglund\\Fr\'ed\'eric Paulin}
\shortauthors{Fr\'ed\'eric Haglund et Fr\'ed\'eric Paulin}
 
\asciiauthors{Frederic Haglund et Frederic Paulin}
\coverauthors{Fr\noexpand\'ed\noexpand\'eric 
Haglund\\Fr\noexpand\'ed\noexpand\'eric Paulin}

\address{Laboratoire de Topologie et Dynamique URA 1169 CNRS\\
Universite Paris-Sud\\
Bat. 425 (Mathematiques)\\91405 ORSAY Cedex, FRANCE}

\email{haglund@math.u-psud.fr, Frederic.Paulin@math.u-psud.fr}

\theoremstyle{definition}
\newtheorem{defi}{D\'efinition}[section]
\theoremstyle{plain}
\newtheorem{prop}[defi]{Proposition}
\newtheorem{theo}[defi]{Th\'eor\`eme}
\newtheorem{conj}[defi]{Conjecture}
\newtheorem{lemm}[defi]{Lemme}
\newtheorem{coro}[defi]{Corollaire}
\theoremstyle{definition}
\newtheorem{rema}[defi]{Remarque}
\newtheorem{exem}[defi]{Exemple}

\newcommand{\bdefi}{\begin{defi}}
\newcommand{\edefi}{\end{defi}}
\newcommand{\bprop}{\begin{prop}}
\newcommand{\eprop}{\end{prop}}
\newcommand{\btheo}{\begin{theo}}
\newcommand{\etheo}{\end{theo}}
\newcommand{\blemm}{\begin{lemm}}
\newcommand{\brema}{\begin{rema}}
\newcommand{\erema}{\end{rema}}
\newcommand{\bexer}{\begin{exem}}
\newcommand{\eexer}{\end{exem}}
\newcommand{\bconj}{\begin{conj}}
\newcommand{\econj}{\end{conj}}
\newcommand{\elemm}{\end{lemm}}
\newcommand{\bcoro}{\begin{coro}}
\newcommand{\ecoro}{\end{coro}}
\newcommand{\dem}{\proof[Preuve]}
\newcommand{\rem}{\medskip\noindent{\bf Remarque}\qua\ignorespaces}

\newcommand{\A}{{\cal A}}

\newcommand{\T}{{\cal T}}

\newcommand{\M}{{\cal M}}
\newcommand{\G}{{\cal G}}

\newcommand{\K}{{\cal K}}
\renewcommand{\H}{{\cal H}}
\renewcommand{\O}{{\cal O}}
\newcommand{\C}{{\cal C}}
\newcommand{\ra}{\rightarrow}

\newcommand{\RR}{\maths{R}}
\newcommand{\NN}{\maths{N}}

\renewcommand{\SS}{\maths{S}}

\newcommand{\HH}{\maths{H}}
\newcommand{\FF}{\maths{F}}
\newcommand{\ZZ}{\maths{Z}}

\newcommand{\XX}{\maths{X}}

\newcommand{\ssi}{si et seulement si }
\newcommand{\psb}{pre\-mi\`e\-re sub\-di\-vi\-sion bary\-cen\-tri\-que}
\newcommand{\cqfd}{\endproof}
\newcommand{\bord}{\partial}

\begin{document}
 
\begin{abstract}
We prove that numerous negatively curved simply connected locally
compact polyhedral complexes, admitting a discrete cocompact group of
automorphisms, have automorphism groups which are locally compact,
uncountable, non linear and virtually {\bf simple}. Examples include
hyperbolic buildings, Cayley graphs of word hyperbolic Coxeter
systems, and generalizations of cubical complexes, that we call {\it
even} polyhedral complexes. We use tools introduced by Tits in the
case of automorphism groups of trees, and Davis--Moussong's geometric
realisation of Coxeter systems.

{\bf R\'esum\'e}\qua 
Nous montrons que de nombreux complexes
poly\'edraux simplement connexes, localement compacts, \`a courbure
n\'egative, admettant un groupe discret cocompact
d'au\-to\-mor\-phismes, ont leur groupe d'automorph\-is\-mes
localement compact, non d\'enombrable, non lin\'eaire et virtuellement
{\bf simple}. Parmi les exemples, certains sont des immeubles
hyperboliques, des graphes de Cayley de syst\`emes de Coxeter
hyperboliques au sens de Gromov, et des g\'en\'eralisations de
complexes cubiques, que nous appelons des complexes poly\'edraux {\it
pairs}.  Nous utilisons des outils dus \`a Tits dans le cas des
groupes d'automorphismes d'arbres, et la r\'ealisation g\'eom\'etrique
de Davis--Moussong des syst\`emes de Coxeter.
\end{abstract}

\asciiabstract{%
We prove that numerous negatively curved simply connected locally
compact polyhedral complexes, admitting a discrete cocompact group of
automorphisms, have automorphism groups which are locally compact,
uncountable, non linear and virtually simple. Examples include hyperbolic
buildings, Cayley graphs of word hyperbolic Coxeter systems, and
generalizations of cubical complexes, that we call even
polyhedral complexes. We use tools introduced by Tits in the case of
automorphism groups of trees, and Davis-Moussong's geometric
realisation of Coxeter systems.}


 
\maketitle

\section{Introduction}

J.~Tits a d\'emontr\'e dans \cite{Tit} que le groupe des
automorphismes (sans inversion) d'un arbre (diff\'erent de la droite)
homog\`ene ou semi-homog\`ene localement fini, est localement compact,
non d\'enombrable et simple. Le but de cet article est de d\'emontrer
la simplicit\'e de groupes d'automorphismes de nombreux complexes
poly\'edraux localement finis, ayant des propri\'et\'es de courbure
n\'egative, comme par exemple des immeubles hyperboliques ou des
complexes cubiques.

Un {\it immeuble hyperbolique} (voir \cite{GP}) est un immeuble de
type un syst\`eme de Coxeter $(W(P),S(P))$ de la forme suivante.  Soit
$P$ un poly\`edre (compact convexe, pas forc\'ement un simplexe) de
l'espace hyperbolique r\'eel $\HH^n$ de dimension $n$, avec $P$ {\it
de Coxeter} (i.e.~ses angles di\`edres sont de la forme
$\frac{\pi}{k}$ avec $k$ un entier au moins $2$). Alors $S(P)$ est
l'ensemble des r\'eflexions (orthogonales) sur les faces de
codimension $1$ de $P$, et $W(P)$ le groupe d'isom\'etries de $\HH^n$
engendr\'e par $S(P)$.

Un premier exemple est l'immeuble de Bourdon $I_{p,q}$ avec $p\geq
5,q\geq 3$, qui est l'unique complexe poly\'edral de dimension $2$,
dont les polygones sont des copies du $p$--gone hyperbolique r\'egulier
\`a angles droits $P_p$, et le link de chaque sommet est isomorphe au
graphe biparti complet \`a $q+q$ sommets (voir \cite{Bourd2}).  Il
existe une num\'erotation des ar\^etes de $I_{p,q}$ (unique une fois
num\'erot\'ees les ar\^etes d'un polygone fix\'e) par $I=\{1,\cdots
p\}$ de sorte que le long du bord de chaque polygone les ar\^etes
apparaissent avec l'ordre cyclique ou l'ordre inverse.  L'ensemble des
polygones de $I_{p,q}$ est alors un syst\`eme de chambres sur $I$,
deux chambres \'etant $i$--adjacentes si et seulement si les polygones
correspondants se rencontrent le long d'une ar\^ete num\'erot\'ee
$i$. Il est facile (voir \cite{GP}) de montrer que $I_{p,q}$ est un
immeuble de type $(W(P_p), S(P_p))$.

\btheo\label{theo:intro_bourdon} 
Le groupe des automorphismes
pr\'eservant le type de l'im\-meu\-ble de Bourdon $I_{p,q}$ est un
groupe localement compact, non d\'enombrable, non lin\'e\-ai\-re au
moins si $p$ est multiple de $4$, et simple.  
\etheo

Dans \cite{Hag} sont construits de nombreux autres exemples.  Soit $L$
un $m$--gone g\'en\'eralis\'e fini \'epais classique (i.e.~un graphe
biparti complet \`a $p+q$ sommets avec $p,q\geq 3$ si $m=2$, ou si
$m\geq 3$, l'immeuble sph\'erique de rang $2$ d'un groupe de Chevalley
fini $\underline G(\FF_q)$, avec $\underline G$ un groupe alg\'ebrique
simple, de groupe de Weyl le groupe di\'edral $D_{2m}$ d'ordre
$2m$). Par exemple, $L$ peut \^etre l'immeuble des drapeaux du plan
projectif sur le corps fini $\FF_q$, avec $m=3$.  Soit $k$ un entier
pair au moins $6$. Alors dans \cite{Hag} est construit un $2$--complexe
poly\'edral $A_{k,L}$, dont les polygones sont des copies du $k$--gone
hyperbolique $P_{k,m}$ r\'egulier \`a angles $\frac{\pi}{m}$, et le
link de chaque sommet est isomorphe au graphe biparti $L$. L'ensemble
de ses polygones poss\`ede aussi une structure naturelle d'immeuble de
type $(W(P_{k,m}), S(P_{k,m}))$ (voir \cite{GP}).

\btheo\label{theo:intro_Akl}
Le groupe des automorphismes de l'im\-meu\-ble 
hyperbolique $A_{k,L}$ est un groupe lo\-ca\-le\-ment compact, non 
d\'enombrable, non lin\'eaire au moins si $k$ est multiple de $4$, 
et virtuellement simple. 
\etheo

En fait, $A_{k,L}$ est la r\'ealisation g\'eom\'etrique au sens de
Davis--Moussong (voir \cite{Mou}) du syst\`eme de Coxeter $W(k,L)$, dont
la matrice de Coxeter est la matrice d'adjacence du graphe $L$, o\`u
les $1$ et $0$ ont \'et\'e remplac\'es par des $\frac{k}{2}$ et
$\infty$ respectivement. Le $1$--squelette de la r\'ealisation
g\'eom\'etrique de Davis--Moussong d'un syst\`eme de Coxeter $(W,S)$
s'identifie au graphe de Cayley de $(W,S)$, et nous montrons que tout
automorphisme du 1--squelette s'\'etend \`a cette r\'ealisation
g\'eom\'etrique (voir section \ref{sect:groupe_auto_polyedre_pair}).
Appelons {\it mur} du graphe de Cayley l'ensemble des points fixes 
d'un conjugu\'e d'un \'el\'ement de $S$. Un mur est {\it propre}
si aucune des deux composantes du compl\'ementaire du mur ne reste \`a 
distance born\'ee du mur. Un automorphisme du graphe de Cayley
fixe {\it strictement} un mur s'il fixe le mur et n'\'echange pas les
deux composantes de son compl\'ementaire.

Un syst\`eme de Coxeter est dit {\it rigide} s'il n'existe pas
d'automorphisme non trivial de son diagramme qui fixe les ar\^etes de
poids fini issues d'un de ses sommets. G\'en\'eralisant \cite{Hag},
nous montrons que le groupe des automorphismes (de graphe) du graphe
de Cayley de $(W,S)$ est non d\'enombrable si et seulement si $(W,S)$
est non rigide, et nous le calculons exactement dans le cas rigide
(voir th\'eor\`eme \ref{theo:coxeter_non_rigide}).

\btheo\label{theo:intro_coxeter}
Si $(W,S)$ est un syst\`eme de Coxeter, avec $W$ ne contenant pas de
sous-groupe isomorphe \`a $\ZZ+\ZZ$, alors le quotient, par son 
sous-groupe distingu\'e des \'el\'ements fixant l'infini, du 
sous-groupe $G^+$ des automorphismes du graphe de Cayley de $(W,S)$ 
engendr\'e par les fixateurs stricts de murs propres, est simple. 
 Il est non trivial, donc
non d\'enombrable, si et seulement si $(W,S)$ n'est pas rigide.
\etheo

Un {\it complexe cubique} de dimension $n$ est un complexe poly\'edral
$P$, dont les poly\`edres sont des cubes euclidiens
$[-\frac{1}{2},\frac{1}{2}]^k$, tout cube de $P$ \'etant contenu dans
un cube de dimension (maximale) $n$.  Il est dit CAT$(0)$ s'il est
simplement connexe, et si pour tout cube $c$ de $P$, le link $lk(c)$
de $c$ v\'erifie la condition suivante: tout cycle d'ar\^etes dans
$lk(c)$ est de longueur au moins $3$, et si de longueur $3$, borde un
simplexe de $lk(c)$.  Pour toute ar\^ete $a$ de $P$, il existe un
unique sous-complexe (de la subdivision barycentrique) de $P$,
appel\'e {\it mur} (``geometric hyperplane'' par M.~Sageev
\cite{Sag}), rencontrant $a$ en son milieu, et dont toute intersection
non triviale avec un cube de dimension $n$ de $P$ est un hyperplan
$[-\frac{1}{2},\frac{1}{2}]^k\times\{0\}\times
[-\frac{1}{2},\frac{1}{2}]^{n-k-1}$ de ce cube.
Par exemple, si $n=1$, alors $P$ est un arbre, et un mur est le milieu
d'une ar\^ete.

Nous introduisons une notion de {\it poly\`edre pair} (section
\ref{sect:polyedre_pair}) et donc de {\it complexe poly\'edral pair}
(i.e.~dont tous les poly\`edres sont pairs), g\'en\'eralisant
strictement celle de cube et complexe cubique, avec ses murs. Un
poly\`edre d'un espace \`a courbure constante est pair s'il est
sym\'etrique par rapport \`a l'hyperplan m\'ediateur de chacune de ses
ar\^etes, et si un tel hyperplan ne passe pas par un de ses sommets.
Nous donnons en section \ref{sect:polyedre_pair} la construction
explicite de tous les poly\`edres pairs euclidiens ou hyperboliques,
\`a partir des syst\`emes de Coxeter finis, ainsi que la liste
compl\`ete des poly\`edres hyperboliques pairs de dimension 2 et 3 qui
sont eux-m\^emes des poly\`edres de Coxeter. M.~Davis nous a signal\'e
que nos poly\`edres pairs sont, du point de vue combinatoire,
exactement les {\it zonotopes de Coxeter} (aussi appel\'es ``Coxeter
cell'' dans \cite{Dav3}), i.e.~les poly\`edres duaux de l'arrangement
d'hyperplans form\'e par les hyperplans fixes des conjugu\'es des
r\'eflexions d'un syst\`eme de Coxeter fini. Nos complexes
poly\'edraux pairs sont donc, du point de vue combinatoire, des cas
particuliers de ``zonotopal cell complex'' au sens de \cite{DJS}.
Notons qu'il existe des poly\`edres pairs non isom\'etriques ayant
m\^eme combinatoire.

\btheo \label{theo:intro_pair}
Soit $P$ un complexe poly\'edral pair (par exemple cubique),
localement fini, {\rm CAT}$(0)$, admettant un groupe discret cocompact
d'auto\-mor\-phis\-mes qui est hyperbolique au sens de Gromov.  Alors 
le groupe d'automorphismes $G^+$ de $P$ engendr\'e par les fixateurs 
stricts de murs propres est presque simple (au sens que
tout \'eventuel sous-groupe distingu\'e propre est relativement compact).
Si $P$ est {\rm CAT}$(-1)$ et tout point de $P$ appartient \`a une droite  
g\'eod\'esique, alors $G^+$ est simple, et non d\'enombrable si non 
trivial.  
\etheo

Bien s\^ur, $G^+$ peut \^etre trivial. Pour tout type de poly\`edre
euclidien pair possible, nous construisons (section
\ref{sect:exemples_pairs}) un complexe poly\'edral pair CAT$(-1)$,
dont les cellules maximales sont de ce type, et dont le groupe $G^+$
est non d\'enombrable.  Un arbre homog\`ene ou semi-homog\`ene
localement fini admet un groupe discret cocompact d'automorphismes qui
est libre, donc hyperbolique au sens de Gromov (voir section
\ref{sect:rappel} pour des rappels sur cette notion.)  Nous retrouvons
ainsi le r\'esultat de J.Tits. La condition de locale finitude n'est
pas vraiment n\'ecessaire (voir section \ref{sect:applications}). La
condition d'hyperbolicit\'e n'est sans doute pas optimale. Mais comme
le montre le cas du produit de deux arbres homog\`enes, il faut une
hypoth\`ese d'irr\'eductibilit\'e sur $P$. Nous renvoyons \`a
\cite{BM} pour un crit\`ere ing\'enieux de simplicit\'e sur les
groupes discrets d'automorphismes du produit de deux arbres.

Une g\'en\'eralisation imm\'ediate du th\'eor\`eme B de Niblo--Reeves 
\cite{NR} est la suivante.

\btheo\label{theo:Niblo-Reeves}
Soit $P$ un complexe poly\'edral pair {\rm CAT}$(0)$ de dimension finie.
Toute action poly\'edrale sur $P$ d'un groupe ayant la propri\'et\'e {\rm (T)} 
de Kazhdan a un point fixe global.
\etheo

Pour g\'en\'eraliser la situation des exemples ci-dessus, nous
introduisons (section \ref{sect:espace_a_murs}) une notion abstraite
d'ensemble discret $X$ muni d'un syst\`eme de {\it murs}, mod\'elisant
les propri\'et\'es de l'ensemble des sommets d'un complexe poly\'edral
cubique (ou pair) CAT$(0)$ et de la famille de ses hyperplans
m\'ediateurs des ar\^etes, ou d'un groupe de Coxeter $W$ muni de sa
famille de murs (voir \cite[page 14]{Ron}).  Dans les sections
\ref{sect:espace_amur_associe_pair} \`a
\ref{sect:hyperbolicite_complexe_polyedral}, nous \'etudions l'espace
\`a murs canoniquement associ\'e \`a un complexe poly\'edral pair.

Sous des hypoth\`eses d'hyperbolicit\'e au sens de Gromov (voir
section \ref{sect:non_elem} pour les propri\'et\'es que nous
utiliserons) du graphe d'incidence de cette famille de murs, nous
montrons (section \ref{sect:simple}) un th\'eor\`eme de simplicit\'e
sur des groupes de bijections de $X$ pr\'eservant le syst\`eme de
murs, v\'erifiant une condition (P) analogue \`a celle introduite par
J.~Tits \cite{Tit} dans le cas des arbres. Le lemme clef \ref
{lem:commutateur_diabolique} sur les commutateurs est analogue au
lemme 4.3 de \cite{Tit}.  Enfin, en section \ref{sect:applications},
nous appliquons ce th\'eor\`eme de simplicit\'e \`a nos exemples.

\medskip {\small Nous remercions F.~Choucroun, pour son expos\'e sur
l'article de J.~Tits, qui a servi de point de d\'epart \`a ce travail,
ainsi que S.~Mozes et M.~Davis.}

\section{Rappels sur les espaces m\'etriques hyperboliques}
\label{sect:rappel}

Nous renvoyons \`a \cite{Gro,GH} pour les d\'efinitions,
r\'ef\'erences, historiques et preuves des propri\'et\'es rappel\'ees
ci-dessous des espaces m\'etriques hyperboliques au sens de Gromov,
\`a \cite{BH,Bourd1} pour celles des espaces m\'etriques CAT$(\chi)$
au sens d'Alexandroff--Topogonov et \`a \cite{Bri} pour celles des
complexes poly\'edraux.  Le lecteur connaisseur peut se ramener
directement \`a la proposition \ref{prop:non_elementaire}.

\subsection{D\'efinitions diverses}

Une {\it g\'eod\'esique} d'un espace m\'etrique $X$ est une
isom\'etrie d'un intervalle $I$ de $\RR$ dans $X$. On parle de {\it
segment, rayon} ou {\it droite g\'eod\'esique} si $I$ est de la forme
$[a,b], [a,+\infty[$ ou $\RR$.  Un espace m\'etrique est {\it
g\'eod\'esique} si par deux de ses points passe un segment
g\'eod\'esique.

Un espace g\'eod\'esique est {\it hyperbolique} (au sens de Gromov)
s'il existe une constante $\delta\geq 0$ (dite {\it constante
d'hyperbolicit\'e}) telle que tout point de tout c\^ot\'e de tout
triangle g\'eod\'esique est \`a distance au plus $\delta$ d'un point
de l'un des deux autres c\^ot\'es.  Un groupe de type fini $G$, muni
d'une partie g\'en\'eratrice $S$, est {\it hyperbolique} (au sens de
Gromov) si le graphe de Cayley de $G$ pour $S$, muni de sa m\'etrique
naturelle, est hyperbolique. Une application $f\co X\ra Y$ entre deux
espaces m\'etriques est une {\it quasi-isom\'etrie} s'il existe des
constantes $\lambda\geq 1,c,c'\geq 0$ telles que pour tous $x,y$ dans
$X$ et $z$ dans $Y$:
$$\frac{1}{\lambda}d(x,y)-c\leq d(f(x),f(y))\leq \lambda d(x,y) + c
\;\;\; {\rm et}\;\;\; d(z,f(X))\leq c'.$$ Un espace g\'eod\'esique 
quasi-isom\'etrique \`a un espace hyperbolique est encore
hyperbolique, donc l'hyperbolicit\'e d'un groupe ne d\'epend pas de la
partie g\'en\'eratrice fix\'ee.

Deux rayons g\'eod\'esiques sont {\it asymptotes} si leur distance de
Hausdorff est finie. Ceci d\'efinit une relation d'\'equivalence sur
l'ensemble des rayons g\'eod\'esiques dans $X$. L'en\-sem\-ble des
classes d'\'equivalence est appel\'e le {\it bord} (ou espace \`a
l'infini) de $X$, et not\'e $\partial X$. Il existe une topologie
naturelle sur $\overline{X}=X\cup\partial X$, m\'etrisable compacte
lorsque $X$ est hyperbolique, localement compact, complet.  Toute
quasi-isom\'etrie entre deux espaces hyperboliques s'\'etend
contin\^ument en un hom\'eomorphisme de $\partial X$ sur $\partial Y$.

Soit $X$ un espace g\'eod\'esique et $\chi\in\RR$. Soit
${\XX}^2_{\chi}$ le plan riemannien complet simplement connexe \`a
courbure constante $\chi$ (${\XX}^2_{\chi}$ est le plan hyperbolique,
le plan euclidien, la sph\`ere de dimension $2$ si $\chi=-1,0,1$).
Soit $\triangle=[xy]\cup[yz]\cup[zx]$ un triangle g\'eod\'esique dans
$X$. Soit $\overline\triangle=[\overline x\,\overline y ]\cup
[\overline y\,\overline z ]\cup[\overline z\,\overline x ]$ un
triangle g\'eod\'esique dans ${\XX}^2_{\chi}$ ayant m\^emes longueurs
des c\^ot\'es que $\triangle$.  Si $s\in\triangle$, le point sur le
c\^ot\'e correspondant de $\overline\Delta$, \`a la m\^eme distance
des extr\'emit\'es que $s$, est not\'e $\overline s$. Un triangle
g\'eod\'esique $\triangle$ dans $X$ est CAT$(\chi)$ s'il est plus
``pinc\'e'' que le triangle correspondant de l'espace mod\`ele,
i.e.~si, pour tous points $s,t\in\triangle$, on a
$$d_X(s,t)\leq d_{{\XX}^2_{\chi}}(\overline s,\overline t).$$ 
Un espace g\'eod\'esique est CAT$(\chi)$ si tout triangle g\'eod\'esique de
$X$ est CAT$(\chi)$. Si $\chi<0$, un espace CAT$(\chi)$ est
hyperbolique au sens de Gromov.

Un {\it complexe poly\'edral} $P$ est un complexe cellulaire (voir par
exemple \cite{Whi}) dont les cellules sont des poly\`edres (compacts
convexes) d'un espace \`a courbure constante, et dont les applications
d'attachements sont cellulaires et localement isom\'etriques sur
chaque cellule ouverte.  Un {\it complexe polygonal} est un complexe
poly\'edral de dimension $2$.  Un complexe poly\'edral, dont les
poly\`edres sont des simplexes ne se rencontrant qu'au plus en une
face, est pr\'ecis\'ement (la r\'ealisation g\'eom\'etrique d') un
complexe simplicial.

Un {\it automorphisme} de complexe poly\'edral de $P$ est un
automorphisme du complexe cellulaire $P$.  Nous identifierons deux
automorphismes qui envoient chaque cellule ouverte sur une m\^eme
cellule ouverte.  Un automorphisme est dit {\it isom\'etrique} (ou une
{\it isom\'etrie poly\'edrale}) si sa restriction \`a chaque
poly\`edre est isom\'etrique. Par exemple, si $P$ est un rectangle
euclidien non carr\'e, alors $P$ admet $4$ isom\'etries poly\'edrales,
et $8$ automorphismes. Si $P$ est muni de la topologie faible usuelle,
le groupe des automorphismes de $P$ sera muni de la topologie
compacte--ouverte. Si $P$ est localement fini, alors Aut$\;G$ est
localement compact, et le fixateur de tout poly\`edre de $P$ est un
groupe compact profini.

Si $P$ n'a qu'un nombre fini de classe d'isom\'etrie de poly\`edres,
alors (voir \cite{Bri}) il existe une m\'etrique $d$ (naturelle pour
les automorphismes de $P$) g\'eo\-d\'e\-si\-que et compl\`ete, ainsi
d\'efinie.  Une {\it g\'eod\'esique bris\'ee} $\gamma$ de $P$ est une
courbe qui, par morceaux, est contenue et g\'eod\'esique dans un
poly\`edre de $P$.  Sa longueur $\ell(\gamma)$ est la somme des
longueurs des morceaux g\'eod\'esiques pr\'ec\'edents. Alors $d(x,y)$
est la borne inf\'erieure des longueurs des g\'eod\'esiques bris\'ees
entre $x$ et $y$.

Sauf mention explicite du contraire, tout complexe poly\'edral sera
muni de cette distance. Toute isom\'etrie poly\'edrale est une
isom\'etrie pour cette distance. La topologie faible et la topologie
induite par cette distance co\"{\i}ncident si et seulement si $P$ est
localement fini. Voir \cite{Gro} pour l'\'equivalence, dans le cas des
complexes cubiques, entre la d\'efinition ci-dessus de CAT$(0)$ et
celle donn\'ee en introduction.

Si $C$ est un complexe poly\'edral n'ayant qu'un nombre fini de types
d'iso\-m\'e\-trie de cellules, et $x\in C$, nous noterons $lk(x,C)$
l'espace des germes de segments g\'eod\'esiques issus de $x$. Il
poss\`ede une structure naturelle de complexe poly\'edral, dont les
cellules sont des poly\`edres sph\'eriques.

Un {\it graphe} est un 1--complexe simplicial connexe. En identifiant
chaque ar\^ete \`a $[-\frac{1}{2},\frac{1}{2}]$, on obtient un
complexe poly\'edral.  Sa m\'etrique est l'unique m\'etrique
g\'eod\'esique rendant chaque ar\^ete isom\'etrique \`a $[0,1]$.  Un
{\it arbre} est un graphe simplement connexe. Un arbre est
CAT$(-\infty)$, i.e.~CAT$(\chi)$ pour tout $\chi\in\RR$.

\subsection{Groupes d'isom\'etries non \'el\'ementaires}
\label{sect:non_elem}

Soit $Y$ un espace m\'etrique complet, g\'eod\'esique et hyperbolique,
tel que par deux points de $Y\cup\partial Y$ passe un segment, rayon
ou droite g\'eod\'esique (cette derni\`ere condition est toujours remplie si
$Y$ est localement compact).  On note $\partial^2 Y$ l'espace des
couples de points distincts de $\partial Y$. On note $\overline Z$
l'adh\'erence dans $Y\cup\partial Y$ d'une partie $Z$ de $Y$, et
$\partial Z=\overline{Z}\cap\partial Y$.

Une isom\'etrie $g$ de $Y$ est dite {\it hyperbolique} si pour un
(donc pour tout) point $x$ dans $Y$, l'application de $\ZZ$ dans $Y$
qui \`a $k$ associe $g^kx$ est une quasi-isom\'etrie sur son image.
En particulier, $g$ admet alors exactement deux points fixes dans
$\partial Y$.

Soit $G$ un sous-groupe du groupe des isom\'etries de $Y$ (n'agissant
peut-\^etre pas proprement discontinument).  D\'efinissons l'{\it
ensemble limite} $\Lambda G$ de $G$ comme l'adh\'erence dans $\partial
Y$ de l'ensemble des points fixes dans $\partial Y$ des \'el\'ements
hyperboliques de $G$.  Le groupe $G$ est dit {\it non \'el\'ementaire}
si son ensemble li\-mite contient au moins trois points et ne
contient pas de point fixe global (cette derni\`ere condition est
toujours remplie si $Y$ est localement compact et $G$ discret). Si $G$ est non
\'el\'ementaire, $\Lambda G$ est non d\'enombrable et sans point
isol\'e; c'est l'ensemble d'accumulation dans $\partial Y$ de l'orbite
par $G$ de tout point de $Y$; c'est le plus petit ferm\'e non vide
invariant par $G$ dans $\partial X$; l'orbite par $G$ de tout point de
$\Lambda G$ est dense dans $\Lambda G$.  On note $\Lambda^2 G$
l'ensemble des couples de points distincts de $\Lambda G$.

\rem 
Par exemple, si $Y$ est localement compact, si $G$ contient un
sous-groupe agissant proprement discontin\^ument avec quotient compact
sur $Y$, alors $G$ est non \'el\'ementaire et $\Lambda G=\partial Y$.

\bprop \label{prop:non_elementaire} 
Si $G$ est non \'el\'ementaire, alors l'ensemble des couples des point
fixes des \'el\'ements hyperboliques de $G$ est dense dans $\Lambda^2 G$.

Soit $H$ un sous-groupe distingu\'e non trivial de $G$.  Si $G$ est
non \'el\'ementaire, alors ou bien $H$ est contenu dans le noyau de
l'action de $G$ sur $\Lambda G$, ou bien $H$ est non \'el\'ementaire,
d'ensemble limite \'egal \`a celui de $G$.  
\eprop

\dem 
La premi\`ere assertion est due \`a \cite[Corollaire 8.2.G]{Gro}.

Pour la seconde assertion, supposons que $h\in H$ n'agisse pas
trivialement sur l'ensemble limite de $G$. Montrons tout d'abord que
$H$ contient au moins un \'el\'ement hyperbolique.

Soit $a\in \Lambda G$ tel que $ha\neq a$. Par invariance, $ha$ est
dans $\Lambda G$.  Soit $\delta$ une constante d'hyperbolicit\'e de
$X$. Soit $U$ un voisinage ouvert suffisamment petit de $a$ dans
$Y\cup \partial Y$, de sorte que $U$ et $hU$ soient disjoints, et
s\'epar\'es d'une distance grande devant $\delta$. Soit $g$ un
\'el\'ement hyperbolique de $G$, dont les points fixes r\'epulsif
$g^-\in \Lambda G$ et attractifs $g^+\in \Lambda G$ sont dans $U$ et
$hU$ respectivement.  Soit $\gamma$ une g\'eod\'esique entre $g^-$ et
$g^+$. Soit $y$ un point de $\gamma\cap U$. En particulier, $hy$
appartient \`a $hU$.

\begin{figure}[htbp]
\cl{\relabelbox\small
\epsfysize 5cm\epsfbox{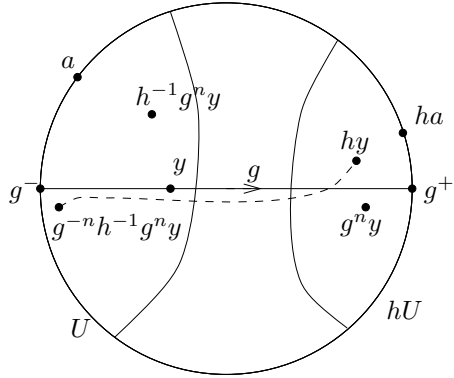}
\relabel {a}{$a$}
\relabel {ha}{$ha$}
\relabel {hgy}{$h^{-1}g^ny$}
\relabel {y}{$y$}
\adjustrelabel <0pt, 1pt> {g}{$g$}
\adjustrelabel <0pt, 1pt> {hy}{$hy$}
\adjustrelabel <-2pt, 0pt> {gm}{$g^-$}
\relabel {gp}{$g^+$}
\adjustrelabel <0pt, -2pt> {gy}{$g^ny$}
\relabel {hU}{$hU$}
\relabel {U}{$U$}
\relabel {ghgy}{$g^{-n}h^{-1}g^ny$}
\endrelabelbox}
\caption{Construction d'un \'el\'ement hyperbolique
  dans $H$}
\end{figure}


Si $n$ est assez grand, alors $g^ny$ est proche de $g^+$, donc
appartient \`a $hU$.  Donc $h^{-1}g^n y$ appartient \`a $U$, et si $n$
est assez grand, $g^{-n} h^{-1}g^n y$ est beaucoup plus proche de
$g^-$ que $y$.  Donc il existe une constante $K$ (ne d\'ependant que
de $\delta$) telle que $y$ est \`a distance au plus $K$ d'un segment
g\'eod\'esique entre $g^{-n} h^{-1}g^n y$ et $hy$. Quitte \`a avoir
pris $U$ suffisamment petit, on a $$\inf\{d(y, hy), d(y,g^{-n}
h^{-1}g^n y)\} > 2K +1000 \delta.$$  Par \cite[Lemma 8.1.A]{Gro}, on en
d\'eduit que $h (g^{-n}h^{-1}g^n)^{-1}$ est hyperbolique. Comme $H$
est distingu\'e, ceci montre notre affirmation pr\'eliminaire.

Maintenant, comme les conjugu\'es d'un \'el\'ement hyperbolique $h$ de
$H$ sont encore dans $H$, que l'orbite par $G$ d'un point fixe de $h$
est contenue et dense dans $\Lambda G$, on en d\'eduit que $\Lambda
H=\Lambda G$. En particulier $\Lambda H$ contient au moins trois
point. Si $H$ fixait un point $a$ de $\Lambda H$, celui-ci serait
unique \cite[8.2.D]{Gro}.  Comme $H$ est distingu\'e dans $G$, le
point $a$ serait fixe par $G$, ce qui est impossible.  \endproof

\blemm\label{lemm:compact_fidele} 
Supposons $\partial Y$ non vide sans point isol\'e.  Si $Y$ est
localement compact, le noyau de l'action de $G$ sur le bord de $Y$ est
relativement compact dans le groupe des isom\'etries de $Y$ (donc
compact si $G$ est ferm\'e dans le groupe des isom\'etries de $Y$).
Si $Y$ est {\rm CAT}$(-1)$ et tout point de $Y$ appartient \`a une droite
g\'eod\'esique, alors $G$ agit fid\`element sur le bord.  
\elemm

\dem 
Pour la premi\`ere assertion, soient $x,y,z$ trois points
distincts de $\partial Y$ et $p$ une quasi-projection de $x$ sur une
g\'eod\'esique entre $y$ et $z$.  Une isom\'etrie de $Y$ qui fixe
(point par point) le bord de $Y$ bouge $p$ d'une distance inf\'erieure
\`a une constante. Le r\'esultat d\'ecoule alors du th\'eor\`eme
d'Ascoli.
 
Pour la seconde assertion, soit $g\in G$ fixant le bord de $Y$.  Soit
$x\in Y$ et $a,b\in \partial Y$ les extr\'emit\'es d'une droite
g\'eod\'esique $D$ passant par $x$. Soient $a', b'$ deux points
proches et distincts de $a,b$ respectivement. Soit $p,p'$ l'unique
projection de $a',b'$ sur $D$. Alors par unicit\'e, $p$ et $p'$ sont
fixes par l'isom\'etrie $g$, et $x\in [p,p']$ aussi, par unicit\'e du
segment g\'eod\'esique entre deux points.  \endproof

\section{Espaces \`a murs}
\label{sect:espace_a_murs}

Soit $X$ un ensemble. Un {\it mur} de $X$ est une partition de $X$ en
deux sous-ensembles, appel\'es les {\it demi-espaces} d\'efinis par le
mur. Un mur {\it s\'epare} deux points $x$ et $y$ de $X$ \ssi\ $x$
appartient \`a l'un des demi-espaces d\'efinis par le mur et $y$
appartient \`a l'autre. Un {\it syst\`eme de murs sur $X$} est un
ensemble $\M$ de murs de $X$ tel que:
\begin{center}
  \parbox{12cm}{(M)\qua Pour tous $x$ et $y$ distincts dans $X$,
    l'ensemble $\M(x,y)$ des murs de $\M$ s\'eparant $x$ et $y$ est
    fini non vide.}
\end{center}
Un {\it espace \`a murs} est un couple $(X,\M)$, o\`u $X$ est un
ensemble et $\M$ un syst\`eme de murs sur $X$. Tout singleton de $X$
est alors l'intersection des demi-espaces qui le contiennent.

Dans un espace \`a murs $(X,\M)$, on dit qu'un point $z$ est {\it
entre} deux points $x$ et $y$ si $\M(x,y)$ est la r\'eunion
(n\'ecessairement disjointe) de $\M(x,z)$ et $\M(z,y)$. Le {\it graphe
associ\'e} \`a $(X,\M)$ est le graphe ayant $X$ pour ensemble de
sommets, et une ar\^ete entre deux sommets $x$ et $y$ \ssi\ les seuls
points de $X$ entre $x$ et $y$ sont $x$ et $y$.  On note $\G=\G(X,\M)$
ce graphe, qui est connexe d'apr\`es l'axiome (M). Un mur $M$ de $X$
est dit {\it transverse} \`a une ar\^ete de $\G(X,\M)$ lorsqu'il
s\'epare ses extr\'emit\'es.

Un espace \`a murs $(X,\M)$ est dit {\it hyperbolique} si son
graphe associ\'e est un espace m\'etrique hyperbolique au sens de
Gromov, et s'il v\'erifie la condition (H) suivante de non
trivialit\'e et de compatibilit\'e entre la structure m\'etrique de
$\G$ et le syst\`eme de demi-espaces d\'efini par $\M$:
\begin{center}
  \parbox{12cm}{(H)\qua
    Pour tout $\xi\in\partial \G$, l'ensemble des parties de
    $\G\cup\partial \G$ de la forme $\overline A$, o\`u $A$ est un
    demi-espace de $(X,\M)$ tel que $\overline A$ contient $\xi$ dans
    son int\'erieur, est une base de voisinages de $\xi$ dans
    $\G\cup\partial \G$.  }
\end{center}

\subsection{Automorphismes d'espaces \`a murs et propri\'et\'e (P) de Tits}
\label{sect:auto_P}

Soit $(X,\M)$ un espace \`a murs.  Un {\it automorphisme} $\phi$ de
$(X,\M)$ est une bijection de $X$ pr\'eservant $\M$. Il induit un
automorphisme du graphe $\G$, encore not\'e $\phi$. Si $(X,\M)$ est
hyperbolique, alors $\phi$ induit un hom\'eomorphisme du bord
hyperbolique $\partial \G$ de $\G$, toujours not\'e $\phi$.

Si Aut$(X,\M)$ est le groupe des automorphismes de $(X,M)$, et
Aut$(\G)$ le groupe des automorphismes de graphe de $\G$, alors
l'application $\phi\mapsto\phi$ est une injection de Aut$(X,\M)$ dans
Aut$(\G)$, en g\'en\'eral non surjective (voir toutefois la preuve du
th\'eor\`eme \ref{theo:meme_groupe_automorphisme}). Nous identifierons
Aut$(X,\M)$ avec son image dans Aut$(\G)$. Lorsque $\G$ est localement
fini, nous munirons Aut$(\G)$ de la topologie compacte--ouverte et
Aut$(X,\M)$ de la topologie induite.

Un automorphisme {\it fixe strictement} un mur $M$ s'il fixe les
sommets de toute ar\^ete transverse \`a $M$.  Un automorphisme d'un
espace \`a murs {\it fixe strictement} un demi-espace $A$ s'il fixe
$A$ et fixe strictement le mur $M=\{\A,X\setminus A\}$.

\blemm \label{lem:fixe_strictement} 
Un automorphisme fixant strictement un mur $M$ pr\'eserve chacun des
demi-espaces de $X$ d\'efinis par $M$.  
\elemm

\dem 
Remarquons d'abord que si $M$ s\'epare deux points $x,y$, alors
tout chemin entre $x$ et $y$ dans $\G$ contient une ar\^ete de $\G$
transverse \`a $M$.

Notons $M=\{A,X\setminus A\}$ et $V(M)$ l'ensemble des sommets
d'ar\^etes de $\G$ transverses \`a $M$.  Si $x$ appartient au
demi-espace $A$, soit $p$ un point de $V(M)$ \`a distance minimale de
$x$.  Par minimalit\'e, $p$ est dans $A$.  Si $\phi$ fixe strictement
$M$, alors il fixe point par point $V(M)$. Il envoie un chemin
$\gamma$ de longueur minimale entre $x$ et $p$ sur un chemin de m\^eme
longueur entre $\phi(p)=p$ et $\phi(x)$.  Si $\phi(x)$ n'est pas dans
$A$, alors le chemin $\phi(\gamma)$ doit contenir une ar\^ete
transverse \`a $M$, ce qui contredit le fait que $\phi$ pr\'eserve la
distance combinatoire \`a $V(M)$.  \endproof

On appelle {\it cha\^{\i}ne} une suite $(A_i)_{i\in \ZZ}$ de
demi-espaces qui est strictement d\'ecrois\-sante pour l'inclusion.  Un
automorphisme {\it fixe strictement} cette cha\^{\i}ne s'il fixe
strictement chaque mur $M_i=\{A_i,X\setminus A_i\}$. Par le lemme
pr\'ec\'edent, il pr\'eserve alors chaque demi-espace $A_i$.

Soit $G$ un groupe d'automorphismes de $(X,\M)$. Si $M=\{A,X\setminus
A\}$ est un mur de $\M$, soit $G_M$ le sous-groupe de $G$ fixant
strictement $M$.  Par le lemme pr\'ec\'edent, le groupe $G_M$
pr\'eserve les ensembles $X\setminus A$ et $A$. Nous notons $G_A$
(resp.~$G_{X\setminus A}$) le groupe des permutations de $A$
(resp.~$X\setminus A$) induit par $G_M$. Le produit des restrictions
donne un morphisme injectif
$$G_M\ra G_{A}\times G_{X\setminus A}.$$ Soit $C=(A_i)_{i\in \ZZ}$ une
cha\^{\i}ne. Soit $G_C$ le sous-groupe de $G$ fixant strictement $C$.
Pour tout $i$, le groupe $G_C$ pr\'eserve l'ensemble $A_{i}\setminus
A_{i+1}$, et nous notons $G_{C,i}$ le groupe des permutations de cet
ensemble induit par $G_C$. Le produit direct des restrictions $G_C\ra
G_{C,i}$ est un morphisme
$$G_C\ra\prod_{i\in\ZZ} G_{C,i}.$$ 

\blemm\label{lem:cemorphisme_injectif}
Ce morphisme est injectif.
\elemm

\dem 
Il suffit de montrer que pour toute cha\^{\i}ne $C=(A_i)_{i\in \ZZ}$
de $(X,\M)$, la r\'eunion des $A_i\setminus A_{i+1}$ vaut tout $X$.
Supposons par l'absurde qu'il existe un point $x$ n'appartenant pas \`a
cette r\'eunion. Supposons que $x$ appartient \`a $A_0$ (si $x\in
X\setminus A_0$, le raisonnement est le m\^eme, quitte \`a renverser
l'ordre de $\ZZ$).  Soit $x_0\in A_0\setminus A_1$. Alors $x_0$
appartient \`a $X\setminus A_i$ et $x$ appartient \`a $A_i$ pour tout
$i\geq 1$.  Donc le mur $M_i=\{X\setminus A_i,A_i\}$ s\'epare $x_0$ et
$x$ pour tout $i\geq 1$, ce qui contredit la finitude de $\M(x,x_0)$.
\endproof

La d\'efinition suivante est alors analogue \`a la propri\'et\'e
homonyme de \cite{Tit}.

\bdefi 
On dit qu'un groupe $G$ d'automorphismes de $(X,\M)$ {\it
  v\'erifie la propri\'et\'e} $(P)$ si pour tout mur $M$ et toute
cha\^{\i}ne $C$, les morphismes pr\'ec\'edents sont surjectifs,
i.e.~des isomorphismes.  
\edefi

\blemm\label{fixateur_mur_stabilisateur_demiespace} 
Soit $G$ un groupe d'automorphismes d'un espace \`a murs, ayant la
propri\'et\'e {\rm(P)}.  Alors le sous-groupe de $G$ engendr\'e par les
fixateurs stricts de murs co\"{\i}ncide avec le sous-groupe de $G$
engendr\'e par les fixateurs stricts de demi-espaces.  
\elemm

\dem 
Le second groupe est contenu dans le premier, par d\'efinition.
Il suffit donc de montrer que tout \'el\'ement $g$ de $G$ fixant
strictement un mur $M=\{A^-,A^+\}$ est produit de deux \'el\'ements
$g^-,g^+$ fixant strictement les demi-espaces $A^-,A^+$
respectivement.  Par la propri\'et\'e (P), le morphisme $G_M\ra
G_{A^-}\times G_{A^+}$ est surjectif. Il suffit de prendre pour
$g^-,g^+$ des pr\'eimages de $(g|_{A^-},id)$ et $(id,g|_{A^+})$
respectivement.  
\endproof

Consid\'erons la propri\'et\'e suivante d'un espace \`a murs $(X,\M)$.
\begin{center}
  \parbox{12cm}{(${\rm M}'$)\qua Pour tous demi-espaces $A,B$ de $(X,\M)$, avec
    $B$ rencontrant $A$ et son compl\'ementaire, tout automorphisme
    fixant strictement le mur $M=\{A,X\setminus A\}$ pr\'eserve $B$.
    }
\end{center}
Dans le cas d'un arbre, cette condition est vide (donc n'appara\^{\i}t
pas dans \cite{Tit}).

\blemm\label{Mprime_entraine_P} 
Si un espace \`a murs $(X,\M)$ v\'erifie la condition {\rm(${\rm M}'$)}, alors le
groupe de tous ses automorphismes v\'erifie la propri\'et\'e {\rm(P)}.
\elemm

\dem 
Soit $M=\{A^-,A^+\}$ un mur de $(X,\M)$. Soit $h^{\pm}$ la restriction
\`a $A^{\pm}$ d'un automorphisme $\overline{h}^{\,\pm}$ de $(X,\M)$
fixant strictement $M$. Comme $A^-\cup A^+=X$, soit $g$ la bijection
de $X$ valant $h^{\pm}$ sur $A^{\pm}$. Montrons que $g$ pr\'eserve
$\M$, ce qui impliquera la surjectivit\'e de Aut$(X,\M)_M\!\ra
$Aut$(X,\M)_{A^-}\!\times $Aut$(X,\M)_{A^+}$.  Soit
$N=\{B,X\setminus B\}$ un mur de $(X,\M)$. Si $B$ est contenu dans
$A^{\pm}$, alors $g(B)=h^{\pm}(B)\subset A^{\pm}$, donc
$g(B)=\overline{h}^{\,\pm}(B)$ est un demi-espace de $(X,\M)$. D'o\`u
$g(N)$ est encore un mur de $(X,\M)$.  Si $B$ rencontre \`a la fois
$A^{-}$ et $A^{+}$, alors les deux automorphismes $\overline{h}^{-}$
et $\overline{h}^{+}$ pr\'eservent $B$ par la propri\'et\'e (${\rm M}'$). Donc
$\overline{h}^{\,\pm}$ pr\'eserve $B\cap A^{\pm}$.  D'o\`u $g$
pr\'eserve $B$, et $g(N)=N$ est encore un mur de $(X,\M)$.

Soit $C=(A_i)_{i\in \ZZ}$ une cha\^{\i}ne de $(X,\M)$, et soit $h_{i}$
la restriction \`a $A_{i}\setminus A_{i+1}$ d'un automorphisme
$\overline{h}_{i}$ de $(X,\M)$ fixant strictement $M_i$.  
Comme $X=\bigcup_{i\in\ZZ} A_i\setminus A_{i+1}$ (voir la 
preuve du lemme \ref{lem:cemorphisme_injectif}), il existe une
bijection $g$ de $X$ valant $h_{i}$ sur $A_{i}\setminus A_{i+1}$. Soit
$B$ un demi-espace de $(X,\M)$. On montre comme pr\'ec\'edemment que
si $B$ est contenu dans un $A_{i}\setminus A_{i+1}$, alors $g(B)$ est
encore un demi-espace, et que, par la propri\'et\'e (${\rm M}'$), si $B$
rencontre au moins deux $A_{i}\setminus A_{i+1}$, alors $g(B)=B$. Donc
$g$ est un automorphisme de $(X,\M)$. Ceci montre la surjectivit\'e de
Aut$\,(X,\M)_C\ra \prod_{i\in\ZZ}\,$Aut$\,(X,\M)_{C,i}$.  \endproof

Un mur d'un espace \`a murs hyperbolique est dit {\it propre} si le
bord \`a l'infini dans $\overline \G$ de chacun des demi-espaces qu'il
d\'efinit n'est pas \'egal \`a tout $\partial \G$. Une cha\^{\i}ne
$C=(A_i)_{i\in \ZZ}$ est {\it propre} si chaque mur
$M_i=\{A_i,X\setminus A_i\}$ est propre. Dans la condition (H), nous
pouvons de plus supposer que les murs d\'efinissant les demi-espaces
$A$ sont propres. Si $G$ est un groupe d'automorphismes de $(X,\M)$,
nous noterons $G^+$ le sous-groupe de $G$ engendr\'e par les fixateurs
stricts de murs propres.

\blemm\label{lem:nontrivial_nondenombrable} 
Soit $(X,\M)$ un espace \`a murs hyperbolique, de graphe associ\'e
$\G$ localement fini, et $G$ un groupe d'au\-to\-mor\-phis\-mes de
$(X,\M)$, ferm\'e vu comme sous-groupe du groupe des automorphismes de
$\G$, ayant la propri\'et\'e {\rm(P)}, agissant de mani\`ere non
\'el\'ementaire sur $\G$ et d'ensemble limite \'egal \`a $\partial
\G$. Si $G^+$ est non trivial, alors $G^+$ est non d\'enombrable.
\elemm
  
\dem 
Soit $M$ un mur propre, de fixateur strict non trivial. Soit
$G^{+}_M$ le sous-groupe de $G$ fixant strictement le mur $M$ et
fixant l'un des demi-espaces, disons $A$, d\'efinis par $M$. Le
sous-groupe $G^{+}_M$ est ferm\'e dans $G$, donc dans Aut$(\G)$. On en
d\'eduit que $G^{+}_M$ est localement compact. Pour montrer qu'il est
non d\'enombrable, il suffit de montrer qu'il n'a pas de point
isol\'e, et comme c'est un groupe topologique, que l'identit\'e n'est
pas isol\'ee.

Soit $g$ un \'el\'ement non trivial de $G^{+}_M$, qui existe par la
propri\'et\'e (P) quitte \`a \'echanger $A$ et $X\setminus A$, et $K$
une partie compacte arbitraire de $\G$. Puisque $M$ est propre, soit
$x$ un point de $\partial X\setminus \partial A$. Soit $U$ un ouvert,
contenu dans $\overline{X}\setminus(\overline {A}\cap K)$, contenant
$x$.  Puisque $G$ est non \'el\'ementaire, il existe un \'el\'ement
hyperbolique $h$ dans $G$ dont le point fixe attractif est contenu
dans $U$ et le point fixe r\'epulsif dans $\partial X\setminus
\partial(X\setminus A)$.  Si $n$ est assez grand, alors
$h^n(\overline{X\setminus A})$ est contenu dans $U$. Posons
$g_n=h^ngh^{-n}$, qui appartient \`a $G$ et m\^eme \`a $G^{+}_M$.
Comme $g$ vaut l'identit\'e sur $A$, l'\'el\'ement $g_n$ vaut
l'identit\'e sur $h^n(A)$, donc sur $K$.  Puisque $g$ est non trivial,
$g_n$ l'est aussi. On en d\'eduit que l'identit\'e n'est pas isol\'ee
dans $G^{+}_M$.  
\cqfd

\subsection{L'exemple classique des syst\`emes de Coxeter}
\label{sect:exemple_syssteme_coxeter}

Adoptons un premier point de vue alg\'ebrique (on trouvera dans
\cite[Chapitre~IV, Section 1, Exemple~16]{Bourb}, \cite{Ron} toutes les
justifications des affirmations ci-dessous).  Soient $(W,S)$ un
syst\`eme de Coxeter, $\T$ l'ensemble de ses {\it r\'eflexions}
(i.e.~des conjugu\'es dans $W$ des \'el\'ements de $S$), et $\ell (w)$
la longueur minimale d'une \'ecriture de $w\in W$ comme mot sur $S$.
Pour $t\in \T$, posons:
$$A^+_t=\{ w\in W, \ell (w) <\ell (tw)\} \ {\rm ~et~}\ A^-_t=\{ w\in
W, \ell (w) >\ell (tw)\} .$$
Alors $A^+_t$ contient $1_W$ et $A^-_t$
contient $t$. De plus, $\ell (w)$ et $\ell (tw)$, n'ayant pas la
m\^eme parit\'e, sont toujours diff\'erents.  Donc $\{ A^+_t,A^-_t\}$
est un mur de $W$ (les demi-espaces $A^\pm_t$ sont appel\'ees
moiti\'es dans \cite{Bourb}). Notons $\M(W,S)=\M$ l'ensemble des murs
ainsi obtenus (en correspondance biunivoque avec $\T$).  Montrons que
$\M(W,S)$ v\'erifie l'axiome (M).

Pour $w',w''\in W$, l'ensemble des murs s\'eparant $w'$ de $w''$
correspond \`a l'ensem\-ble des r\'eflexions $t$ de la forme $s_1\ldots
s_{i-1}s_is_{i-1}\ldots s_1$, pour une \'ecriture g\'eod\'es\-ique
fix\'ee $w'^{-1}w''=s_1\ldots s_n$, avec $n=\ell (w'^{-1}w'')$.  Il
y a $n$ telles r\'eflexions, autrement dit card$\,\M(w',w'')$ $=\ell
(w'^{-1}w'')$. En particulier, l'axiome $(M)$ est v\'erifi\'e, et le
graphe de l'espace \`a murs $(W,\M)$ s'identifie au graphe de Cayley
de $(W,S)$. Cette identification est $W$--\'equivariante (l'image par
$w$ de $A^+_t$ est $A^{\varepsilon}_{t'}$, o\`u $t'=w^{-1}tw$ et
$\varepsilon =+$ si $w\in A^+_t$, $\varepsilon =-$ sinon).

On peut aussi d\'efinir le syst\`eme de murs $\M$ sur $W$ en
consid\'erant diverses actions de $W$ sur des complexes poly\'edraux.

Si $W$ agit sur un espace $P$ et si $t$ est une r\'eflexion, appelons
{\it mur} de $t$ dans $P$, et notons $M(t,P)$, l'ensemble des points
fixes de $t$ dans $P$. Pour $P$, prenons successivement le graphe de
Cayley de $(W,S)$ (not\'e $\G(W,S)$), la r\'ealisation g\'eom\'etrique
standard de $(W,S)$ (not\'ee $|W|$, voir \cite{Ron}), et enfin sa
r\'ealisation g\'eom\'etrique au sens de Davis--Moussong (not\'ee
$|W|_0$).  Chacun de ces trois complexes est un ``appartement'' au
sens de Davis, voir \cite{Dav} pour les d\'efinitions et
propri\'et\'es concernant ces espaces $W$--homog\`enes; le complexe
$|W| _0$ est introduit dans \cite{Dav}, et muni d'une m\'etrique
CAT$(0)$ dans \cite{Mou}.

Notons que $|W|$ est un complexe simplicial de dimension card$\,S-1$
sur lequel $W$ agit, de mani\`ere simplement transitive sur les
simplexes de dimension maximale.  On identifie les \'el\'ements de $W$
aux centres de ces simplexes maximaux.

(Rappelons bri\`evement la construction de $|W| _0$. Soit $\Delta_S$ le
simplexe standard d'en\-sem\-ble de sommets $S$, dont les faces
s'identifient aux parties de $S$. Si $T$ est une partie de $S$, on
note $W_T$ le {\it sous-groupe sp\'ecial} de $W$ engendr\'e par $T$.
Soit $N=N(W,S)$ le sous-complexe simplicial de $\Delta_S$, appel\'e 
{\it nerf fini} de $(W,S)$, dont les simplexes sont les parties $T$ de
$S$ telles que $W_T$ soit fini. En particulier, $N$ contient tous les
sommets de $\Delta_S$.  Soit $C(W,S)=x_0\ast N'$ le c\^one simplicial
(de sommet $x_0$) sur la subdivision barycentrique $N'$ de $N$. Pour
tout sommet $s$ de $N$, on note $F_s$ l'\'etoile de $s$ dans $N'$,
naturellement contenu dans $C(W,S)$. On consid\`ere alors le quotient
$$W\times C(W,S)/\sim$$ o\`u $\sim$ est la relation d'\'equivalence
engendr\'ee par $(w,x)\sim(w',x')$ s'il existe $s\in S$ tel que
$w'=ws$ et $x'=x\in F_s$. On montre (voir \cite{Mou}) que ce quotient
admet une structure de subdivision barycentrique d'un complexe
poly\'edral euclidien CAT$(0)$ $|W|_0$, d'ensemble de sommets l'image
de $W\times\{\ast\}$, que l'on identifie avec $W$.)

Pour chacune des trois actions consid\'er\'ees,
\begin{itemize}
\item 
  le mur $M$ d'une r\'eflexion de $W$ s\'epare $P$ en deux
  composantes connexes, appel\'ees {\it demi-espaces} de $P$ d\'efinis
  par $M$;
\item 
  dans $P$, il y a un plongement $W$--\'equivariant de $\G(W,S)$,
  \'etendant celui de $W$ (c'est le $1$--squelette de $|W|_0$ par
  construction, et le graphe dual de $|W|$);
\item 
  si $t$ est une r\'eflexion, son mur dans $P$ \'evite $W$, et
  deux \'el\'ements de $W$ sont dans une m\^eme composante connexe de
  $P-M(t,P)$ si et seulement s'ils le sont dans
  $\G(W,S)-M(t,\G(W,S))$.
\end{itemize}

C'est pourquoi, pour chaque r\'eflexion $t$, les intersections de $W$
avec les deux demi-espaces de $P$ d\'efinis par $M(t,P)$ donnent un
mur de $W$ ind\'ependant de $P$. D'autre part, on v\'erifie que, si
$P=\G(W,S)$, l'ensemble de murs ainsi obtenu est $\M(W,S)$.

Puisque le graphe de l'espace \`a murs $(W,\M)$ s'identifie
au graphe de Cayley de $(W,S)$, il est hyperbolique (au sens de
Gromov) si et seulement si $W$ est un groupe hyperbolique.  Nous
v\'erifierons dans la section suivante que la condition (H) est
satisfaite.  Pour information, par un th\'eor\`eme de G.~Moussong
\cite{Mou}, les conditions suivantes sont \'equivalentes: 
\begin{enumerate}
\item 
  $W$ est un groupe hyperbolique;
\item 
  $W$ ne contient pas de sous-groupe isomorphe \`a $\ZZ\times\ZZ$;
\item 
  il n'existe pas de partie $T$ de $S$ telle que $(W_T,T)$ soit un
  syst\`eme de Coxeter affine de rang au moins $3$, ni de paires de
  parties $T_1,T_2$ de $S$, disjointes, avec $W_{T_1},W_{T_2}$
  commutants et infinis.
\end{enumerate}

\noindent{\bf Cas particuliers}\qua 
(Complexes de Benakli--Haglund, voir \cite{Ben,Hag}) Soit $k$ un entier
pair au moins $4$, et $L$ un graphe fini (sans boucle ni ar\^ete
double), de maille (i.e.~la plus petite longueur d'un cycle) au moins
$5$ si $k=4$, et $4$ si $k=6$.  Soit $(W(k,L), S(k,L))$ le syst\`eme
de Coxeter de matrice de Coxeter la matrice d'adjacence du graphe $L$,
avec les $1$ et $0$ remplac\'es respectivement par $\frac{k}{2}$ et
$\infty$.  Il v\'erifie clairement la condition (3) ci-dessus.

Nous noterons $A(k,L)$ la r\'ealisation g\'eom\'etrique au sens de
Davis--Moussong de ce syst\`eme de Coxeter. Alors (voir \cite{Hag})
$A(k,L)$ est un complexe polygonal CAT$(-1)$, dont les polygones sont
des $k$--gones hyperboliques, le link de chaque sommet \'etant isomorphe
\`a $L$.

Si $p$ est un entier pair et $L_{q,q}$ est le graphe biparti complet
sur $q+q$ sommets, alors l'immeuble de Bourdon $I_{p,q}$ est
isomorphe, en tant que complexe polygonal, \`a $A(p,L_{q,q})$.

\begin{figure}[htbp]
\cl{\relabelbox\small
\epsfysize 5cm\epsfbox{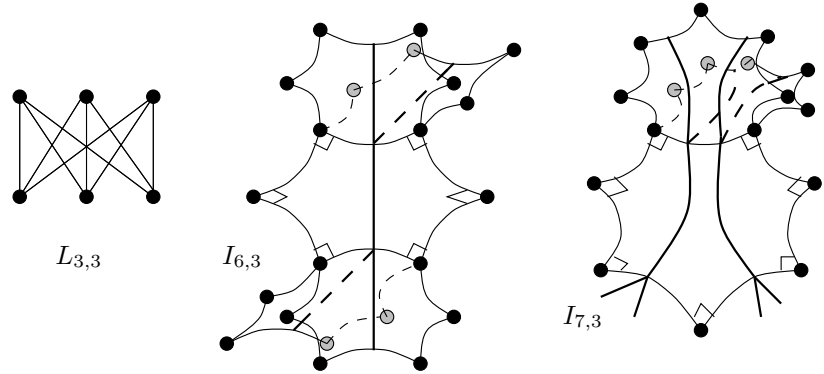}
\relabel {L}{$L_{3,3}$}
\relabel {I}{$I_{7,3}$}
\relabel {II}{$I_{6,3}$}
\endrelabelbox}
\caption{L'immeuble de Bourdon: son link, cas $p$
  pair, cas $p$ impair}
\end{figure}

\medskip
\noindent{\bf Un autre exemple d'espace \`a murs}\qua 
Par contre, si $p=2m+1$ est impair et $q\geq 5$, l'immeuble de Bourdon
$I_{p,q}$ n'est isomorphe ni \`a un complexe polygonal $A(k,L)$ ni \`a
un complexe cubique (sauf \`a passer \`a une subdivision). Supposons
$p\geq 7$. Pour chaque c\^ot\'e fix\'e $A$ du $p$--gone r\'egulier \`a
angles droits $P$, num\'erotons cycliquement $A=A_1,A_2,\cdots,A_p$
les c\^ot\'es de $P$. Consid\'erons les deux segments de
perpendiculaire commune aux paires de c\^ot\'es respectivement $A,
A_{p-1/2}$ et $A, A_{p+1/2}$. Notons $\alpha_1, \alpha_2$ ces
segments.

Nous appellerons {\it mur} de $I_{p,q}$ toute partie $M$ de $I_{p,q}$
ainsi obtenue. Pour tout $i=1,2$ et pour toute identification
isom\'etrique d'un polygone de $I_{p,q}$ avec $P$, on consid\`ere la
r\'eunion $M$ de toutes les g\'eod\'esiques de $I_{p,q}$ passant par
le segment $\alpha_i$.  Nous notons $X_{{p,q}}$ l'ensemble des sommets
de $I_{p,q}$, et {\it mur} de $X_{{p,q}}$ la partition de $X_{{p,q}}$
obtenue en prenant l'intersection de $X_{{p,q}}$ avec les deux
composantes connexes du compl\'ementaire d'un mur de $I_{p,q}$.
(Comme $I_{p,q}$ est simplement connexe, et qu'un mur s\'epare
localement en deux composantes connexes, il s\'epare globalement en
deux composantes connexes.)

Il est facile de montrer que l'espace \`a murs
$(X_{{p,q}},\M_{{p,q}})$ ainsi d\'efini v\'erifie l'axiome (M). Le
graphe associ\'e $\G$ s'identifie avec le $1$--squelette de $I_{p,q}$,
mais les deux sommets de chaque ar\^ete de $I_{p,q}$ sont s\'epar\'es
par exactement deux murs. Ce syst\`eme de demi-espace est diff\'erent
de celui obtenu par subdivision en complexe cubique.  Comme $I_{p,q}$
est CAT$(-1)$, son $1$--squelette est un espace m\'etrique
hyperbolique, de m\^eme bord que $I_{p,q}$. La condition (H) est
facile \`a v\'erifier.

Le groupe Aut$\;I_{p,q}$ des automorphismes de complexe polygonal de
l'im\-meu\-ble de Bourdon $I_{p,q}$ s'identifie naturellement \`a
Aut$\,(X_{{p,q}},\M_{{p,q}})$. En effet, tout automorphisme de
$I_{p,q}$ est une isom\'etrie pour la distance de $I_{p,q}$, et donc
envoie tout segment de perpendiculaire commune entre deux ar\^etes \`a
distance cyclique $q-1/2$ ou $q+1/2$ sur le bord d'un polygone de
$I_{p,q}$ sur un tel autre segment.  Donc il pr\'eserve l'ensemble des
sommets $X_{{p,q}}$ de $I_{p,q}$, ainsi que l'ensemble $\M_{{p,q}}$
des demi-espaces, et Aut$\;I_{p,q}$ est contenu dans
Aut$\,(X_{{p,q}},\M_{{p,q}})$.

Comme les seuls cycles de longueur $p$ dans le $1$--squelette de
$I_{p,q}$ sont les bords des polygones, il en d\'ecoule que
Aut$\;I_{p,q}$ est \'egal \`a Aut$\,(X_{{p,q}},\M_{{p,q}})$

\section{Complexes poly\`edraux pairs \`a courbure n\'egative ou nulle}
\label{sect:complexe_polyedral_pair}

 \subsection{Poly\`edres pairs}
\label{sect:polyedre_pair}

\begin{figure}[htbp]
\cl{\relabelbox\small
\epsfysize 10cm\epsfbox{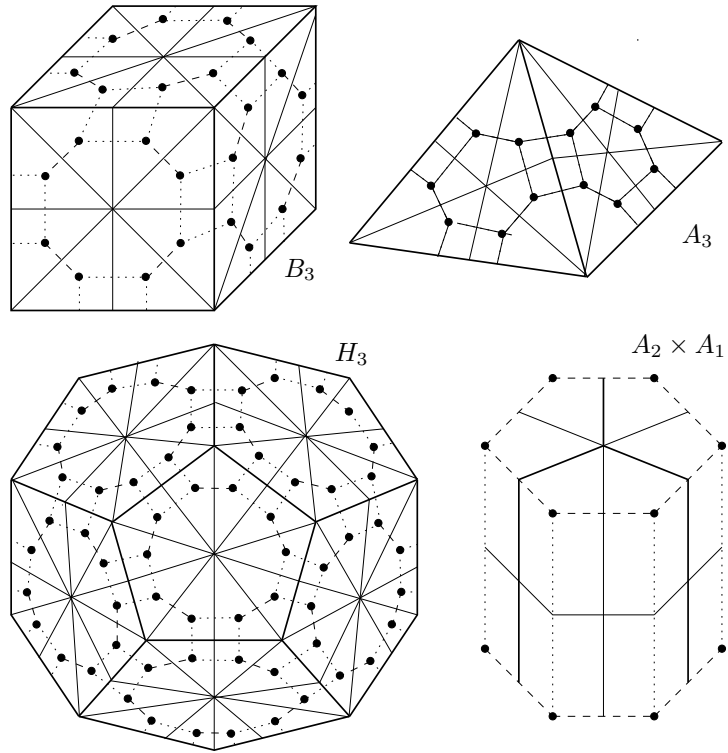}
\relabel {A}{$A_3$}
\relabel {B}{$B_3$}
\relabel {H}{$H_3$}
\relabel {AA}{$A_2\times A_1$}
\endrelabelbox}
\caption{Poly\`edres de Coxeter
  hyperboliques pairs de dimension 3}\label{poly_cox}
\end{figure}

Un poly\`edre (compact convexe) $C$ d'une vari\'et\'e 
riemanienne (compl\`ete, simplement
connexe) \`a courbure constante $\leq 0$ est {\it pair} si
\begin{center}
  \parbox{12cm}{pour toute ar\^ete $a$ de $C$, l'unique r\'eflexion
    $\sigma_{a,C}$ de l'espace ambiant \'echangeant les extr\'emit\'es
    de $a$ pr\'eserve $C$, mais ne fixe aucun sommet de $C$.}
\end{center}
Par exemple, si $C$ est un polygone r\'egulier, il est pair si et
seulement s'il a un nombre pair de c\^ot\'es. Un cube euclidien
r\'egulier de dimension quelconque est pair. Plus g\'en\'eralement, le
produit de deux poly\`edres euclidiens pairs est un poly\`edre
euclidien pair. Voir figure \ref{poly_cox} pour d'autres exemples.  
Nous donnons
ci-dessous une caract\'erisation constructive de tous les poly\`edres
pairs.

Soit $\XX_\kappa$ l'espace \`a courbure constante $\kappa\leq 0$ de
dimension $n$.  Si $\kappa=0$, nous prendrons $\XX_\kappa=\RR^n$. Si
$\kappa<0$, nous utiliserons le mod\`ele de la boule de Poincar\'e
pour l'espace hyperbolique $\XX_\kappa$ \`a courbure constante
$\kappa$. Le groupe des isom\'etries de $\XX_\kappa$ fixant l'origine
s'identifie alors avec $O(n)$. Notons $\phi\co\RR^n\ra\XX_\kappa$
l'exponentielle riemannienne en l'origine (l'identit\'e si $\kappa=0$).
Soit $W$ un groupe fini engendr\'e par des r\'eflexions sur des
hyperplans vectoriels de $\RR^n$.  L'application $\phi$ permet alors
de d\'efinir les notions de {\it chambres, murs ...  dans
$\XX_\kappa$} pour l'action isom\'etrique de $W$ sur $\XX_\kappa$.

\bprop \label{prop:construction_polyedre_pair}
Un poly\`edre (compact convexe) $C$ d'un espace $\XX_\kappa$ \`a
courbure constante $\kappa\leq 0$ est pair si et seulement s'il existe
un point $x$ dans $\XX_\kappa$, un syst\`eme de Coxeter fini $(W,S)$
et une repr\'esentation (injective, envoyant chaque \'el\'ement de $S$
sur une r\'eflexion) $\rho$ de $W$ dans le groupe des isom\'etries de
$\XX_\kappa$ fixant $x$ telle que $C$ est l'enveloppe convexe de
l'orbite par $W$ d'un point $y$ de l'int\'erieur d'une chambre.  De
plus, le $1$--squelette de $C$ est isomorphe au graphe de Cayley de
$(W,S)$.  
\eprop

\dem 
Supposons tout d'abord que $C$ est pair. Notons $W$ le groupe
engendr\'e par les r\'eflexions dans $\XX_\kappa$ par rapport aux
hyperplans m\'ediateurs des ar\^etes de $C$. Puisque $C$ est invariant
par $W$, le groupe $W$ est fini et admet au moins un point fixe, le
{\it centre m\'etrique} $x$ de la cellule $C$ (c'est le centre de
l'unique plus petite boule de $\XX_\kappa$ contenant $C$).  Nous
supposerons que $x$ est l'origine de $\XX_\kappa$.

Fixons $y$ un sommet de $C$, et notons $S$ l'ensemble des r\'eflexions
dans $\XX_\kappa$ par rapport aux hyperplans m\'ediateurs des ar\^etes
de $C$ ayant $y$ pour sommet. Par connexit\'e du $1$--squelette de $C$,
le groupe $W$ est engendr\'e par $S$. Puisque c'est vrai au niveau de
l'espace tangent en $x$ (voir \cite{Bourb} par exemple), le groupe $W$
agit simplement transitivement sur les chambres dans $\XX_\kappa$ (qui
sont les composantes connexes du compl\'ementaire des hyperplans
m\'ediateurs des ar\^etes). Tout sommet de $C$ est contenu dans une
chambre, et la chambre contenant $y$ ne contient pas d'autre sommet de
$C$. Donc le groupe $W$ agit simplement transitivement sur les sommets
de $C$. Le sommet $y$ de $C$ est joint par une ar\^ete pr\'ecis\'ement
aux sommets $sy$ avec $s$ dans $S$. Par d\'efinition du graphe de
Cayley, le $1$--squelette de $C$ s'identifie donc au graphe de Cayley
de $(W,S)$. Comme $C$ est l'enveloppe convexe de ses sommets, $C$ est
bien l'enveloppe convexe de l'orbite de $y$ par $W$.

\smallskip 
R\'eciproquement, soit $C$ l'enveloppe convexe de l'orbite par $W$
d'un point $y$ de l'in\-t\'e\-rieur d'une chambre pour une
repr\'esentation comme dans l'\'enonc\'e d'un syst\`eme de Coxeter
fini $(W,S)$.  Montrons que $C$ est pair. Puisque toutes les images de
$y$ par $W$ sont \`a la m\^eme distance de $x$, par convexit\'e
stricte des sph\`eres, les sommets de $C$ sont exactement les images
de $y$ par $W$.  Le m\^eme argument de convexit\'e stricte montre que
le point $y$ est strictement au-dessus de l'hyperplan affine passant
par les $sy$ pour $s$ dans $S$. Donc les segments de droites entre $y$
et les $sy$ sont des ar\^etes de $C$.  
\endproof

\bprop\label{prop:pair_implique_simple}
Soit $C$ un poly\`edre pair d'un espace \`a courbure constante
n\'egative ou nulle. Alors $C$ est {\it simple}, i.e.~les links de ses
sommets sont des simplexes (sph\'eriques). Si $C$ est euclidien, alors
les longueurs des ar\^etes des links de faces de $C$ sont dans
$[\frac{\pi}{2},\pi]$ (et en particulier ses angles di\`edres sont
obtus).  
\eprop

\dem
Comme le type combinatoire des poly\`edres pairs ne d\'epend pas de la
courbure, nous pouvons supposer $C$ euclidien. Si la dimension $n$ de
$C$ est \'egale \`a celle de l'espace ambiant (ce que nous pouvons
toujours supposer), le groupe fini engendr\'e par des r\'eflexions $W$
construit ci-dessus est {\it essentiel} (i.e.~il ne fixe aucun vecteur
tangent au centre m\'etrique de $C$ non nul).  Si $v$ est un sommet de
$C$, alors les sommets du link de $v$ sont en bijection avec les murs
de la chambre contenant $v$. Or (voir \cite[Ch.~V, section 3,
Prop.~7]{Bourb}) les chambres sont des c\^ones simpliciaux. Donc le
link de $v$ (qui est de dimension $n-1$) a exactement $n-1$ sommets,
et est donc un simplexe.

Si $a,b$ sont deux ar\^etes de $C$, le plan $P$ qui les contient
rencontre perpendi\-culairement les hyperplans m\'ediateurs de $a,b$ en
deux droites $\alpha,\beta$. Les ar\^etes $a,b$ et les droites
$\alpha,\beta$ d\'efinissent un quadrilat\`ere dont deux angles sont
droits et l'un des deux autres est l'angle di\`edre entre les
hyperplans m\'ediateurs de $a,b$.  L'angle di\`edre entre deux
murs d'une m\^eme chambre est dans $[0,\frac{\pi}{2}]$, donc l'angle 
entre deux ar\^etes de $C$ est dans $[\frac{\pi}{2},\pi]$.  La
longueur de toute ar\^ete du link de tout sommet $s$ de $C$ est donc
dans $[\frac{\pi}{2},\pi]$. Par les formules de trigonom\'etrie
sph\'erique, il en d\'ecoule que l'angle en un sommet d'une $2$--face
du link de $s$ est au moins $\frac{\pi}{2}$, donc que la longueur des
ar\^etes des links de face de dimension $1$ est au moins
$\frac{\pi}{2}$. Le r\'esultat en d\'ecoule par r\'ecurrence sur la
dimension de la face.  
\endproof

Nous donnons ci-dessous la liste compl\`ete des poly\`edres
hyperboliques pairs qui sont des poly\`edres de Coxeter, en dimension
2 et 3.  Dans le tableau suivant, $m$ est un entier, avec $m=5$ ou
$m\geq 7$. \`A tout poly\`edre pair $P$ de dimension $n$, et \`a tout
sommet $x_0$ de celui-ci, est associ\'e par la proposition
\ref{prop:construction_polyedre_pair} un syst\`eme de Coxeter fini
$(W,S)$ de rang $n$, dont nous donnons le type et le diagramme de
Coxeter. Les ar\^etes de $P$ issues de $x_0$ sont en bijection avec
les \'el\'ements de $S$. Si $P$ est de dimension $3$, nous donnons les
angles di\`edres $(\alpha_a,\alpha_b,\alpha_c)$ des ar\^etes issues de
$x_0$ correspondant aux \'el\'ements de $S=\{a,b,c\}$. Par la formule
de Gauss--Bonnet, un polygone hyperbolique pair est d\'etermin\'e \`a
isom\'etrie pr\`es par $(m,\alpha,\ell)$ dans
$\NN\setminus\{0,1\}\times]0,\frac{(p-1)\pi}{p}[\times]0+\infty[$,
avec $2m$ son nombre de c\^ot\'es, $\alpha$ l'angle en chacun de ses
sommets, et $\ell$ la longueur d'un de ses c\^ot\'es (et donc des
c\^ot\'es \`a distance paire de celui-ci).

\newcommand{\vsp}[1]{\rule{0pt}{#1}}
$$\begin{array}{|c|c|}
\hline
{\rm Rang~2} & {\rm Rang~3} \\ \hline \vsp{4.5mm}
\begin{array}{lr}
(W,S) & \mbox{}\;\;\;\;\;\;\;\;\;\alpha
\end{array} &
\begin{array}{lr}
(W,S) & \mbox{}\;\;\;\;\;\;\;\;\;(\alpha_a,\alpha_b,\alpha_c)
\end{array}\\[1mm]
\hline
    \begin{array}{cll}
    A_1\times A_1 & \bullet \;\;\;\;\; \bullet & 
                    \frac{\pi}{n}, n\geq 3\\
    A_2  & \bullet \!\!\!\stackrel{\mbox{}}{-\!\!\!-\!\!\!-}\!\! \!\bullet 
& 
                     \frac{\pi}{n}, n\geq 2\\
    B_2  & \bullet \!\!\!\stackrel{4}{-\!\!\!-\!\!\!-} \!\!\!\bullet & 
                    \frac{\pi}{n}, n\geq 2\\
    G_2  & \bullet \!\!\!\stackrel{6}{-\!\!\!-\!\!\!-}\!\!\! \bullet & 
                    \frac{\pi}{n}, n\geq 2\\
    I_2(m)   & \bullet \!\!\!\stackrel{m}{-\!\!\!-\!\!\!-} \!\!\!\bullet & 
                    \frac{\pi}{n}, n\geq 2\\
    \end{array} &
    \begin{array}{cll}\vsp{7mm}
    A_1\times W & \stackrel{a}{\bullet}\;\;\;\;\;  \stackrel{b}{\bullet} 
\!\!\stackrel{m}{-\!\!\!-\!\!\!-}\!\!
\stackrel{c}{\bullet} & 
(\frac{\pi}{2},\frac{\pi}{3},\frac{\pi}{n}), n=3,4,5     \\
 {\rm o\grave{u}~}W=& A_2,B_2,G_2 & {\rm ou~} I_2(m)\\
    A_3 & \stackrel{a}{\bullet}\!\! \stackrel{\mbox{}}{-\!\!\!-\!\!\!-}
    \!\!\stackrel{b}{\bullet} \!\!\stackrel{\mbox{}}{-\!\!\!-\!\!\!-}
\!\!\stackrel{c}{\bullet} & 
(\frac{\pi}{2},\frac{\pi}{n},\frac{\pi}{2}), n\geq 3     \\
\mbox{}& & (\frac{\pi}{2},\frac{\pi}{n},\frac{\pi}{3}), n=3,4,5 \\
    B_3 & \stackrel{a}{\bullet} \!\!\stackrel{\mbox{}}{-\!\!\!-\!\!\!-}
   \!\! \stackrel{b}{\bullet}\!\! \stackrel{4}{-\!\!\!-\!\!\!-}
\!\!\stackrel{c}{\bullet} & 
(\frac{\pi}{2},\frac{\pi}{n},\frac{\pi}{2}), n\geq 3     \\
\mbox{}& & (\frac{\pi}{2},\frac{\pi}{n},\frac{\pi}{3}), n=3,4,5 \\
    H_3 & \stackrel{a}{\bullet} \!\!\stackrel{\mbox{}}{-\!\!\!-\!\!\!-}
    \!\!\stackrel{b}{\bullet}\!\! \stackrel{5}{-\!\!\!-\!\!\!-}
\!\!\stackrel{c}{\bullet} & 
(\frac{\pi}{2},\frac{\pi}{n},\frac{\pi}{2}), n\geq 3     \\
\mbox{}& & (\frac{\pi}{2},\frac{\pi}{n},\frac{\pi}{3}), n=3,4,5 \\
           \vsp{0.5mm} \end{array}\\ 
\hline
\end{array}
$$
\vglue 2mm

\bprop 
\`A isom\'etrie pr\`es, un poly\`edre hyperbolique pair de
dimension 2 ou 3 qui est un poly\`edre de Coxeter est donn\'e \`a
isom\'etrie pr\`es par le tableau pr\'ec\'edent (avec un param\`etre libre
$\ell\in]0,+\infty[$ en rang $2$).
\eprop

\dem 
Soit $(W,S)$ un syst\`eme de Coxeter de rang 3. Soit $Z$ la
cellulation duale de la subdivision barycentrique $\tau$ de la
cellulation de la sph\`ere $\SS^2$ d\'ecrite ci-dessous:
\begin{itemize}
\item 
  la cellulation de la sph\`ere $\SS^2$ par $4,6,8,12,2m$ bigones
  si $(W,S)$ est de type $A_1\times W$ avec $W$ le groupe de Coxeter
  de rang $2$ de type $A_1\times A_1, A_2,B_2,G_2, I_2(m)$
  respectivement;
\item 
  la cellulation bord du t\'etra\`edre, cube, dod\'eca\`edre si
  $(W,S)$ est de type $A_3,B_3,H_3$ respectivement.
\end{itemize}

Notons que si $P$ est un poly\`edre hyperbolique pair construit \`a
partir de $(W,S)$ comme dans la proposition
\ref{prop:construction_polyedre_pair}, alors son bord est isomorphe
\`a la cellulation $Z$.

Par le th\'eor\`eme d'Andr\'eev \cite{And}, si $\alpha$ est une
application de l'ensemble des ar\^etes de $Z$ dans
$]0,\frac{\pi}{2}]$, alors il existe un poly\`edre hyperbolique
(compact), unique \`a isom\'etrie pr\`es, dont la cellulation du bord
est isomorphe \`a $Z$, avec angle di\`edre $\alpha(z)$ le long d'une
ar\^ete $z$ si et seulement si
\begin{enumerate}
\item 
  la somme des angles le long d'un cycle de longueur $3$ dans
  $\tau$ qui ne borde pas un triangle de $\tau$ est strictement
  inf\'erieure \`a $\pi$,
\item 
  la somme des angles le long d'un cycle de longueur $3$ dans
  $\tau$ qui borde un triangle de $\tau$ est strictement sup\'erieure
  \`a $\pi$,
\item 
  la somme des angles le long d'un cycle de longueur $4$ dans
  $\tau$ qui ne borde pas la r\'eunion de deux triangles de $\tau$ est
  strictement inf\'erieure \`a $2\pi$.
\end{enumerate}
Comme il n'existe pas de cycle de longueur $3$ dans $\tau$ qui ne
borde pas un triangle, et que les seuls triangles sph\'eriques de
Coxeter ont pour angles
$\{\frac{\pi}{2},\frac{\pi}{2},\frac{\pi}{n}\}$, $n\geq 2$ ou
$\{\frac{\pi}{2},\frac{\pi}{3},\frac{\pi}{n}\}, n=3,4,5$, le
r\'esultat en d\'ecoule par examination des divers cas possibles.
L'unicit\'e d\'ecoule de l'unicit\'e dans le th\'eor\`eme d'Andreev,
en remarquant que ces poly\`edres ont une sym\'etrie suppl\'ementaire
(i.e.~qui n'est pas dans $W$), par rapport \`a un hyperplan passant
par des sommets.  
\endproof


D\'efinissons maintenant la notion de parall\'elisme d'ar\^etes.  Soit
$C$ un poly\`edre pair de dimension quelconque.  Si $a$ est une
ar\^ete de $C$, nous noterons $M(a,C)$ l'ensemble des points de $C$
fixes par $\sigma_{a,C}$.  C'est un convexe compact de codimension 1
dans $C$, s\'eparant $C$ en deux composantes connexes. Il ne peut
rencontrer une ar\^ete $b$ de $C$ qu'en son milieu, et
perpendiculairement: dans ce cas $M(a,C)=M(b,C)$. Deux ar\^etes $a,b$
de $C$ sont dites {\it parall\`eles dans} $C$ si $M(a,C)=M(b,C)$.  La
relation de parall\'elisme dans $C$ est une relation d'\'equivalence
sur les ar\^etes de $C$.

\subsection{L'espace \`a murs d'un complexe poly\'edral pair}
\label{sect:espace_amur_associe_pair}

Soit $P$ un complexe poly\'edral, n'ayant qu'un nombre fini de types
d'isom\'etrie de cellules.  Nous dirons que $P$ est un {\it complexe
poly\'edral pair} si toute cellule $C$ de $P$ est paire.  Par exemple,
un arbre, ou plus g\'en\'eralement un complexe cubique (voir
\cite{Gro,Sag,NR}) est un complexe poly\'edral pair.

La r\'eunion des relations de parall\'elisme sur les ar\^etes d'une
m\^eme cellule de $P$ engendre une relation d'\'equivalence sur
l'ensemble de toutes les ar\^etes de $P$, que nous appelerons {\it
parall\'elisme entre ar\^etes} dans $P$. (Voir \cite[section~2.4]{Sag} 
pour le cas des complexes cubiques.) D\'efinissons alors le {\it
mur de $P$ transverse \`a une ar\^ete $a$} comme l'union des
$M(b,C')$, avec $b$ une ar\^ete parall\`ele \`a $a$ contenue dans une
cellule (maximale pour l'inclusion) $C'$ de $P$.

Puisque $P$ n'a qu'un nombre fini de types d'isom\'etrie de cellules,
et les compacts $M(b,C')$ ne contenant aucun sommet de $P$, car $C'$
est pair, il vient:
\begin{itemize}
\item 
  tout mur de $P$ est ferm\'e, (localement compact si le link de toute
  cellule de $P$ de dimension $>0$ est compact) et \'evite l'ensemble
  $X$ des sommets de $P$;
\item
  l'ensemble des murs de $P$ est localement fini.
\end{itemize}

Comme dans le cas des complexes cubiques \cite[Theo.~4.10]{Sag}, le 
premier r\'esultat est le suivant.

\blemm\label{lem:mur_geom_separe} 
Soit $P$ un complexe poly\'edral pair {\rm CAT}$(0)$ et $M$ le mur de $P$ 
transverse \`a une ar\^ete $a$. Alors $M$ est 
convexe dans $P$, et s\'epare $P$ en deux composantes connexes.  
\elemm

\dem 
Soit $V(M)$ l'union des cellules de $P$ contenant une ar\^ete
parall\`ele \`a $a$. Donnons d'abord une description du rev\^etement
universel de $V(M)$.

Soit $\C(a)$ l'ensemble des suites de la forme $(a_0,a_1,\ldots
,a_n,C)$, o\`u les $(a_i)_{0\leq i\leq n}$ sont des ar\^etes de $P$,
avec $a_0=a$, $a_i$ parall\`ele \`a $a_{i+1}$ dans une cellule de $P$,
et $C$ est une cellule de $P$ contenant $a_n$. Si $a_i, a_{i+1}$ et
$a_{i+2}$ sont trois ar\^etes parall\`eles \`a $a$ dans une m\^eme
cellule $C'$, nous dirons qu'il y a entre $(a_0,\ldots
,a_i,a_{i+1},a_{i+2},\ldots, a_n,C)$ et $(a_0,\ldots
,a_i,a_{i+2},\ldots ,a_n,C)$ une {\it homotopie \'el\'ementaire} (\`a
extr\'emit\'es fix\'ees).  Les homotopies \'el\'ementaires engendrent
une relation d'\'equivalence sur $\C(a)$: nous noterons $[a_0,\ldots
,a_n,C]$ la classe d'\'e\-qui\-va\-len\-ce de $(a_0,\ldots ,a_n,C)$
pour cette relation.

Soit $\overline{V}(M)$ le complexe poly\'edral obtenu \`a partir de
l'union disjointe des cellules de la forme $[a_0,\ldots ,a_n,C]\times
C$ en identifiant deux points de la forme $([a_0,\ldots ,a_n,C'] ,x')$
et $([a_0,\ldots ,a_n,C''] ,x'')$ lorsque $x'=x''(\in C'\cap C'')$.
Notons $p\co \overline{V}(M)\ra V(M)$ l'application poly\'edrale
naturelle.  Alors $p$ est surjective et un isomorphisme sur chaque
cellule.  Via $p$, le complexe $\overline{V}(M)$ h\'erite d'une
structure de complexe poly\'edral, n'ayant qu'un nombre fini de types
d'isom\'etrie de cellules (pour laquelle $p$ est une isom\'etrie sur
chaque cellule). Montrons que sur $\overline{M} = p^{-1}(M)$,
l'application $p$ est une isom\'etrie.

D'abord, $\overline{M}$ est localement convexe dans $\overline{V}(M)$.
En effet, $\overline{V}(M)$ poss\`ede une r\'eflexion
$\overline{\sigma} _a$ (obtenue sur chaque cellule $\overline C$ de
$\overline{V}(M)$ image de $[a_0,\ldots ,a_n,C]\times C$ en conjugant
$\sigma _{a_n,C}$ par $p_{\vert{} \overline C}$). L'ensemble des
points fixes de l'isom\'etrie $\overline{\sigma} _a$ est
pr\'ecis\'ement $\overline{M}$. Or la m\'etrique de $\overline{V}(M)$
est localement convexe. Il en r\'esulte que $\overline{M}$ est
localement convexe dans $\overline{V}(M)$.

Ensuite, l'image d'une g\'eod\'esique $\overline{\gamma} $ de
$\overline{M}$ par $p$ est une g\'eod\'esique de $P$ contenue dans
$M$. En effet, on remarque d'abord que $\overline{\gamma} $ est une
g\'eod\'esique locale de $\overline{V}(M)$, puis que $p$ est une
isom\'etrie locale au voisinage de $\overline{M}$.  Donc
$p(\overline{\gamma} )$ est une g\'eod\'esique locale de $P$. Mais
comme $P$ est CAT$(0)$, ceci implique que $p(\overline{\gamma})$ est
une g\'eod\'esique globale de $P$.

Puisque $\overline{M}$ est \'evidemment connexe, $M$ est convexe dans
$P$, et $p$ induit une isom\'etrie de $\overline{M}$ sur $M$.

En fait, $p\co \overline{V}(M)\ra V(M)$ est un hom\'eomorphisme. En
effet, notons d'abord qu'un point $\overline x$ de $\overline{V}(M)$
est dans une cellule minimale $\overline C_{\overline x}$ de
$\overline{V}(M)$ rencontrant $\overline{M}$. Si $\overline{x}'$
d\'esigne la projection orthogonale de $\overline{x}$ sur
$\overline{M} \cap \overline{C}_{\overline{x}}$, alors toute
g\'eod\'esique de $\overline{M}$ issue de $\overline{x}'$ fait avec
$[\overline{x}',\overline{x}]$ un angle au moins \'egal \`a
$\frac{\pi}{2}$. Maintenant, si deux points $\overline{x}$ et
$\overline y$ de $\overline{V}(M) - \overline{M}$ sont identifi\'es
par $p$, il appara\^{\i}t dans $P$ un triangle de sommets $p(\overline
x)=p(\overline y)$, $p(\overline{x}')$ et $p(\overline{y}')$, avec des
angles \`a la base sup\'erieurs ou \'egaux \`a $\frac{\pi}{2}$. Comme
$P$ est CAT$(0)$, cela n'est possible que si
$p(\overline{x}')=p(\overline{y}')$. Donc
$\overline{}x'=\overline{y'}$, et $\overline{C}_{\overline x} =
\overline{C}_{\overline y}$. Or $p$ est un plongement sur chaque
cellule: donc $\overline x = \overline y$.

Apr\`es avoir v\'erifi\'e que $\overline{M}$ s\'epare
$\overline{V}(M)$ en deux composantes connexes, on en d\'eduit que $M$
s\'epare $P$ en deux composantes connexes (parce qu'il s\'epare son
voisinage $V(M)$, et que $P$ est simplement connexe).  \endproof

Le r\'esultat suivant d\'ecoule aussi de la preuve du lemme
pr\'ec\'edent.

\blemm \label{lem:prolongement_dans_mur}
Pour toute cellule $C$ de $P$ maximale pour l'inclusion, le mur de
$P$ transverse \`a une ar\^ete $a$ de $C$ est la r\'eunion de tous les
segments g\'eod\'esiques rencontrant $M(a,C)$ en un intervalle
d'int\'erieur non vide.  
\hfill$\Box$ \elemm

Soient $X=X_P$ l'ensemble des sommets de $P$ et $M$ un mur de $P$
transverse \`a une ar\^ete; notons $P^+(M)$ et $P^-(M)$ les deux
composantes connexes de $P-M$. Comme $X\cap M = \emptyset$, la paire
$\{X\cap P^+(M),X\cap P^-(M)\}$ est une partition de $X$. Nous
noterons encore $M$ ce mur de $X$, et $\M=\M_P$ l'ensemble des murs de
$X$ ainsi d\'efini.

\bprop\label{prop:pair_implique_amur}
Soit $P$ un complexe poly\'edral pair {\rm CAT}$(0)$.\nl Alors $(X_P,\M_P)$ est
un espace \`a murs.  
\eprop

\dem 
V\'erifions que $\M$ satisfait l'axiome (M).

Soient $x$ et $y$ deux sommets de $P$, et $\gamma$ la g\'eod\'esique
de $P$ qui les joint. Tout mur de $\M(x,y)$ correspond \`a un mur de
$P$ s\'eparant topologiquement $x$ et $y$, donc coupant $\gamma$.
L'ensemble des murs de $P$ \'etant localement fini, on en d\'eduit que
$\M(x,y)$ est fini.

D'autre part, $\gamma$ part de $x$ par l'int\'erieur d'une (unique)
cellule $C$, elle doit traverser un des murs $M(a,C)$ (avec $a$ une
ar\^ete issue de $x$) avant de retoucher $\partial{C} $: donc
$\M(x,y)$ est non vide.  \endproof

Avant de poursuivre l'\'etude de cet exemple fondamental, il
convient de faire quelques remarques.
 
\medskip
\noindent{\bf Remarque 1}\qua 
Le syst\`eme de murs d'un syst\`eme de Coxeter $(W,S)$ d\'efini dans
la section \ref{sect:espace_a_murs} peut s'obtenir par la pr\'esente
construction, en prenant pour $P$ la r\'ealisation g\'eom\'etrique au
sens de Davis--Moussong $|W| _0$ de $(W,S)$. Par construction m\^eme
(voir \cite{Mou}), une cellule de $|W| _0$ est paire, le groupe
engendr\'e par les r\'eflexions orthogonales le long des ar\^etes de
la cellule \'etant isomorphe \`a un sous-groupe sp\'ecial fini de
$(W,S)$; d'autre part, $|W| _0$ est bien CAT$(0)$ (voir \cite{Mou}).

\medskip
\noindent{\bf Remarque 2}\qua De nombreux complexes poly\`edraux
CAT$(0)$ admettent des subdivisions r\'eguli\`eres cubiques qui
restent CAT$(0)$ lorsqu'on munit les cubes de leurs m\'etriques
euclidiennes standard. Par exemple, si $P$ est un complexe polygonal
CAT$(0)$ sans triangle tel que le link d'un sommet de $P$ ne contient
aucun circuit de longueur $3$, alors la subdivision de chaque $k$--gone
de $P$ en $k$ carr\'es, identifi\'es au carr\'e euclidien unit\'e,
fournit un complexe carr\'e encore CAT$(0)$. Voir par exemple
l'exemple \`a la fin de la section
\ref{sect:exemple_syssteme_coxeter}, o\`u le syst\`eme de murs est
toutefois diff\'erent de celui obtenu par subdivision cubique. Ce
genre de subdivision permet d'appliquer nos r\'esultats de
simplicit\'e \`a des complexes poly\`edraux CAT$(0)$ non
n\'ecessairement pairs (comme l'immeuble de Bourdon avec $p$ impair).

\medskip
\noindent{\bf Remarque 3}\qua On pourrait penser que tout complexe
poly\'edral pair CAT$(0)$ peut \^etre subdivis\'e en cubes, tout en restant
CAT$(0)$, et donc qu'il suffit d'\'etudier les complexes cubiques. 
Mais il n'en est rien, comme le montre l'exemple suivant en dimension 2.
 
Soient $\ell$ et $m$ deux entiers sup\'erieurs ou \'egaux \`a 3.
Consid\'erons un ensemble $S_{\ell ,m}$ de $\ell m$ points, r\'epartis
en $\ell$ colonnes de $m$ points chacune. Relions deux points de
$S_{\ell,m}$ si et seulement s'ils n'appartiennent pas \`a la m\^eme
colonne. Nous noterons $K_{\ell,m}$ le graphe ainsi obtenu (dont le
graphe compl\'ementaire est donc une union disjointe de $\ell$ graphes
complets sur $m$ sommets).  Comme $\ell \geq 3$, ce graphe contient
des circuits de longueur 3.  Fixons d'autre part un entier $k\geq 4$.
Nous pouvons consid\'erer le syst\`eme de Coxeter
$(W_{k,\ell,m},S_{\ell ,m})$ dont le graphe de Coxeter a des ar\^etes
de poids infini entre points d'une m\^eme colonne, et des ar\^etes de
poids $k$ entre points n'appartenant pas \`a la m\^eme colonne. Alors
la r\'ealisation g\'eom\'etrique de Davis--Moussong de $(W_{k,\ell
,m},S_{\ell ,m})$ est (la subdivision barycentrique d')un complexe
polygonal $W_{k,\ell ,m}$--homog\`ene $X$, dont les polygones sont
hyperboliques r\'eguliers \`a $2k$ c\^ot\'es, d'angle aux sommets
$\frac{2\pi}{3}$, et tel que le link de chaque sommet est isomorphe
\`a $K_{\ell ,m}$ (voir \cite{Hag}). Donc $X$ est un complexe
poly\'edral pair CAT$(0)$. Une subdivision en carr\'es de $X$ donne
alors des angles aux sommets \'egaux \`a
$\frac{\pi}{2}$, donc des circuits de longueur \'egale
\`a $\frac{3\pi}{2}<2\pi$ dans le link m\'etrique des sommets, ce qui
emp\^eche $X$ d'\^etre CAT$(0)$.

\subsection{Le graphe associ\'e \`a l'espace \`a murs d'un complexe
 poly\'edral pair}
\label{sect:graphe_associe_complexe_polyedral}

Soit \`a nouveau $P$ un complexe poly\'edral pair CAT$(0)$ et 
$\M=\M_P$ son syst\`eme
de murs sur l'ensemble $X=X_P$ de ses sommets. Nous \'etudions
maintenant les g\'eod\'esiques du 1--squelette $\G$ de $P$, pour la
m\'etrique g\'eod\'esique sur $\G$ rendant chaque ar\^ete isom\'etrique au
segment unit\'e (qui n'est pas forc\'ement celle induite par
$P$). Nous allons voir que cette m\'etrique sur $\G$ v\'erifie des
propri\'et\'es analogues \`a la m\'etrique des mots d'un syst\`eme
de Coxeter.

Si $c=(a_0,a_1,\ldots ,a_n)$ est un chemin combinatoire de $\G$ 
empruntant les $n+1$ ar\^etes $a_0,a_1,\ldots ,a_n$, nous noterons 
$M(c)$ la suite $M(a_0),M(a_1),\ldots ,M(a_n)$ des murs travers\'es 
par $c$.

\blemm\label{lem:sans_repet_geodesique}
 Soit $c$ un chemin combinatoire de $\G$ d'extr\'emit\'es $x$ et $y$.
\begin{itemize}
\item[\rm a)]
Un mur $M$ s\'epare $x$ de $y$ si et seulement s'il appara\^{\i}t 
un nombre impair de fois dans la suite $M(c)$.
\item[\rm b)] 
Si la suite $M(c)$ est sans r\'ep\'etition, alors $c$ 
est une g\'eod\'esique de $\G$.
\end{itemize}
\elemm

\dem 
a)\qua D'une part, tout mur s\'eparant $x$ de $y$ est travers\'e par
$c$. D'autre part, si un mur $M$ est travers\'e un nombre pair de fois
par $c$, c'est donc que $x$ et $y$ sont dans la m\^eme composante
connexe de $P-M$.

 b)\qua Il r\'esulte du a) que, pour un tel chemin, l'ensemble des murs
travers\'es par $c$ est $\M(x,y)$, et la longueur de $c$ est le cardinal
de $\M(x,y)$. Si $c'$ est un autre chemin d'extr\'emit\'es $x$ et $y$, sa
longueur est \'egale au nombre de murs qu'il traverse, donc au moins
\'egale au nombre de murs qu'il traverse un nombre impair de fois. Donc
$c'$ est au moins aussi long que $c$.  
\endproof

\blemm
Soit $\O$ l'ouvert des points de $P$ qui ne sont sur aucun mur de $P$.
Alors toute composante connexe de $\O$ contient un et un seul sommet 
de $P$. 
\elemm

\dem
Puisque deux sommets distincts de $P$ sont toujours s\'epar\'es
par un mur, il y a au plus un sommet de $P$ par composante connexe de 
$\O$.
 
Pour la r\'eciproque, il suffit de consid\'erer le cas o\`u $P$ est
r\'eduit \`a une cellule $C$.  Par la proposition
\ref{prop:construction_polyedre_pair}, ceci d\'ecoule du fait qu'un
groupe de Coxeter fini agit simplement transitivement sur ses chambres.
\endproof

Le r\'esultat suivant montre qu'on peut accompagner une g\'eod\'esique
de $P$ par une g\'eod\'esique de son 1--squelette.

\blemm \label{lem:pistage_geodesique} 
Soient $x'$ et $y'$ deux points du poly\`edre $P$ n'appartenant \`a
aucun mur de $P$, et $\gamma$ la g\'eod\'esique qui les joint dans
$P$. Alors il existe un chemin combinatoire $c$ de $\G$ tel que $M(c)$
est sans r\'ep\'etition, contenu dans la r\'eunion $V(\gamma )$ des
cellules de $P$ touchant $\gamma$, et d'extr\'emit\'es $x$ et $y$
d\'efinis par: $x'$ et $x$ (resp.~$y'$ et $y$) sont dans la m\^eme
composante connexe de $\O$.  
\elemm

\dem 
D'abord, par le lemme \ref{lem:prolongement_dans_mur}, la
g\'eod\'esique $\gamma$ rencontre un nombre fini de murs, en des
points distincts $z_1,z_2,\ldots ,z_n$. Pour prouver le lemme, il
suffit de l'\'etablir lorsque $n=1$. En effet, pour $n$ quelconque, on
d\'ecoupe $\gamma$ en $n$ segments g\'eod\'esiques successifs $\gamma
_i$, contenant $z_i$, d'extr\'emit\'es $x'_i$ et $x'_{i+1}$ contenues
dans aucun mur de $P$. On applique le lemme pour $n=1$ \`a chacun de
ces segments, ce qui fournit $n$ chemins combinatoires $c_1,\ldots
,c_n$, les extr\'emit\'es de $c_i$ \'etant $x_i$ et $x_{i+1}$, seuls
sommets de $P$ appartenant \`a la m\^eme composante connexe de $\O$ que
$x'_i$ et $x'_{i+1}$ respectivement. Ainsi, les $c_i$ se raccordent
pour former un chemin $c$ de $x$ \`a $y$.  De plus, la suite des murs
travers\'es par $c_i$ est sans r\'ep\'etition.  En effet, d'apr\`es le
lemme \ref{lem:sans_repet_geodesique} b), l'ensemble des murs
travers\'es est l'ensemble des murs s\'eparant $x_i$ de $x_{i+1}$, ou
encore l'ensemble des murs s\'eparant $x'_i$ de $x'_{i+1}$,
c'est-\`a-dire pr\'ecis\'ement l'ensemble des murs passant par
$z_i$. Un mur \'etant convexe, il ne peut contenir $z_i$ et $z_j$ pour
$i\neq j$ (sinon, il contiendrait tous les points entre $z_i$ et
$z_j$). Ceci ach\`eve de prouver que la suite des murs travers\'es par
$c$ est sans r\'ep\'etition. Enfin, $c\subset V(\gamma _1)\cup \ldots
\cup V(\gamma _n)\subset V(\gamma )$.

Consid\'erons donc une g\'eod\'esique $\gamma $ entre deux points $x'$
et $y'$ n'appartenant \`a aucun mur, de sorte que $\gamma$ quitte $\O$
en un seul point $z$.

Juste avant $z$ (resp.~juste apr\`es $z$), la g\'eod\'esique $\gamma$
est dans l'int\'erieur d'une unique cellule $C_-$ (resp.~$C_+$). Les
points de $\gamma$ avant (resp.~apr\`es) $z$ sont dans une m\^eme
composante connexe de $\O$, celle de $x$ (resp.~de $y$).  Donc $x\in
C_-$ et $y\in C_+$. En revanche, $x'$ n'est pas n\'ecessairement dans
$C_-$ (ni $y'$ dans $C_+$).

 Si $C$ d\'esigne la plus petite cellule contenant $z$, on a $C\subset
C_-$ (resp: $C\subset C_+$), mais pas n\'ecessairement
\'egalit\'e. Cependant, nous allons montrer que $x\in C$ et $y\in C$
(m\^eme lorsque $C$ est une face stricte de $C_-$ ou $C_+$).

Raisonnons par r\'ecurrence sur $dim(C_-)-dim(C)$. Si ce nombre est
nul, il n'y a rien \`a prouver. Sinon $C$ est contenu dans le bord de
$C_-$, et nous pouvons projeter radialement \`a partir du centre
m\'etrique de $C_-$ sur $\partial C_-$ la partie de $\gamma$ contenue
dans $C_-$. Nous obtenons une g\'eod\'esique par morceaux $\gamma _z$
de $\partial C_-$ aboutissant \`a $z$. Mis \`a part $z$, aucun point
de $\gamma _z$ n'est sur un mur de $C_-$, sinon, par convexit\'e des
murs, le point de $\gamma$ correspondant serait sur le m\^eme mur. La
partie de $\gamma _z$ juste avant $z$ (not\'ee $\gamma _z^- $) est une
g\'eod\'esique aboutissant \`a $z$ dans une face stricte de $C_-$: on
peut lui appliquer l'hypoth\`ese de r\'ecurrence, assurant que
l'unique sommet $x_z$ de $P$ contenu dans la composante connexe de
$\O$ contenant $\gamma _z^- $ est un sommet de $C$. D'autre part, un
point de $\gamma \cap C_-$ (diff\'erent de $z$) et sa projection sur
$\gamma _z$ ne sont s\'epar\'es par aucun mur (par convexit\'e, un tel
mur, qui passe par le centre m\'etrique de $C_-$, devrait contenir le
point de $\gamma$). Ce qui prouve que $x=x_z$ et ach\`eve la
r\'ecurrence.
 
Les deux sommets $x$ et $y$ appartenant \`a une m\^eme cellule $C$
(rencontrant $\gamma$ en $z$, et engendr\'ee par ce point), nous
pouvons consid\'erer une g\'eod\'esique $c$ du 1--squelette de $C$
entre $x$ et $y$. Nous avons d\'ej\`a $c\subset V(\gamma )$.  Il reste
\`a prouver que la suite des murs de $P$ travers\'es par $c$ est sans
r\'ep\'etition. Raisonnons par l'absurde: si c'est le cas, il existe
deux ar\^etes $a$ et $b$ de $c$ d\'efinissant un m\^eme mur $M$ de
$C$, et dont les milieux sont les extr\'emit\'es d'une composante
connexe $c_0$ de $c-M$. Alors, en rempla\c cant $c_0$ par
$\sigma_{a,C} (c_0)$, on obtient un chemin du 1--squelette de $C$ de
m\^eme longueur et m\^emes extr\'emit\'es que $c$, mais avec deux
allers-retours dans les ar\^etes $a$ et $b$. Ceci contredit le fait
que $c$ est g\'eod\'esique.  
\endproof

Le corollaire suivant nous permet d'identifier par la suite le
1--squelette de $P$ au graphe de l'espace \`a murs $(X_P,\M_P)$.

\bcoro \label{coro:graphe_egal_unsquelette}
Deux sommets de $P$ sont li\'es par une ar\^ete de $P$ 
si et seulement s'ils sont li\'es dans $\G(X_P,M_P)$.
\ecoro

\dem 
La condition est bien s\^ur n\'ecessaire. R\'eciproquement, soient $x$
et $y$ deux sommets de $P$ \`a distance combinatoire $n>1$. Il s'agit
de montrer que $x$ et $y$ ne sont pas li\'es dans $\G(X,M)$, autrement
dit qu'il existe un sommet $z$ de $P$ entre $x$ et $y$ (au sens des
murs). Consid\'erons la g\'eod\'esique de $P$ entre $x$ et
$y$. Appliquons-lui le lemme \ref{lem:pistage_geodesique}. Nous
trouvons un chemin $c$ de $\G$ entre $x$ et $y$, tel que la suite des
murs travers\'es par $c$ est sans r\'ep\'etition. En particulier,
d'apr\`es le lemme \ref{lem:sans_repet_geodesique} b), le chemin $c$
est g\'eod\'esique.  Comme $n>1$, le chemin $c$ contient un point $z$
diff\'erent de ses extr\'emit\'es, qui d\'ecoupe $c$ en deux
sous-chemins $c_-$ et $c_+$. Si $c$ est constitu\'e des ar\^etes
$a_1,\ldots ,a_n$, avec $z=a_i\cap a_{i+1},i<n$, on obtient, gr\^ace
au lemme \ref{lem:sans_repet_geodesique} b):
$$
\begin{array}{c}
\M(x,y)=\{  M(a_1),\ldots ,M(a_i),M(a_{i+1}),\cdots,M(a_n)\} , \\
\M(x,z)=\{  M(a_1),\ldots ,M(a_i)\}  \\
\M(z,y)=\{  M(a_{i+1}),\ldots ,M(a_n)\}.
\end{array}
$$
Donc $\M(x,y)$ est bien l'union (disjointe) de $\M(x,z)$ 
et $\M(z,y)$: le point $z$ est entre $x$ et $y$ dans $(X,\M)$.
\endproof

Le r\'esultat suivant est analogue \`a celui des complexes de Coxeter
(voir \cite{Ron}) et des complexes cubiques (voir \cite{Sag}).

\bprop \label{prop:geodesique_murs_traverses_sans_repetition}
Un chemin combinatoire du $1$--squelette est une g\'eod\'esique 
si et seulement si la suite des murs qu'il traverse est sans 
r\'ep\'etition.
\eprop

\dem 
Compte tenu du lemme \ref{lem:sans_repet_geodesique} b), il ne reste
que le sens ``seulement si'' \`a d\'emontrer.  Commen\c cons par un
analogue combinatoire de la convexit\'e des murs de $P$.

\blemm\label{lem:remplacement_voisinage_mur} 
Soient $M$ un mur de $P$, $V(M)$ la r\'eunion des cellules touchant
$M$, $x$ et $y$ deux sommets de $V(M)$. Alors il existe une
g\'eod\'esique de $\G$ d'extr\'e\-mi\-t\'es $x$ et $y$ contenue dans
$V(M)$.  
\elemm

\dem
D'abord, d'apr\`es les hypoth\`eses de finitude sur les types 
d'isom\'etrie des cellules de $P$, il existe un $\varepsilon >0$ 
tel que toute cellule de $P$ passant \`a distance inf\'erieure
ou \'egale \`a $\varepsilon$ de $M$ coupe $M$.

Soit alors $C_x$ une cellule de $P$ contenant $x$ et touchant $M$.  Le
centre m\'etrique $\hat{C}_x$ de $C_x$ est dans $M$, mais le segment
de $x$ \`a $\hat{C}_x$ ne touche aucun mur de $P$ entre ses
extr\'emit\'es (sinon $x$ serait dans ce mur). Nous pouvons donc
trouver sur ce segment un point $x'$ distinct de $\hat{C}_x$, mais
$\varepsilon$ proche de celui-ci, donc $\varepsilon$ proche de $M$. Il
faut noter que $x$ et $x'$ sont dans la m\^eme composante connexe de
$\O$. De m\^eme, il existe un point $y'$ n'appartenant \`a aucun mur,
dans la m\^eme composante connexe de $\O$ que $y$, et $\varepsilon$
proche de $M$. Par convexit\'e (de l'espace $P$ et de $M$ dans $P$),
la g\'eod\'esique $\gamma$ de $x'$ \`a $y'$ reste \`a distance
inf\'erieure ou \'egale \`a $\varepsilon$ de $M$. Par d\'efinition de
$\varepsilon$, cela entra\^{\i}ne que $V(\gamma )\subset V(M)$. Donc
la g\'eod\'esique de $\G$ fournie par le lemme 2 entre $x$ et $y$
reste dans $V(M)$.  
\endproof

\smallskip 
Pour montrer la proposition, consid\'erons un chemin $c$ qui traverse
(au moins) deux fois un mur $M$ de $P$, et prouvons que $c$ n'est pas
g\'eod\'esique. Nous pouvons trouver un sous-chemin $c_0$ de $c$ qui
ne traverse pas $M$, mais dont les extr\'emit\'es sont des sommets $x$
et $y$ d'ar\^etes $a$ et $b$ transverses \`a $M$ et contenues dans
$c$.

\begin{figure}[htbp]
\cl{\relabelbox\small
\epsfysize 4cm\epsfbox{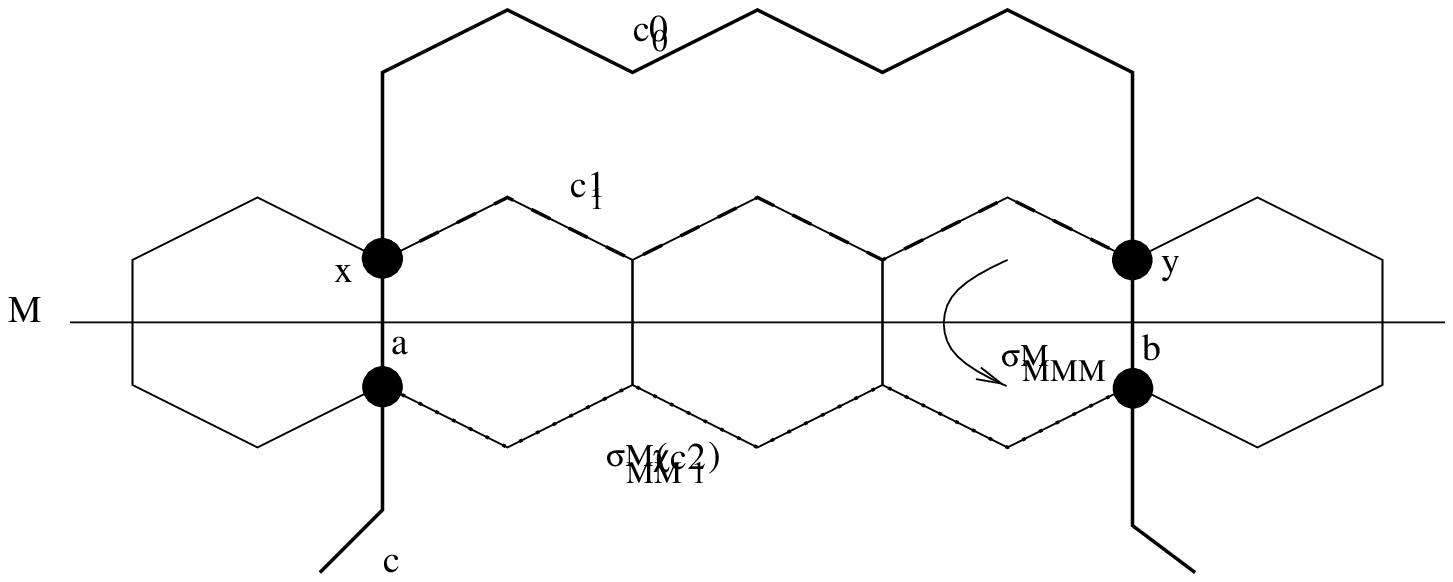}
\relabel {c}{$c$}
\relabel {c1}{$c_1$}
\adjustrelabel <0pt, 1pt> {c0}{$c_0$}
\adjustrelabel <0pt, -1pt> {a}{$a$}
\adjustrelabel <-1pt, -1pt> {y}{$y$}
\adjustrelabel <0pt, -1pt> {b}{$b$}
\adjustrelabel <-1pt, -1pt> {x}{$x$}
\relabel {M}{$M$}
\relabel {sM}{$\sigma_M$}
\adjustrelabel <0pt, -3pt> {sMc}{$\sigma_M(c_1)$}
\endrelabelbox}
\caption{Comment raccourcir les chemins par 
r\'eflexion}
\end{figure}

En appliquant le lemme \ref{lem:remplacement_voisinage_mur}, nous
rempla\c cons $c_0$ par une g\'eod\'esique $c_1$ de $\G$ contenue dans
$V(M)$ et d'extr\'emit\'es $x$ et $y$. Le chemin $c'$ ainsi obtenu a
les m\^emes extr\'emit\'es que $c$, il n'est pas plus long, et il
contient comme sous-chemin $(a,c_1,b)$. Or $V(M)$ poss\`ede une
r\'eflexion $\sigma _M$ par rapport \`a $M$: le chemin $\sigma _M
(c_1)$ a les m\^emes extr\'emit\'es que $(a,c_1,b)$, mais il est plus
court de deux unit\'es. Ceci prouve que ni $c'$, ni a fortiori $c$, ne
sont g\'eod\'esiques.  
\endproof

Compte tenu de la proposition
\ref{prop:geodesique_murs_traverses_sans_repetition}, la preuve du
th\'eor\`eme \ref{theo:Niblo-Reeves} est exactement la m\^eme que
celle du th\'eor\`eme $B$ de \cite{NR}.

\subsection{Hyperbolicit\'e de l'espace \`a murs d'un complexe 
poly\'edral pair}
\label{sect:hyperbolicite_complexe_polyedral}

Soit $P$ un complexe poly\'edral pair CAT$(0)$, dont la m\'etrique
est hyperbolique au sens de Gromov (par exemple, $P$ est CAT$(-1)$).
Comme $P$ n'a qu'un nombre fini de types d'isom\'etrie de cellules, le
diam\`etre des cellules est uniform\'ement major\'e. Donc l'inclusion
du $1$--squelette $\G$ dans $P$ est une quasi-isom\'etrie
(quasi-surjective), et $\G$ est hyperbolique.

Nous allons montrer que la condition (H) est remplie dans
$(X_P,\M_P)$, en \'etablissant son analogue dans $P$.  Comme
d'habitude, nous notons $\overline{P}$ le compactifi\'e de Gromov de
$P$ (donc $\overline{P} = P\cup \partial{P} $), et si $E$ est une
partie de $P$, nous notons $\overline{E}$ son adh\'erence dans
$\overline{P}$. Compte tenu des lemmes \ref{lem:prolongement_dans_mur}
et \ref{lem:mur_geom_separe}, le premier lemme suivant est clair.
 
\blemm \label{lem:non_tangente} 
Soient $M$ un mur de $P$, $x$ un point
de $M$ et $p_x\co \overline{P}\ra lk(x,P)$ la projection qui \`a un rayon
de $P$ d'origine $x$ associe la direction qu'il d\'efinit en partant de
$x$. Alors $\overline{M}$ s\'epare $\overline{P}$ en deux composantes
connexes, images r\'eciproques par $p_x$ des deux composantes connexes
de $lk(x,P)-lk(x,M)$. 
\hfill$\Box$ 
\elemm

\blemm\label{lem:loin_implique_separes} 
Il existe une constante $D>0$
telle que deux points de $P$ \`a distance sup\'e\-rieu\-re ou \'egale
\`a $D$ sont s\'epar\'es par au moins un mur de $P$.  
\elemm

\dem 
Puisque $P$ n'a qu'un nombre fini de types d'isom\'etrie de cellules,
il existe un entier $N$ bornant le nombre de murs susceptibles de
traverser une cellule donn\'ee de $P$. Soient $x$ et $y$ deux points
quelconques de $P$, et consid\'erons deux sommets $x_0$ et $y_0$
contenus dans une m\^eme cellule que $x$ et $y$ respectivement.  Le
nombre des murs s\'eparant $x_0$ de $x$ ou $y_0$ de $y$ est
inf\'erieur \`a $2N$. D'autre part, d'apr\`es l'\'etude de la distance
combinatoire sur $\G$, nous savons que le nombre de murs s\'eparant
$x_0$ de $y_0$ vaut la distance entre $x_0$ et $y_0$ dans $\G$. Cette
distance tend vers l'infini avec la distance dans $P$ entre $x$ et
$y$, par quasi-isom\'etrie entre $P$ et $X$, et puisque le diam\`etre
des cellules est uniform\'ement born\'e.  En particulier, il existe un
nombre $D>0$ tel que, si $d_P(x,y) >D$, alors $x_0$ et $y_0$ sont
s\'epar\'es par au moins $2N+1$ murs de $P$. L'un de ces murs ne
s\'epare ni $x_0$ de $x$, ni $y_0$ de $y$. Donc il s\'epare $x$ de
$y$.  
\endproof

Si $x_0$ est un point base de $P$, $\xi$ un point de $\partial{P}$
et $r_0$ l'unique rayon g\'eod\'esique de $P$ joignant $x_0$ \`a 
$\xi$, nous notons $\M(r_0)$ l'ensemble des murs $M$ de $P$ 
tels que $\overline{M}$ s\'epare $x_0$ de $\xi$.

\blemm\label{lem:murs_sur_rayon_infinis}
 Pour tout rayon g\'eod\'esique $r$ de $P$, l'ensemble $\M(r)$ est infini.
\elemm

\dem 
Consid\'erons la suite de points $(x_k)_{k\geq 0}$ du rayon $r$
d\'efinie par: $x_0$ est l'origine de $r$, et $x_k$ est le point de
$r$ \`a distance $kD$ de $x_0$ --- o\`u $D$ est la constante du lemme
\ref{lem:loin_implique_separes} pr\'ec\'edent.  Il existe donc pour
$k>0$ un mur $M_k$ s\'eparant $x_{k-1}$ de $x_k$. Pour $k<\ell$, on a
n\'ecessairement $M_k\neq M_{\ell}$ (sinon, par convexit\'e, ce mur
contiendrait les points $x_{k+1}$ et $x_{\ell -1}$).

Le mur $M_k$ et le point $\xi$ ne sont pas adh\'erents. En effet, si
$m_k$ d\'esigne le point d'intersection de $M_k$ avec le sous-segment
de $r$ entre $x_{k-1}$ de $x_k$, la projection de $\xi$ dans
$lk(m_k,P)$ correspond \`a la g\'eod\'esique $[ m_k,x_k]$, non
tangente \`a $M$. Le lemme \ref{lem:non_tangente} entra\^{\i}ne bien
que $\xi \not\in \overline{M}$.  La portion de $r$ de $x_0$ \`a
$x_{k-1}$ ne coupe pas $M_k$ (sinon, par convexit\'e, $M_k$
contiendrait $x_{k-1}$). De m\^eme, la portion de $r$ de $x_k$ \`a
l'infini ne coupe pas $M_k$. Mais $M_k$ s\'epare $x_{k-1}$ de
$x_k$. Donc $M_k\in \M(r)$.  
\endproof

 Soient $M$ un mur de $P$ et $\xi$ un point de $\partial P$ non
adh\'erent \`a $M$ dans $\overline{P}$. Nous noterons $V(\xi ,M)$ la
composante connexe de $\overline{P} - \overline{M}$ contenant $\xi$.

\bprop\label{prop:base_de_voisinage}
La famille $\bigl( V(M,\xi )\bigr)_{M\in \M(r_0)}$ est une
base de voisinages de $\xi$ dans $\overline{P}$.
\eprop

\dem 
Remarquons tout d'abord que $V(M,\xi )$ est bien un voisinage de
$\xi$. Pour montrer que la famille est une base de voisinages,
raisonnons par l'absurde. Par d\'efinition de la topologie de
$\G\cup\partial\G$, supposons que les distances $d_P\bigl(x_0,$ $V(M,\xi
)\bigr)$ restent born\'ees pour $M\in \M(r_0)$.

En fait, $x_0$ n'appartient \`a aucun des voisinages $V(M,\xi )$ pour
$M\in \M(r_0)$. Donc $d_P\bigl(x_0,V(M,\xi )\bigr)$ est atteinte sur
le bord de $V(M,\xi )$, c'est-\`a-dire sur $M$.  Nous sommes donc en
train de supposer que tous les murs de $\M(r)$ rencontrent une
certaine boule ferm\'ee $B$ de centre $x_0$ et de rayon $R$.

Si $P$ est suppos\'e localement compact, nous obtenons imm\'ediatement
une contradiction entre la locale finitude de l'ensemble des murs et
le fait que $\M(r)$ est infini.

Donnons un raisonnement g\'en\'eral, o\`u l'on ne suppose plus les
links de sommet de $P$ compacts.  Dans ce cas les boules de $P$ de
rayons trop grands peuvent rencontrer une infinit\'e de
murs. Cependant, par finitude du nombre de types d'isom\'etrie de
cellules de $P$, il existe un $\varepsilon _0 >0$ (qu'on peut choisir
strictement inf\'erieur \`a $R$) et un entier $N_0>0$ tels que toute
$\varepsilon _0$--boule ferm\'ee de $P$ rencontre un nombre de murs
strictement inf\'erieur \`a $N_0$.

Pour un entier $N\geq N_0$, posons $t_N=ND\bigl({{R}\over{\varepsilon
_0}}+1\bigr)$ et $s_N=ND{{R}\over{\varepsilon _0}}$ (le nombre $D$ est
celui qui appara\^{\i}t dans le lemme
\ref{lem:loin_implique_separes}).  Appelons $x_N$ (resp.~$y_N$) le
point du rayon $r$ \`a distance $s_N$ (resp.~$t_N$) de l'origine
$x_0$.  Montrons tout d'abord que toute g\'eod\'esique $\gamma$ de $P$
joignant un point $u$ de la boule $B$ \`a un point $v$ de $r$ entre
$x_N$ et $y_N$ passe par la $\varepsilon _0$--boule ferm\'ee de $P$ de
centre $x_N$.

\begin{figure}[htbp]
\cl{\relabelbox\small
\epsfysize 3cm\epsfbox{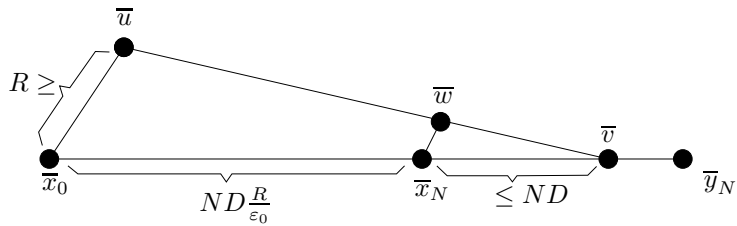}
\relabel {w}{$\overline w$}
\adjustrelabel <0pt, -4pt> {v}{$\overline v$}
\adjustrelabel <0pt, -2pt> {u}{$\overline u$}
\adjustrelabel <0pt, -1pt> {x}{$\overline x_0$}
\relabel {xb}{$\overline x_N$}
\relabel {y}{$\overline y_N$}
\adjustrelabel <-8pt, -2pt> {Rg}{$R\ge$}
\adjustrelabel <-2pt, -1pt> {lND}{$\le ND$}
\adjustrelabel <0pt, -2pt> {ND}{$ND{R\over\varepsilon_0}$}
\endrelabelbox}
\caption{Triangle de comparaison}
\end{figure}

En effet, consid\'erons le triangle g\'eod\'esique de $P$ dont les
sommets sont $x_0$ et les extr\'emit\'es $u,v$ de $\gamma$.  Soient
$\overline{x_0}, \overline{u},\overline{v}$ les sommets correspondants
d'un triangle euclidien de comparaison.  Si $\overline{x_N}\in
[\overline{x_0}, \overline{v}]$ est le point correspondant \`a $x_N$,
alors:
$$d(\overline{x_0}, \overline{u})\leq R, d(\overline{x_0},
\overline{v})\geq ND\frac{R}{\varepsilon_0}, d(\overline{x_N},
\overline{v})\leq ND.$$
Soit $\overline{w}$ le point de
$[\overline{u}, \overline{v}]$ situ\'e sur la parall\`ele au c\^ot\'e
$[\overline{x_0}, \overline{u}]$ passant par $\overline{x_N}$. Alors
par le th\'eor\`eme de Thal\`es, il vient
$$\frac{d(\overline{w},\overline{x_N})}{ND}\leq
\frac{R}{ND\frac{R}{\varepsilon_0}}.$$
Par l'in\'egalit\'e CAT$(0)$,
la distance de $x_N$ au point $w$ de $\gamma$ correspondant \`a
$\overline{w}$ est donc inf\'erieure \`a $\varepsilon_0$.

Pour achever la d\'emonstration de la proposition, d\'ecoupons le
sous-segment de $r$ entre $x_N$ et $y_N$ en $N$ intervalles de
longueur $D$. Par le lemme \ref{lem:loin_implique_separes}, on trouve
$N$ murs deux \`a deux distincts s\'eparant les extr\'emit\'es de ces
intervalles. Ces $N$ murs sont dans $\M(r_0)$ (voir preuve du lemme
\ref{lem:murs_sur_rayon_infinis}).  Ils passent par un point du
sous-segment de $r$ entre $x_N$ et $y_N$, et d'autre part ils coupent
la boule $B$ par hypoth\`ese.  Par convexit\'e des murs et ce qui
pr\'ec\`ede, chacun de ces murs coupe la $\varepsilon _0$--boule
ferm\'ee de $P$ de centre $x_N$.  Ainsi cette boule est coup\'ee par
$N$ murs, avec $N\geq N_0$, en contradiction avec les d\'efinitions de
$\varepsilon _0$ et $N_0$.  
\endproof

L'image r\'eciproque par l'inclusion canonique de $\G$ dans $P$ est
une quasi-isom\'e\-trie, se prolongeant en un hom\'eomorphisme entre les
bords. On obtient ainsi un plongement $\overline{i}\co \overline{\G}\ra
\overline{P}$.  De plus l'image r\'eciproque par $\overline{i}$ d'un
voisinage d'un point de $\partial{P} $ est un voisinage du point
correspondant sur $\partial{\G} $. La proposition pr\'ec\'edente
entra\^{\i}ne donc que $(X,\M)$ v\'erifie l'axiome (H).

Nous r\'esumons les r\'esultats \ref{prop:pair_implique_amur},
\ref{coro:graphe_egal_unsquelette}, \ref{prop:base_de_voisinage} dans
l'\'enonc\'e suivant.

\btheo \label{theo:resume_prop_complexe_pair} 
Soit $P$ un complexe poly\'edral pair {\rm CAT}$(0)$, hyperbolique au sens
de Gromov.  Alors $(X_P,\M_P)$ est un espace \`a murs hyperbolique,
dont le graphe associ\'e est le 1--squelette de $P$.  
\etheo

Un mur $M$ de $P$ est dit {\it propre} si $\partial P\setminus
\partial A$ est non vide pour chacune des composantes connexes $A$ de
$P\setminus M$. Ceci \'equivaut au fait que le mur correspondant de
l'espace \`a murs $(X_P,\M_P)$ est propre.

\blemm \label{lemm:tout_mur_propre} 
Supposons que chaque ar\^ete de $P$ soit contenue dans une droite
g\'eod\'esique. Alors tout mur de $P$ est propre.  
\elemm

\dem 
Soit $M$ un mur transverse \`a une ar\^ete $d$, et $A,B$ les deux
composantes connexes de $P\setminus M$. Soit $D$ une droite
g\'eod\'esique contenant $d$, et $a,b$ l'extr\'emit\'e du rayon
g\'eod\'esique $D\cap A,D\cap B$ respectivement. Alors puisque $P$ est
CAT$(0)$ et que l'angle entre $M$ et $d$ est droit au point
d'intersection, le point $a$ n'appartient pas \`a $\partial B$, ni $b$
\`a $\partial A$.  Donc $M$ est propre.  
\endproof

\section{Groupes d'automorphismes d'un complexe poly\-\'ed\-ral pair}
\label{sect:groupe_auto_poly}

Nous fixons $P$ un complexe poly\'edral pair CAT$(0)$. Nous notons
$(X,\M)=(X_P,M_P)$ son espace \`a murs associ\'e et $\G$ le
$1$--squelette de $P$.

\subsection{Automorphismes de l'espace \`a murs d'un complexe 
  poly\-\'ed\-ral pair}
\label{sect:groupe_auto_polyedre_pair}

Le but de cette section est de montrer que le groupe des
automorphismes de $P$ et celui de $(X,\M)$ co\"{\i}ncident.

Si $f$ est un automorphisme isom\'etrique de $P$, $C$ une cellule de
$P$ et $a$ une ar\^ete de $C$, alors $f(M(a,C))=M(f(a),f(C))$. Aussi,
tout automorphisme isom\'etrique de $P$ agit sur l'ensemble des murs
de $P$.  Plus g\'en\'eralement, un isomorphisme (non n\'ecessairement
isom\'etrique) entre deux cellules paires pr\'eserve le parall\'elisme
entre ar\^etes. En effet, deux ar\^etes $a$ et $b$ d'une cellule paire
$C$ sont parall\`eles si et seulement s'il existe une g\'eod\'esique
combinatoire $\gamma$ du 1--squelette de $C$ joignant une extr\'emit\'e
de $a$ \`a une extr\'emit\'e de $b$, de sorte que $a$ suivie de
$\gamma$, ainsi que $\gamma$ suivie de $b$, soit encore
g\'eod\'esique, mais $(a,\gamma ,b)$ n'est plus g\'eod\'esique.
D'autre part, deux sommets $x$ et $y$ sont du m\^eme c\^ot\'e d'un mur
$M$ \ssi\ une g\'eod\'esique de $x$ \`a $y$ ne contient pas d'ar\^ete
transverse \`a $M$.

Ainsi, le parall\'elisme des ar\^etes est une notion ne faisant appel
qu'\`a la combinatoire de $C$, et m\^eme seulement de son 1--squelette.
Si $f$ est un isomorphisme (poly\'edral) d'une cellule paire $C$ sur
une autre cellule paire $C'$, et si $M$ est un mur de $C$, alors les
ar\^etes de $C'$ images par $f$ des ar\^etes de $C$ transverses \`a
$M$ sont toutes transverses \`a un m\^eme mur de $C'$, qu'on notera
$f(M)$. Et deux sommets $x$ et $y$ de $C$ sont du m\^eme c\^ot\'e de
$M$ si et seulement si $f(x)$ et $f(y)$ sont du m\^eme c\^ot\'e de
$f(M)$.

Les r\'esultats pr\'ec\'edents restent valables pour $P$ tout entier.
Il y a donc un morphisme canonique (d'ailleurs clairement injectif) du
groupe Aut$(P)$ des automorphismes (poly\`edraux) de $P$ dans
Aut$(X,\M)$.

\btheo\label{theo:meme_groupe_automorphisme} Soit $P$ un complexe
poly\'edral pair {\rm CAT}$(0)$. Alors le morphisme de {\rm Aut}$(P)$ dans
{\rm Aut}$(X_P,\M_P)$ ci-dessus est un isomorphisme.  \etheo

\dem Si $\G(X,\M)$ est le graphe associ\'e \`a $(X,\M)$, alors nous
avons d\'efini un morphisme injectif Aut$(X,\M)\ra {\rm
  ~Aut~}\G(X,\M)$.  Comme $\G(X,\M )$ s'identifie avec le 1--squelette
combinatoire $\G$ de $P$, si $\rho\co {\rm ~Aut~}P\ra {\rm ~Aut~}\G$ est
l'applic\-ation de restriction d'un automorphisme de $P$ \`a son
1--squelette, alors le diagramme suivant est commutatif:
$$\begin{array}{ccccc}
  & & {\rm ~Aut~}(X,\M) & & \\
  & \nearrow & & \searrow & \\
  {\rm ~Aut~} P & & \stackrel{\rho}{\longrightarrow} & & {\rm ~Aut~}
  \G
\end{array}
$$
Pour \'etablir que tous ces morphismes injectifs sont des
isomorphismes, il suffit de montrer que $\rho$ est surjective,
i.e.~que l'on peut construire un automorphisme (poly\'edral) de $P$
\`a partir d'un automorphisme de son $1$--squelette $\G$.

\blemm\label{lem:lemme_un} 
Soient $C$ une cellule de $P$ et $a$ une
ar\^ete de $P$ telle que l'inter\-sec\-tion $a\cap C$ est r\'eduite \`a un
sommet $x_0$. Alors le mur transverse \`a $a$ ne coupe pas $C$.
\elemm

\dem 
Supposons, par l'absurde, qu'il existe une cellule $C$, une ar\^ete
$a$ et un sommet $x_0$ tels que $a\cap C=\{ x_0\}$ et $M=M(a)$ coupe
$C$. Soient $y_0$ le sommet de $C$ sym\'etrique de $x_0$ par rapport
\`a $M$, et $p$ le point o\`u la g\'eod\'esique qui joint $x_0$ \`a
$y_0$ (dans $C$) coupe $M$. Alors $p$ est le point de $M\cap C$ le
plus proche de $x_0$.

En fait, pour toute cellule $D$ dont $C$ est une face, $p$ est encore
le point de $M\cap D$ le plus proche de $x_0$: donc $p$ est un minimum
local (strict) pour la fonction qui \`a un point $q$ de $M$ associe sa
distance \`a $x_0$ dans $P$.  Mais il en va de m\^eme pour le point
$p'$, milieu de l'ar\^ete $a$. Or $p'\neq p$, puisque $a\not\subset
C$, ce qui donne deux minimaux locaux sur $M$ \`a la fonction
``distance \`a $x_0$'', en contradiction avec la convexit\'e de cette
fonction et celle de $M$ dans $P$.  
\endproof

\bcoro 
Le $1$--squelette d'une cellule $C$ est convexe dans $\G$, le
1--squelette de $P$.  
\ecoro

\dem 
Soient $x_0,y_0$ deux sommets de $C$, et $\gamma$ un chemin de $\G$
entre $x_0$ et $y_0$, qui sort de $C$.  D'apr\`es le lemme
pr\'ec\'edent, la suite des murs travers\'es par $\gamma$ contient un
mur ne coupant pas $C$.  Or l'ensemble des murs qui s\'epare $x_0,y_0$
est contenu dans l'ensemble des murs coupant $C$. Donc, d'apr\`es la
proposition \ref{prop:geodesique_murs_traverses_sans_repetition},
$\gamma$ ne peut \^etre g\'eod\'esique.  
\endproof

Notons ${\cal E}$ l'ensemble des cellules de $P$ et ${\cal F}$
l'ensemble des sous-graphes convexes de $\G$ isomorphes au graphe de
Cayley d'un syst\`eme de Coxeter fini.

Comme le $1$--squelette d'une cellule paire est le graphe de Cayley
d'un syst\`eme de Coxeter fini (Proposition
\ref{prop:construction_polyedre_pair}), le corollaire ci-dessus montre
que l'application $i\co C\mapsto C\cap \G$ est une application
(injective) de ${\cal E}$ dans ${\cal F}$.  Pour retrouver les
cellules de $P$ \`a partir de son 1--squelette, nous allons montrer que
$i ({\cal E})={\cal F}$.

\blemm \label{lem:lemme_trois} 
Soient $K$ un \'el\'ement de ${\cal F}$
et $a$ une ar\^ete de $P$ telle que l'in\-ter\-sec\-tion $a\cap K$ est
r\'eduite \`a un sommet $x_0$.  Alors le mur $M$ transverse \`a
l'ar\^ete $a$ ne recoupe pas $K$.  
\elemm

\dem
Raisonnons par l'absurde. Soit $b$ une ar\^ete de $K$ transverse \`a
$M$.  Notons $y_0$ l'extr\'emit\'e de $b$ du m\^eme c\^ot\'e de $M$
que $x_0$, puis $x'_0$ et $y'_0$ les images de $x_0$ et $y_0$ par la
r\'eflexion $\sigma _M$ du voisinage $V(M)$ de $M$.  Comme le
1--squelette de $V(M)$ est g\'eod\'esique dans le 1--squelette de $P$
(voir lemme \ref{lem:remplacement_voisinage_mur}), il existe une
g\'eod\'esique $\gamma$ de $\G$ entre $x_0$ et $y_0$ contenue dans
$V(M)$. Mais comme $K$ est convexe dans $\G$, on a $\gamma \subset K$.
$$\begin{array}{c}
\mbox{\relabelbox\small\epsfysize 3cm\epsfbox{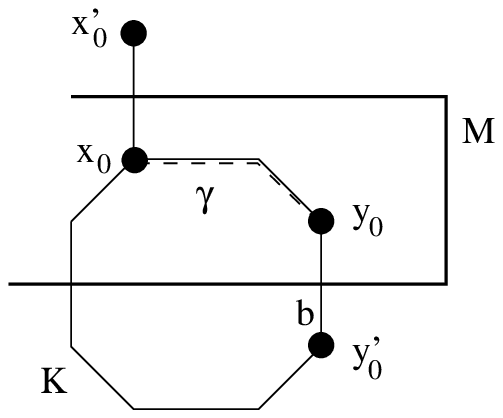}
\relabel {M}{$M$}
\adjustrelabel <-2pt, 0pt> {K}{$K$}
\relabel {g}{$\gamma$}
\adjustrelabel <-1pt, 0pt> {b}{$b$}
\adjustrelabel <-1.5pt, 0pt> {x}{$x_0$}
\adjustrelabel <-2pt, 0pt> {y}{$y_0$}
\adjustrelabel <-1.5pt, 0pt> {x'}{$x'_0$}
\adjustrelabel <-2pt, 0pt> {y'}{$y'_0$}
\endrelabelbox}
\end{array}$$
Puisque $x_0$ et $y_0$ ne sont pas s\'epar\'es par $M$, le chemin
$\gamma$ ne coupe pas le mur $M$. Donc $(\gamma ,b)$ est une
g\'eod\'esique de $\G$ entre $x_0$ et $y'_0$ (voir lemme
\ref{prop:geodesique_murs_traverses_sans_repetition}).  Le chemin
$(a,\sigma _M(\gamma ))$ a les m\^emes extr\'emit\'es et la m\^eme
longueur, mais il passe par $x'_0\not\in K$: ceci contredit la
convexit\'e de $K$ dans $\G$.  
\endproof

\blemm\label{lem:lemme_quatre}
 Soient $K$ un \'el\'ement de ${\cal F}$ et $M$ un mur coupant une
ar\^ete $a$ de $K$. Alors chaque ar\^ete de $K$ touchant $a$ est
contenue dans $V(M)$, et l'ensemble de ces ar\^etes est invariant par
$\sigma _M$.  
\elemm

\dem
Soient $x_0$ et $y_0$ les extr\'emit\'es de $a$, et $b$ une ar\^ete de
$K$ distincte de $a$, contenant $y_0$.  Il s'agit de montrer que $b$
est dans $V(M)$, et que $\sigma _M(b)$ est dans $K$.

Soit $(W,S)$ le syst\`eme de Coxeter de graphe de Cayley $\G(W,S)$
isomorphe \`a $K$. Puisque $W$ est transitif sur les sommets de
$\G(W,S)$, on peut trouver un isomorphisme $\varphi$ de $\G(W,S)$ sur
$K$ envoyant 1 sur $x_0$. Soient $s$ et $w$ les \'el\'ements de $W$
dont l'image par $\varphi$ sont $y_0$ et $z_0$, la deuxi\`eme
extr\'emit\'e de $b$. D'abord, $s\in S$, puisqu'il est li\'e \`a 1
dans $\G(W,S)$ par l'ar\^ete $\varphi ^{-1}(a)$. Ensuite, il existe
$t\neq s$, $t\in S$ tel que $w=st$. Consid\'erons $\G_{s,t}$, le
sous-graphe plein de $\G(W,S)$ dont les sommets sont
$1,s,st,sts,\ldots ,1$. C'est un graphe hom\'eomorphe \`a un cercle,
contenant $2m_{s,t}$ ar\^etes, o\`u $m_{s,t}$ d\'esigne l'ordre du
produit $st$ dans $W$. Ce sous-graphe est une maille de $\G (W,S)$, au
sens suivant: une {\it maille} est un circuit de longeur $2m$
totalement g\'eod\'esique dans $\G(W,S)$, tel que si deux de ses
sommets sont \`a distance strictement inf\'erieure \`a $m$, il y a une
unique g\'eod\'esique de $\G(W,S)$ les joignant (alors
n\'ecessairement contenue dans le circuit).

L'image de $\G_{s,t}$ dans $K$ est une maille $K_{s,t}$ de $K$; par
convexit\'e de $K$ dans $\G$, c'est aussi une maille de $\G$.
L'ar\^ete $a'$ de $K_{s,t}$ la plus \'eloign\'ee de $a$ est
caract\'eris\'ee par l'existence d'un sous-segment $c$ de $K_{s,t}$,
tel que $c$ joint $y_0$ \`a une extr\'emit\'e $y'_0$ de $a'$, $(a,c)$
et $(c,a')$ sont g\'eod\'esiques, mais $(a,c,a')$ ne l'est pas.

Des trois derni\`eres propri\'et\'es et de la proposition
\ref{prop:geodesique_murs_traverses_sans_repetition}, il r\'esulte que
$M(a)=M(a')$.

D'apr\`es le lemme \ref{lem:remplacement_voisinage_mur}, il existe une
g\'eod\'esique de $y_0$ \`a $y'_0$ contenue dans $V(M)$. Mais comme
$K_{s,t}$ est une maille, cette g\'eod\'esique est $c$.  Alors $\sigma
_M(c)$ est une g\'eod\'esique entre deux points de $K_{s,t}$ \`a
distance strictement inf\'erieure \`a $m_{s,t}$, donc $\sigma
_M(c)\subset K_{s,t}$.  
\endproof

\blemm\label{lem:lemme_cinq}
Soient $K$ un \'el\'ement de ${\cal F}$ et $C$ une cellule de $P$ dont
le 1--squelette contient un sommet $x_0$ de $K$ tel que
$St(x_0,K)=St(x_0,i(C))$. Alors $K=i(C)$.  
\elemm

\dem
On peut supposer la dimension de $C$ au moins \'egale \`a deux, sinon
il n'y a rien \`a montrer.  Par connexit\'e de $C$, il suffit de
montrer que $K$ est un ouvert de $i(C)$. Par connexit\'e de $K$, il
suffit de montrer que si $x_0$ est un sommet de $K$ tel que
$St(x_0,K)=St(x_0,i(C))$, alors pour tout voisin $y_0$ de $x_0$ dans
$K$, on a encore $St(y_0,K)=St(y_0,i(C))$.

Soit $a$ l'ar\^ete de $K$ d'origine $x_0$ et d'extr\'emit\'e $y_0$, et
$M$ le mur transverse \`a $a$.  Comme $\sigma _M$ pr\'eserve $i(C)$ et
l'ensemble des ar\^etes de $K$ touchant $a$ (d'apr\`es le lemme
\ref{lem:lemme_quatre}), on a:
$$ St(y_0,K)=St(\sigma _M(x_0),K)=\sigma _M (St(x_0,K))=$$ $$\sigma _M
(St(x_0,i(C)))=St(y_0,i(C)).$$\vglue-.9cm\endproof

\bprop
L'application $i\co {\cal E}\ra{\cal F}$ est surjective.
\eprop

\dem
On raisonne par r\'ecurrence sur le rang du syst\`eme de Coxeter dont $K$
est le graphe de Cayley (cela correspond
au degr\'e du graphe r\'egulier $K$). Il n'y a rien \`a dire en rang 1.

Soit $K\in {\cal F}$ de rang $r+1$ sup\'erieur ou \'egal \`a $2$.
Consid\'erons un sommet $x_0$ de $K$, et soient $a_0,a_1,\ldots ,a_r$
les ar\^etes issues de $x_0$; nous noterons $M_i$ le mur transverse
\`a l'ar\^ete $a_i$. Alors il existe un syst\`eme de Coxeter fini
$(W,S=\{ s_0, s_1,\ldots ,s_r\} )$ et un isomorphisme de son graphe de
Cayley $\G(W,S)$ sur $K$ envoyant 1 sur $x_0$ et l'ar\^ete issue de 1
pr\'eserv\'ee par $s_i$ sur $a_i$.  Consid\'erons maintenant $V$, le
sous-groupe sp\'ecial de $(W,S)$ engendr\'e par $T=\{ s_1,\ldots
,s_r\}$. Il y a une unique copie de son graphe de Cayley contenue dans
le graphe de Cayley de $(W,S)$ et passant par 1; \`a ce sous-graphe
correspond un sous-graphe $L$ de $K$.

Un r\'esultat classique sur les sous-groupes sp\'eciaux
(cf.~\cite{Bourb}) entra\^{\i}ne que $L$ est convexe dans $\G$. On
peut donc appliquer l'hypoth\`ese de r\'ecurrence \`a $L$, et trouver
une cellule $C$ de $P$ dont $L$ est le 1--squelette.

Comme $(W,S)$ est fini, il poss\`ede un unique \'el\'ement $w_0$ de
longueur maximale: soit $x'_0$ le sommet correspondant de $K$.

\medskip
\noindent{\bf Fait 1}\qua 
Notons d'abord que les ar\^etes de $K$ issues de $x'_0$ sont travers\'ees
par les murs $M_i$, qui de plus s\'eparent $x'_0$ de $x_0$.
$$\begin{array}{c}
\mbox{\relabelbox\small\epsfysize 6.5cm\epsfbox{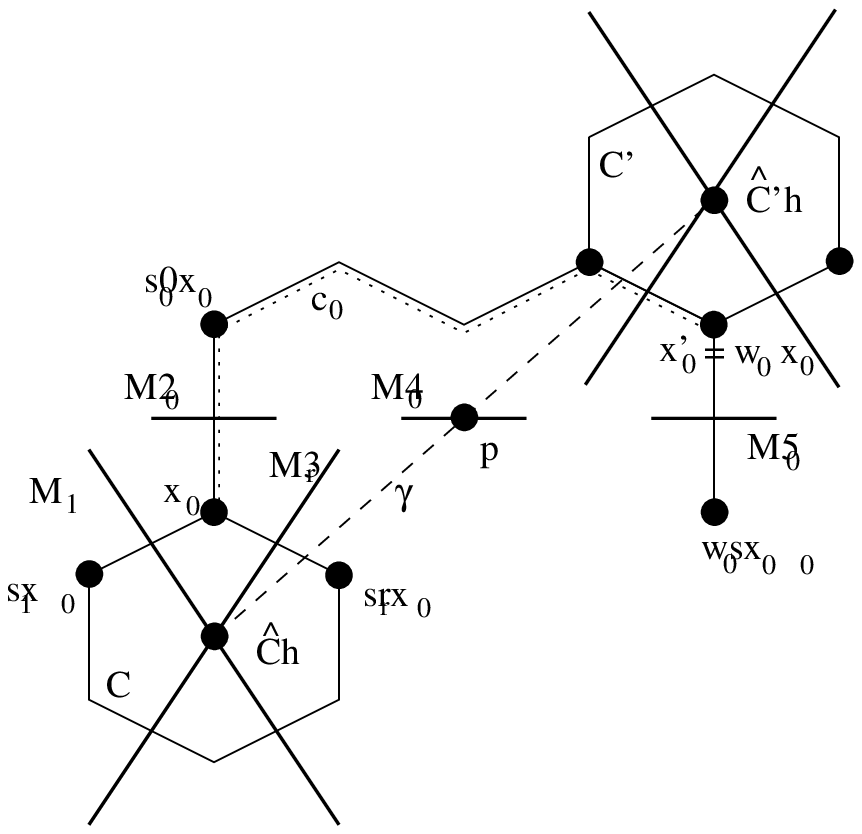}
\relabel {M}{$M_1$}
\relabel {M2}{$M_0$}
\adjustrelabel <-2pt, 0pt> {M3}{$M_r$}
\relabel {M4}{$M_0$}
\relabel {M5}{$M_0$}
\relabel {C}{$C$}
\relabel {Ch}{$\hat C$}
\relabel {C'}{$C'$}
\relabel {C'h}{$\hat C'$}
\relabel {c}{$c_0$}
\relabel {p}{$p$}
\relabel {g}{$\gamma$}
\adjustrelabel <-1pt, 0pt> {x}{$x_0$}
\adjustrelabel <-3pt, 0pt> {x'}{$x'_0\!=\!w_0x_0$}
\adjustrelabel <-1pt, -1pt> {sx}{$s_1x_0$}
\adjustrelabel <-2pt, 0pt> {srx}{$s_rx_0$}
\relabel {s0x}{$s_0x_0$}
\adjustrelabel <-5pt, 0pt> {wsx}{$w_0s_0x_0$}
\endrelabelbox}
\end{array}$$

 \dem 
Pour tout $s_i$ de $S$, l'\'el\'ement $s_iw_0$ doit \^etre li\'e
\`a $w_0$ dans $\G(W,S)$, ce qui signifie qu'il existe un $s_j\in S$
tel que $s_iw_0=w_0s_j$. Soit $c_i$ une g\'eod\'esique de $\G(W,S)$ de
1 \`a $w_0$ commen\c cant par l'ar\^ete de 1 \`a $s_i$: la suite des
murs de $\G(W,S)$ travers\'es par $c_i$ est sans r\'ep\'etition. Alors
le chemin $\gamma _i$ form\'e de $c_i$ suivi de l'ar\^ete de $w_0$ \`a
$s_i.w_0=w_0.s_j$ n'est pas g\'eod\'esique, car le premier mur qu'elle
traverse est $M(s_i)=M(w_0.s_j.w_0^{-1})$, donc \'egal au dernier. On
en d\'eduit que le chemin $c'_i$ tel que l'ar\^ete de 1 \`a $s_i$
suivie de $c'_i$ \'egale $\gamma _i$ est g\'eod\'esique.

En prenant les images de ces trois chemins dans $K$, en utilisant la
convexit\'e de $K$ dans $\G$ et la caract\'erisation des
g\'eod\'esiques combinatoires par la suite des murs travers\'es, on
voit que le mur $M_i$ coupe une ar\^ete issue de $x'_0$, et s\'epare
$x_0$ de $x'_0$. 
\endproof

Comme ci-dessus, il y a une unique copie convexe de $L$ dans $K$
passant par $x'_0$, coup\'ee par les murs $M_1,\ldots M_r$: nous la
noterons $L'$ et $C'$ sera la cellule de $P$ dont le $1$--squelette est
$L'$.
 
Les centres m\'etriques des cellules $C$ et $C'$ sont des points
$\hat{C}$ et $\hat{C'}$ de $M_1\cap \ldots \cap M_r$; par convexit\'e,
la g\'eod\'esique $\gamma$ qui les joint est aussi dans cette
intersection.  D'autre part, la g\'eod\'esique joignant $\hat{C}$ \`a
$x_0$ ne coupe que les murs de $C$, donc pas $M_0$ (d'apr\`es le lemme
\ref{lem:lemme_un}). Un r\'esultat analogue \'etant vrai pour $x'_0$,
et $M_0$ s\'eparant $x_0$ de $x'_0$, la g\'eod\'esique $\gamma$ doit
couper $M_0$.  Comme $\gamma \not\subset M_0$, l'intersection de
$\gamma$ avec $M_0$ ne contient qu'un point $p$.

Soit $D$ la cellule de $P$ engendr\'ee par $\gamma$ juste apr\`es
$\hat{C}$; comme $\gamma\cap C=\hat{C}$, la cellule $D$ contient $C$
comme face stricte.

\medskip
\noindent{\bf Fait 2}\qua Le point $p$ appartient \`a $D$.

\dem 
Par l'absurde, supposons que $p$ n'est pas dans $D$.  Alors $\gamma$
ressort de $D$ par un point $q$ de son bord; ce point est dans
$M_1\cap \ldots \cap M_r$. Comme $\gamma$ reste dans l'int\'erieur de
$D$ entre $\hat{C}$ et $q$, ces deux points ne peuvent \^etre sur une
m\^eme face du bord de $D$. Soit $F$ la face stricte de $D$
engendr\'ee par $q$; cette cellule paire est coup\'ee par les murs
$M_i,1\leq i\leq r$, donc invariante par les r\'eflexions $\sigma
_{M_i}$, tout comme $C$.

Montrons que $F$ est disjointe de $C$. Si $F$ contenait un sommet de
$C$, elle contiendrait toutes ses images par le groupe d'isom\'etrie
de $D$ engendr\'ee par les r\'eflexions $\sigma _{M_i},1\leq i\leq
r$. Mais ce groupe est (simplement) transitif sur l'ensemble des
sommets de $C$.  Donc $F$ contiendrait tous les sommets de $C$,
autrement dit $C$ elle-m\^eme. Mais alors $\hat{C}$ et $q$ seraient
dans une m\^eme face $F$ du bord de $D$, ce qui n'est pas.

Consid\'erons une g\'eod\'esique combinatoire $c$ de $x_0\in C$ \`a un
sommet de $F$, de longueur minimale.  Par convexit\'e de $i(D)$, on a
$c\subset D$.  Comme $C\cap F =\emptyset$, la longueur de $c$ est non
nulle: donc $c=(b,\ldots )$, o\`u $b$ est une ar\^ete de $D$ issue de
$x_0$. Comme $\sigma _{M_i}$ pr\'eserve $F$, il est \'evident, par
minimalit\'e, que $b\not\subset C$.  Donc, d'apr\`es le lemme
\ref{lem:lemme_un}, $M(b)$ ne peut pas couper $C$.

Maintenant, le mur $M(b)$ ne peut pas non plus couper $F$: sinon en
appliquant $\sigma _{M(b)}$ au sous-segment de $c$ apr\`es $b$, on
trouverait une g\'eod\'esique de $x_0$ \`a un sommet de $F$, de
longueur inf\'erieure \`a celle de $c$, en contradiction avec la
minimalit\'e de celle-ci.
 
Il en r\'esulte que $M(b)$ s\'epare les deux cellules $C$ et $F$, donc
en particulier les deux points $\hat{C}$ et $q$. Alors $M(b)$ s\'epare
$\hat{C}$ et $\hat{C'}$, $C$ et $C'$, donc $x_0$ et $x'_0$.

Le mur $M(b)$ n'est pas le mur $M_0$: car celui-l\`a coupe la
g\'eod\'esique $\gamma$ dans $D$, alors que celui-ci la coupe en $p$,
suppos\'e ext\'erieur \`a $D$. Nous nous retrouvons avec un
\'el\'ement $K$ de ${\cal F}$ et une ar\^ete $b$ de $P$ contenant un
sommet de $K$, mais non contenue dans $K$, telle que $M(b)$ s\'epare
deux points de $K$: une contradiction avec le lemme
\ref{lem:lemme_trois}. Cette absurdit\'e prouve que $p\in D$.
\endproof

\smallskip 
Puisque l'ar\^ete $a_0$ est issue d'un sommet $x_0$ de $D$
et que le mur $M_0=M(a_0)$ recoupe $D$ (en $p$), le lemme
\ref{lem:lemme_un} entra\^{\i}ne que $a_0\subset D$. Alors la
sous-cellule $E$ de $D$ engendr\'ee par les ar\^etes $a_0,a_1,\ldots
a_r$ v\'erifie $St(x_0,K)=St(x_0,i(E))$, donc $K=i(E)$ d'apr\`es le
lemme \ref{lem:lemme_cinq}.  
\endproof

\bcoro 
Le morphisme de restriction de {\rm Aut}$(P)$ dans {\rm Aut}$(\G)$ est un
isomorphisme.  
\ecoro

\dem
Il suffit de montrer la surjectivit\'e. Si $\varphi$ est un
automorphisme de $\G$, d\'efinissons un automorphisme
$\overline{\varphi}$ de $P$ de la fa\c con suivante.  Pour une cellule
$C$ de $P$, consid\'erons l'\'el\'ement $K'$ de ${\cal F}$ d\'efini
par $K'=\varphi (i(C))$.  D'apr\`es la proposition pr\'ec\'edente, il
existe une (unique) cellule $C'$ dont le 1--squelette est $K'$.  Alors
il existe un unique isomorphisme poly\'edral de $C$ sur $C'$
prolongeant $\varphi |_{i(C)}$.

La collection d'isomorphismes poly\'edraux locaux $\overline{\varphi} _C$
ainsi obtenue se recolle pour donner l'automorphisme $\overline{\varphi}$.
\endproof

Ce corollaire termine la preuve du th\'eor\`eme
\ref{theo:meme_groupe_automorphisme}.  
\endproof

\subsection{Existence d'automorphisme non trivial fixant 
strictement un mur propre}
\label{sect:fixateur_strict_propre_nontrivial}

Un automorphisme de $P$ {\it fixe strictement} un mur $M$ de $P$ si 
et seulement s'il fixe $M$ (point par point) et pr\'eserve chacune des 
deux composantes connexes de $P\setminus M$.  

Le but de cette section est de donner des exemples de $P$ dont le
groupe Aut$^+(P)$, sous-groupe de Aut$(P)$ engendr\'e par les 
stabilisateurs stricts de murs propres est tr\`es gros.

\rem
(1)\qua L'automorphisme $f$ fixe strictement le mur $M$ si et seulement s'il
fixe point par point $M\cup a$, o\`u $a$ est une ar\^ete transverse
\`a $M$. Une condition \'equivalente est que $f$ fixe $V(M)$ point par
point.  Et un automorphisme de $P$ fixe strictement un mur $M$ \ssi\
l'automorphisme correspondant de $(X_P,\M_P)$ fixe strictement le mur
correspondant \`a $M$.

(2)\qua Soient $P^+$ et $P^-$ les adh\'erences des deux composantes
connexes de $P\setminus M$. Alors le sous-groupe de Aut$(P)$ form\'e
des automorphismes fixant strictement $M$ est le produit direct de
Fix$(P^+)$ et de Fix$(P^-)$.

\blemm\label{lemm:espace_amur_associe_verifie_Mprime}
Soit $P$ un complexe poly\'edral pair {\rm CAT}$(0)$. Alors
son espace \`a murs $(X_P,M_P)$ v\'erifie la propri\'et\'e {\rm(${\rm M}'$)}.
\elemm

\dem
Soit $f$ un automorphisme de $P$ fixant strictement un mur $M$ et $A$
une des deux moiti\'es de $X$ d\'efinies par $M$.  Soit $B$ une
moiti\'e de $X$ telle que $A\cap B$ et $(X\setminus A)\cap B$ sont non
vides. Notons $N$ le mur de $P$ dont le mur associ\'e sur $X$ est
$(B,X\setminus B)$. Alors on voit que $N$ contient des points
s\'epar\'es par $M$. Donc, par convexit\'e, $M\cap N$ est non vide. En
particulier, il existe une cellule $C$ de $P$ coup\'ee par $M$ et
$N$. Puisque $f$ fixe strictement $M$, elle vaut l'identit\'e sur $C$.
Donc $f$ fixe une ar\^ete transverse \`a $N$: $f$ pr\'eserve
globalement $N$, ainsi que les deux composantes connexes de
$X\setminus N$.  
\endproof

 Nous allons \'etudier le cas o\`u $P$ est la r\'ealisation
g\'eom\'etrique de Davis--Mous\-song d'un syst\`eme de Coxeter.

Soit $(W,S)$ un syst\`eme de Coxeter. Nous noterons $N=N(W,S)$ le nerf
fini de $(W,S)$. Nous munissons la \psb\ $N'$ de $N$ d'une fonction
$m$, d\'efinie sur l'ensemble des milieux $\hat a$ des ar\^etes $a$ de
$N$ par la formule: $m(\hat{a})$ est l'ordre du produit $st$, avec $s$
et $t$ les r\'eflexions de $S$ correspondant aux extr\'emit\'es de
$a$.  Il est alors imm\'ediat que les automorphismes du graphe de
Coxeter de $(W,S)$ correspondent aux automorphismes de $N'$ qui
proviennent d'un automorphisme de $N$ et pr\'eservent la fonction $m$.

Notons $P=|W|_0$ la r\'ealisation de Davis--Moussong de $(W,S)$. On a
$P'= (W\times{(x_0*N')})/\sim$ (voir section
\ref{sect:exemple_syssteme_coxeter}), et nous noterons $[w,x]$ la
classe de $(w,x)$.  Les sommets de $P$ sont les points $[w,x_0]$ pour
$w\in W$. Nous identifierons un point $x$ de $x_0*N'$ avec son image
$[id,x]$ dans $P'$. En particulier, le link de $x_0$ dans $P'$
s'identifie avec $N'$.  L'action \`a gauche de $W$ sur le produit
passe au quotient, en une action simplement transitive sur les sommets
$wx_0$ de $P$.  Mais on peut aussi construire, \`a partir de $(W,S)$,
des \'el\'ements de Aut$(P)$ fixant le sommet $x_0$.

Soit $G(W,S)$ le groupe des automorphismes du diagramme de Coxeter de
$(W,S)$. Tout \'el\'ement $f$ de $G(W,S)$ agit sur $N'$ (en
pr\'eservant $m$), donc nous pouvons consid\'erer son prolongement
conique \`a $x_0*N'$, encore not\'e $f$. D'autre part, $f$ induit
naturellement un automorphisme du groupe $W$ (permutant $S$), que nous
noterons $\overline{f}$. Alors l'application $(w,x)\mapsto
(\overline{f}(w),f(x))$ est compatible avec $\sim$, donc induit un
automorphisme $\hat{f}$ de $P'$. On a $\hat{f}([w,x] )=[
\overline{f}(w),f(x)]$, donc $\hat{f}(x_0)=x_0 $, et $\hat{f}$ agit sur le
link de $x_0$ comme $f$ sur $N'$.  Enfin, $\hat{f}$ provient d'un
automorphisme de $P$ (car $f$ provient d'un automorphisme de $N$).

Nous obtenons ainsi une repr\'esentation fid\`ele de $G(W,S)$ dans
Aut$(P)$, d'image contenue dans le stabilisateur de $x_0$.  D'apr\`es
la formule $\hat{f} (w\cdot[w',x] )=\overline{f}(w)\cdot\hat{f}([w',x])$,
si $\overline{f}$ fixe point par point un sous-ensemble $T$ de $S$, alors
$\hat{f}$ commute avec l'action sur $P$ du
sous-groupe sp\'ecial engendr\'e par $T$.

\bdefi\label{def:facette_bloc}
 Soient $Q$ un complexe poly\'edral pair CAT$(0)$, et $a$ une ar\^ete
de $Q$.  La {\it facette} de $Q$ transverse \`a $a$ est la r\'eunion
des simplexes de $Q'$ (la \psb\ de $Q$) qui contiennent le milieu de
l'ar\^ete $a$, mais aucune de ses extr\'emit\'es. Nous la noterons
$\phi (a)$. Si $x_0$ est un sommet de $Q$, le {\it bloc de centre $x_0$}
est l'\'etoile de $x_0$ dans $Q'$. 
\edefi

\blemm
Le mur transverse \`a l'ar\^ete $a$ est la r\'eunion des facettes
$\phi (b)$, avec $b$ parall\`ele \`a $a$. Les deux blocs centr\'es sur
les extr\'emit\'es d'une ar\^ete $a$ ont pour intersection la facette
$\phi (a)$.  
\elemm

\dem
V\'erifions d'abord que $\phi (a)\subset M(a)$. Soit $\Delta$ un simplexe
de $\phi  (a)$. Consid\'erons la plus petite
cellule $C$ contenant $\Delta$: les sommets de $\Delta$ sont les centres
m\'etriques de certaines faces de $C$
contenant l'ar\^ete $a$. Donc chacun de ces sommets est invariant par
$\sigma (a,C)$: autrement dit $\Delta\subset M(a,C)$.

Pour achever de montrer la premi\`ere assertion, il suffit de prouver
que si $C$ est un poly\`edre pair et $a$ une ar\^ete de $C$, alors
$M(a,C)$ est contenu dans l'union des facettes $\phi (b)$, avec $b$
parall\`ele \`a $a$ dans $C$.  On raisonne par r\'ecurrence sur
$dim(C)$, la propri\'et\'e \'etant \'evidente en dimension $1$.

Soit $x$ un point de $M(a,C)$. Si $x=\hat{C}$ (le centre m\'etrique de
$C$), alors $x$ est dans toutes les facettes de $C$, en particulier
dans $\phi (a)$. Si $x\neq \hat{C}$, nous pouvons consid\'erer la
g\'eod\'esique de $\hat{C}$ \`a $x$, et la prolonger jusqu'au bord de
$C$, qu'elle touche en un point $y$. Comme $x$ et $\hat{C}$ sont dans
$M(a,C)$, le point $y$ est aussi dans ce mur. Cela signifie que $D$,
la face stricte de $C$ engendr\'ee par $y$, est coup\'ee par
$M(a,C)$. On applique alors l'hypoth\`ese de r\'ecurrence \`a $y\in
D$: il existe un simplexe $\Delta _y$ de la facette d'une ar\^ete $b$
de $D$ telle que $M(b,D)=M(a,C)\cap D$ qui contient $y$.  Alors $b$
est parall\`ele \`a $a$, et $x$ est dans $\Delta _x$, le joint de
$\Delta _y$ avec $\hat{C}$. Ceci conclut, car $\Delta _x$ est dans la
facette de $b$ dans $C$.

Pour la seconde assertion, soit $a$ une ar\^ete de $Q$
d'extr\'emit\'es $x_0$ et $y_0$. Un sommet de $Q'$ est joignable aux
extr\'emit\'es de $a$ si et seulement s'il est le centre m\'etrique
d'une face $C$ contenant $x_0$ et $y_0$. Ceci, par convexit\'e,
\'equivaut \`a dire que $C$ contient $a$, autrement dit $\hat{C}\in
\phi (a)$.
\cqfd

Si $(W,S)$ est un syst\`eme de Coxeter, nous appellerons {\it facette
de $(W,S)$ (au sens de Davis--Moussong)} l'\'etoile dans $N'$ d'un
sommet de $N$; si ce sommet correspond \`a la r\'eflexion $s$, nous
noterons $\phi _s$ cette facette. Le syst\`eme de Coxeter est dit {\it
rigide} si le fixateur dans $G(W,S)$ de toute facette de $(W,S)$ est
trivial. Tous les blocs de la r\'ealisation g\'eom\'etrique de
Davis--Moussong $P$ sont isomorphes au c\^one sur $N'$; l'intersection
de deux blocs centr\'es sur des sommets voisins de $P$ est donc une
facette de $(W,S)$.

\btheo\label{theo:coxeter_non_rigide}
Si $(W,S)$ est rigide, alors {\rm Aut}$(P)$ est
discret: c'est le produit semi-direct de $W$ et de $G(W,S)$. 

Supposons que $(W,S)$ n'est pas rigide, et que $(W,S)$ est
hyperbolique au sens de Gromov. Alors, pour tout automorphisme non
trivial $f$ de $G(W,S)$ fixant une facette $\phi _s$, le mur $M_s$
passant par $\phi _s$ est propre. De plus, il existe un automorphisme
$\varphi$ de $P$, dont la restriction \`a l'\'etoile de $x_0$ dans $P$
est $\hat{f}$, et qui fixe strictement $M_s$.  En particulier,
{\rm Aut}$^+(P)\neq\{1\}$, et {\rm Aut}$(P)$ est non discret.  
\etheo

\dem
Supposons d'abord que $(W,S)$ est rigide. Il s'agit de montrer que le
stabilisateur de $x_0$ est $G(W,S)$.

D'abord le fixateur de $St(x_0,P')$ dans Aut$(P)$ est trivial: car si
$F\in$ Aut$(P)$ fixe l'\'etoile d'un sommet dans $P'$, alors par
rigidit\'e $F$ fixe l'\'etoile de tout sommet voisin. 

Ensuite, si $F\in$ Aut$(P)$ fixe $x_0$, il induit un automorphisme du
link de $x_0$ dans $P'$ (isomorphe \`a $N'$), provenant d'un
automorphisme de $N$, et pr\'eservant la fonction $m$. En effet, cette
fonction a une interpr\'etation g\'eom\'etrique: $m(x)$ est simplement
le diam\`etre combinatoire du bord de la 2--face dont $x$ est le
centre. Il existe donc un $f\in G(W,S)$ tel que $\hat{f}$
co\"{\i}ncide avec $F$ sur le bloc de $P$ de centre $x_0$. D'apr\`es
la premi\`ere partie, $\hat{f}=F$.

\blemm
Soit $(W,S)$ un syst\`eme de Coxeter hyperbolique. Alors
l'en\-sem\-ble des $w\in W$ qui agissent trivialement au bord de $W$ est
un sous-groupe sp\'ecial fini $W_F$ tel que le syst\`eme
$(W_{S\setminus F},S\setminus F)$ est irr\'eductible, et tout
\'el\'ement de $F$ commute avec tout \'el\'ement de $S\setminus F$.
En particulier, si $W$ est irr\'eductible, alors $W$ agit fid\`element
sur son bord.
\elemm

\dem
Soit $G$ le sous-groupe de $W$ agissant trivialement sur $\bord{W}$ :
c'est un sous-groupe distingu\'e fini de $W$ (voir \cite{Cha}). En
tant que sous-groupe fini, $G$ est contenu dans un conjugu\'e d'un
sous-groupe sp\'ecial fini $W_T$. Mais comme $G$ est distingu\'e, on a
$G\subset W_T$, avec toujours $G$ distingu\'e dans $W$. En prenant
l'intersection des sous-groupes sp\'eciaux finis contenant $G$, on
trouve un sous-groupe sp\'ecial fini $W_F$ contenant $G$, et tel que,
pour tout $t\in F$, il existe $g\in G$ tel que $t$ apparaisse dans une
\'ecriture de longueur minimale de $g$. Pour $s$ n'appartenant pas \`a
$F$ et $g\in G$, on a $s.g.s\in G\subset W_F$. Donc $s$ commute avec
tous les \'el\'ements $t$ de $F$ apparaissant dans une \'ecriture
g\'eod\'esique de $g$. On en d\'eduit que tout \'el\'ement de $F$
commute avec tout \'el\'ement de $S\setminus F$.

Il reste \`a montrer que $(W_{S\setminus F},S\setminus F)$ est
irr\'eductible. Supposons que $S\setminus F=T_1\cup T_2$, avec
$T_1\cap T_2 = \emptyset$ et tout \'el\'ement de $T_1$ commute avec
tout \'el\'ement de $T_2$.  On ne peut avoir $W_{T_1}$ et $W_{T_2}$
infinis, puisque $W$ est hyperbolique et contient $W_{T_1}\times
W_{T_2}$.  Si par exemple $W_{T_1}$ est fini, il commute \`a
$W_{T_2\cup F}$, donc agit trivialement au bord: d'o\`u $T_1\subset
W_{T_1}\subset G\subset W_F$, et donc $T_1=\emptyset$.  
\cqfd

Supposons maintenant que $(W,S)$ est non rigide. Soient $f\in
G(W,S)$ et $s\in S$ tels que $f$ est non trivial et $f$ fixe $\phi _s$
(point par point). Le fait que $M_s$ soit un mur propre de $P$
r\'esulte du lemme pr\'ec\'edent. En effet, si $M_s$ n'est pas propre,
comme $s$ permute les deux demi-espaces d\'efinis par $M_s$, le bord
de $M_s$ est \'egal \`a tout le bord de $P$, donc $s$ agit
trivialement sur le bord de $P$. Par le lemme, $s$ appartient \`a $F$,
et son \'etoile est \'egale \`a tout le nerf fini de $(W,S)$, ce qui
contredit la non-trivialit\'e de $f$.

Comme $f$ fixe la facette $\phi _s$, $\overline{f}$ fixe tous
les $t\in S$ tels que $m_{s,t}<\infty$. Donc $\hat{f}$ commute \`a
tout produit de telles r\'eflexions. Comme d'autre part $\hat{f}$ fixe
l'ar\^ete transverse \`a $M_s$ passant par $x_0$, c'est donc que
$\hat{f}$ fixe toutes les ar\^etes de la forme $w.a_s$, avec $w\in
W_{T_s}$, o\`u $T_s=\{t\in S\,/\,m_{s,t}<\infty\}$, et $a_s$ est
l'ar\^ete de $P$ entre $x_0$ et $sx_0$.

Soient $\bord{\H}$ la r\'eunion de ces ar\^etes, $\H$ la r\'eunion des
chemins d'origine $x_0$ dans le $1$--squelette $\G$ de $P$, qui ne
traversent pas $\partial \H$, et ${\H}^c$ le sous-graphe de $\G$ r\'eunion
des ar\^etes non dans $\H$. La proposition
\ref{prop:technique_Coxeter} de la section suivante dit que ${\H}^c$
contient le demi-espace $A$ de $W$ d\'efini par $M_s$ et contenant
$s$. L'automorphisme $\hat{f}$ vaut l'identit\'e sur $\bord{\H }$ , et
$\H \cap {\H}^c$ est contenu dans $\bord{\H}$. Donc on peut d\'efinir
un automorphisme $\varphi$ de $\G $ qui co\"{\i}ncide avec
l'identit\'e sur ${\H}^c$, et avec $\hat{f}$ sur $\H$.  Comme $A$ est
contenu dans ${\H}^c$, $\varphi$ vaut l'identit\'e sur le demi-espace
$A$, donc fixe strictement le mur $M_s$.  Enfin, $\hat{f}$ agit sur
l'ensemble $ \bigcup _{t\in S\setminus\{s\}}a_t$ des ar\^etes \`a la
fois dans $\H$ et dans l'\'etoile de $x_0$, comme $f$ sur
$S\setminus\{s\}$. Donc $\varphi$, qui co\"{\i}ncide avec $\hat{f}$
sur l'\'etoile de $x_0$, est non trivial.  
\cqfd

\subsection{Un r\'esultat technique sur les groupes de Coxeter}
\label{sect:technique_Coxeter}

Soit $(W,S)$ un syst\`eme de Coxeter, notons $1$ son \'el\'ement
neutre et $\G=\G(W,S)$ son graphe de Cayley.

Si $t\in S$ et $T\subset S$, nous noterons $a_t$ l'ar\^ete de
$\G(W,S)$ entre 1 et $t$ et $\G _T$ le sous-graphe de $\G$ r\'eunion
des ar\^etes reliant deux sommets de $\G $ appartenant au sous-groupe
sp\'ecial $W_T$ engendr\'e par $T$. Alors $\G _T$ est isomorphe \`a
$\G (W_T,T)$, et c'est un sous-graphe convexe de $\G $ (voir
\cite{Bourb}).  On peut aussi voir $\G _T$ comme la r\'eunion des
chemins de ${\cal G}$ d'origine $1$, et dont toutes les ar\^etes ont
un label dans $T$ (i.e.~sont de la forme $wa_t$, avec $t\in T$).

Pour $s\in S$ quelconque, soit $T_s$ la partie de $S$ form\'ee des
r\'eflexions $t$ telles que $m_{s,t}<\infty$.  Notons alors $\H =\H_s$
la r\'eunion des chemins $c$ de ${\cal G}(W,S)$ d'origine 1, et
n'empruntant que des ar\^etes de la forme $wa_t$, avec $t\neq s$, ou
de la forme $wa_s$, avec $w\not\in W_{T_s}$.  Introduisons enfin
$A=A_s$, l'ensemble des \'el\'ements de $W$ s\'epar\'es de $1$ par le
mur $M_s$ de $s$.

Notre but est de montrer que $\H$ et $A$ sont disjoints, ce qui est le
r\'esultat voulu dans la preuve du th\'eor\`eme
\ref{theo:coxeter_non_rigide}.

Commen\c{c}ons par donner une description plus constructive de $\H$.
Posons
$$\K_1=\{1\},\H _1=\G _{S\setminus \{s\}}, \ldots ,$$ $$\K_{n+1}=\bigcup _
{w\in \H _{n},w\not\in W_{T_s}}w.\G _{T_s}, \H _{n+1}=\bigcup _{w\in
  {\cal K}_{n+1}}w.\G _{{S\setminus \{s\}}}.$$
Alors $\H$ est la r\'eunion croissante des $\H _n$.

\blemm\label{lem:new_lemme_un} 
Soient $w\in \H $ et $n=n(w)$ le plus petit indice $i$ tel que 
$w\in \H _i$. Alors il existe deux suites
$v_1^-,\ldots ,v_{n-1}^-$ et $v_1^+,\ldots ,v_{n-1}^+$ d'\'el\'ements
de $W$ appartenant \`a $\H$, et une suite $M_1,\ldots ,M_{n-1}$ de
murs de $(W,S)$ tels que, pour $1\leq i< n$,
\begin{itemize}
\item 
$n(v_i^-)= i$ et $n(v_i^+)= i+1$;
\item $v_i^-$ et $v_i^+$ sont congrus modulo $W_{T_s}$, $v_{i-1}^+$ et
  $v_i^-$ sont congrus modulo $W_{{S\setminus \{s\}}}$ (en posant
  $v_0^+=1$), et $v_{n-1}^+$ est congru \`a $w$ modulo $W_{{S\setminus
      \{s\}}}$;
\item $M_i$ est transverse \`a l'ar\^ete de type $s$ d'origine
  $v_i^-$, et s\'epare $v_i^-.\G _{{S\setminus \{s\}}}$ de $v_i^+.\G
  _{{S\setminus \{s\}}}$.
\end{itemize}
\elemm

\begin{figure}[htbp]
\cl{\relabelbox\small
\epsfysize 5.25cm\epsfbox{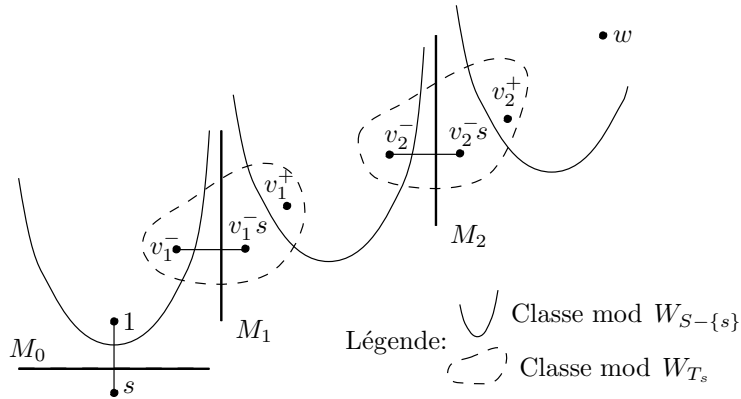}
\adjustrelabel <-6pt, 2pt> {w}{$w$}
\adjustrelabel <-2pt, 0pt> {1}{$1$}
\adjustrelabel <-2pt, 0pt> {s}{$s$}
\relabel {M0}{$M_0$}
\adjustrelabel <0pt, -2pt> {M1}{$M_1$}
\adjustrelabel <0pt, -2pt> {M2}{$M_2$}
\adjustrelabel <0pt, 0pt> {v}{$v_1^-$}
\adjustrelabel <0pt, -1pt> {vs}{$v_1^-\!s$}
\relabel {v1p}{$v_1^+$}
\adjustrelabel <0pt, -2pt> {v2m}{$v_2^-$}
\adjustrelabel <0pt, -2pt> {v2s}{$v_2^-\!s$}
\adjustrelabel <0pt, -2pt> {v2p}{$v_2^+$}
\adjustrelabel <0pt, 6pt> {Lgende}{L\'egende:}
\relabel {ClasseWt}{Classe mod $W_{T_s}$}
\relabel {ClasseWs}{Classe mod $W_{S-\{s\}}$}
\endrelabelbox}
\caption{Description constructive de $\H$}\label{coxtech1}
\end{figure}

\dem 
Par r\'ecurrence sur $n$. Si $n=1$, il n'y a rien \`a
d\'emontrer.

Supposons donc $n>1$. Comme $w\in \H _n$, il existe un $v_{n-1}^+$ de
$\K_n$ congru \`a $w$ modulo $W_{{S\setminus \{s\}}}$; puis il existe
un $v_{n-1}^-$ de ${\cal H}_{n-1}$ auquel $v_{n-1}^+$ est congru
modulo $W_{T_s}$. Quitte \`a multiplier $v_{n-1}^-$ et $v_{n-1}^+$ par
des \'el\'ements convenables de $W_{T_s\setminus \{s\}}$ (ce qui ne
change ni les classes modulo $W_{T_s}$, ni les classes modulo
$W_{{S\setminus \{s\}}}$), on peut supposer que
$d_W(v_{n-1}^-.W_{T_s\setminus \{s\}},v_{n-1}^+.W_{T_s\setminus
  \{s\}})=d_W(v_{n-1}^-,v_{n-1}^+)$.

Si $v_{n-1}^-$ \'etait dans un $\H _{i}$ avec $i<n-1$, $w$ serait dans
$\H _{n-1}$, en contradiction avec la d\'efinition de $n$. De m\^eme,
$v_{n-1}^+\in {\cal H}_{n}\setminus \H _{n-1}$ (en particulier,
$v_{n-1}^-.W_{{S\setminus \{s\}}}\neq v_{n-1}^+.W_{{S\setminus
    \{s\}}}$).  Donc, si on compl\`ete les suites fournies par la
r\'ecurrence appliqu\'ee \`a $v_{n-1}^-$ \`a l'aide de $v_{n-1}^-$ et
$v_{n-1}^+$ d'une part, et d'autre part \`a l'aide du mur $M_{n-1}$
fourni par le lemme suivant d'autre part, on obtient le r\'esultat au
rang $n$.  
\cqfd

\blemm\label{lem:new_lemme_deux}
Supposons que $v^-.W_{T_s}=v^+.W_{T_s}$ et $v^-.W_{{S\setminus
    \{s\}}}\neq v^+.W_{{S\setminus \{s\}}}$.  Si de plus
$d_W(v^-.W_{T_s\setminus \{s\}},v^+.W_{T_s\setminus
  \{s\}})=d_W(v^-,v^+)$, alors le mur transverse \`a l'ar\^ete
$v^-.a_s$ issue de $v^-$ s\'epare $v^-.\G _{{S\setminus \{s\}}}$ de
$v^+.\G _{{S\setminus \{s\}}}$. 
\elemm

\dem
Quitte \`a multiplier par l'inverse de $v^-$, on peut supposer
$v^-=1$. On a alors $v=v^+\in W_{T_s}\setminus W_{{S\setminus
\{s\}}}$, et $v$ est l'\'el\'ement de plus petite longueur dans sa
classe modulo $W_{T_s\setminus \{s\}}$ (cette longueur est non nulle,
sinon $v$ serait dans $W_{T_s\setminus \{s\}}$, donc dans
$W_{{S\setminus \{s\}}}$).  En particulier, toute g\'eod\'esique de $1$
\`a $v$ passe par $a_s$, et le mur $M_s$ transverse \`a $a_s$ s\'epare
$1$ de $v$.

La convexit\'e de $\G _{{S\setminus \{s\}}}$ l'emp\^eche d'\^etre
coup\'ee par le mur $M_s$.  Supposons que $M_s$ soit transverse \`a
une ar\^ete de $v.{\cal G}_{{S\setminus \{s\}}}$. Cela signifie qu'il
existe $t\in S\setminus \{s\}$ et $w\in W_{{S\setminus \{s\}}}$ tels
que $s(vw)=(vw)t$. Donc $v^{-1}sv=wtw^{-1}$: par convexit\'e des
sous-groupes sp\'eciaux, la r\'eflexion $v^{-1}sv$ est donc dans
$W_{T_s}\cap W_{{S\setminus \{s\}}}= W_{T_s\setminus \{s\}}$. Alors
l'\'el\'ement $v'=v(v^{-1}sv)=s.v$ est congru \`a $v$ modulo
$W_{T_s\setminus \{s\}}$, mais il est de longueur 1 de moins que $v$,
puisque toute g\'eod\'esique de $1$ \`a $v$ commence par $s$.  Ceci
contredit la minimalit\'e suppos\'ee de $|v|$.  
\cqfd

Nous allons montrer que, vus dans la r\'ealisation de Davis-Moussong
$P$ de $(W,S)$, les murs apparaissant dans le lemme
\ref{lem:new_lemme_un} sont disjoints, et ne s\'eparent pas deux
points de $A$. Pour cela, nous faisons agir $W$ sur un certain arbre.

Soient $s,t$ dans $S$ tels que $m_{s,t}=\infty$. Alors $W$ est le
produit amalga\-m\'e $W_{{S\setminus \{s\}}} *_{W_{{S\setminus
\{s,t\}}}}W_{{S\setminus \{t\}}}$.  Consid\'erons le graphe biparti
$\T_{s,t}$ ayant un sommet de type $s$ pour chaque classe de $W$
modulo $W_{{S\setminus \{s\}}}$, un sommet de type $t$ pour chaque
classe de $W$ modulo $W_{{S\setminus \{t\}}}$, avec une ar\^ete entre
une classe modulo $W_{{S\setminus \{s\}}}$ et une classe modulo
$W_{{S\setminus \{t\}}}$ lorsque ces deux classes ne sont pas
disjointes. Notons que si $w$ appartient \`a $uW_{{S\setminus
\{s\}}}\cap vW_{{S\setminus \{t\}}}$, alors $wW_{{S\setminus
\{s,t\}}}$ est contenu dans $uW_{{S\setminus \{s\}}}\cap
vW_{{S\setminus \{t\}}}$.  La convexit\'e des sous-groupes sp\'eciaux
entra\^{\i}ne alors que $wW_{{S\setminus \{s,t\}}}$ est \'egal \`a
$uW_{{S\setminus \{s\}}}\cap vW_{{S\setminus \{t\}}}$. Ainsi, les
ar\^etes de $\T_{s,t}$ correspondent bijectivement aux classes
de $W$ modulo $W_{{S\setminus \{s,t\}}}$.

Le groupe $W$ agit sur $\T_{s,t}$ par multiplication \`a gauche.
Cette action est transitive sur les ar\^etes de $\T _{s,t}$, le
stabilisateur de $x_s=W_{{S\setminus \{s\}}}$ est $W_{{S\setminus
\{s\}}}$, le stabilisateur de $x_t=W_{{S\setminus \{t\}}}$ est
$W_{{S\setminus \{t\}}}$, et le stabilisateur de l'ar\^ete joignant
ces deux sommets est $W_{{S\setminus \{s,t\}}}$. Il r\'esulte alors
de la th\'eorie de Bass--Serre \cite{Ser} que $\T _{s,t}$ est un arbre.

\blemm\label{lem:new_lemme_trois} 
Vus dans $P$, les murs $M_{i-1}$ et $M_i$ apparaissant dans le lemme
\ref{lem:new_lemme_un} sont disjoints. Le mur $M_{i}$ ne s\'epare pas
$v^-_{i-1}.a_s$ de $v^-_{i}$. Enfin $M_1\cap M(a_s)=\emptyset$.
\elemm

\dem
Notons que $v_{i-1}^+$ et $v_i^-$ ne peuvent \^etre dans la m\^eme
classe modulo $W_{T_s}$ (sinon $w\in \H _{n-1}$). Par
$W$--homog\'en\'eit\'e, il suffit donc de montrer le r\'esultat suivant
(lequel donne du m\^eme coup la derni\`ere partie du lemme).

Soient $v^-\in W_{{S\setminus \{s\}}}\setminus W_{T_s}$ et $M$ le mur
d'une r\'eflexion $u$ de $W_{T_s}$ ne coupant pas $\G _{{S\setminus
\{s\}}}$. Alors $M$ est disjoint du mur $M'$ transverse \`a $v^-.a_s$.

Puisque $v^-\not\in W_{T_s}$, il existe un $t\in S$ tel que
$m_{s,t}=\infty$ et $v^-\not\in W_{S\setminus \{t\}}$. Nous
raisonnons en consid\'erant l'action de $W$ sur l'arbre $\T _{s,t}$.

D'abord, dire que $M\cap \G _{{S\setminus \{s\}}}= \emptyset$, c'est
dire que $u\not\in W_{S\setminus \{s\}}$.  Autrement dit, $u.x_s\neq
x_s$. De m\^eme, si $u'$ est la r\'eflexion par rapport \`a $M'$, on a
$u'.x_s\neq x_s$.

Soit $y_t$ la classe \`a gauche de $v^-$ modulo $W_{S\setminus \{t\}}$. 
Alors $u'.y_t=y_t$. D'autre part, comme $T_s\subset S\setminus
\{t\}$, on a aussi $u.x_t=x_t$.

Enfin, $y_t\neq x_t$, et $x_s$ est li\'e dans $\T _{s,t}$ \`a $x_t$ et
$y_t$.

\begin{center}
\mbox{\relabelbox\small\epsfysize 1.5cm\epsfbox{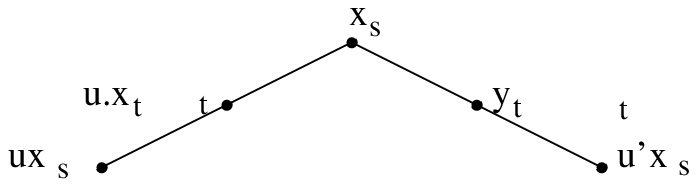}
\relabel {x}{$x_s$}
\relabel {y}{$y_t\!=\!u'.y_t$}
\adjustrelabel <-6pt, 0pt> {u.x}{$u.x_t\!=\!x_t$}
\relabel {ux}{$u.x_s$}
\relabel {u'x}{$u'.x_s$}
\endrelabelbox}
\end{center}

Donc le produit $u'u$ agit comme une translation non triviale de
l'arbre $\T _{s,t}$, et est n\'ecessairement d'ordre infini.  Or, si
les murs des deux r\'eflexions $u$ et $u'$ se coupaient dans $P$, le
produit $u'u$ aurait un point fixe, donc devrait \^etre d'ordre fini
(l'action de $W$ sur $P$ est propre).

En fait, non seulement $M\cap M'=\emptyset$, mais de plus $M'$ ne
s\'epare pas 1 d'une ar\^ete $a$ transverse \`a $M$ et contenue dans
$\G _{T_s}$ (ce qui ach\`eve de prouver le lemme). Car si c'\'etait le
cas, par convexit\'e, $M'$ serait transverse \`a une ar\^ete $a'$ de
$\G _{T_s}$, et $u'$ serait une r\'eflexion de $W_{T_s}$. Donc $u'$
fixerait $x_t$. Comme $u'$ fixe d\'ej\`a $y_t$, elle fixerait l'unique
sommet de type $s$ li\'e \`a la fois \`a $x_t$ et \`a $y_t$, c'est \`a
dire $x_s$. Or nous avons vu que ce n'\'etait pas le cas.  
\cqfd

\bprop\label{prop:technique_Coxeter}
 $\H$ est disjoint de $A$.
\eprop

\dem 
Si $w\in \H $, appliquons le lemme \ref{lem:new_lemme_un} pour trouver
une suite $s=v^-_0,v^-_1,\ldots ,v^-_{n-1},$ $v^-_n=w$ et une suite de
murs $M_0=M_s,M_1,\ldots ,M_{n-1}$ tels que $M_{i}$ s\'epare $v^-_{i}$
de $v^-_{i+1}$ et $M_i$ est transverse \`a $v^-_i.a_s$.  D'apr\`es le
lemme \ref{lem:new_lemme_trois}, les murs $M_i$ et $M_{i+1}$ sont
disjoints, et $M_{i+1}$ ne s\'epare pas $v^-_{i+1}$ de $M_i$ (voir
figure \ref{coxtech1}).  Soit $A_i$ la moiti\'e de $W$ d\'efinie par
$M_i$ contenant $v^-_i$.  Il est maintenant imm\'ediat que la suite
des moiti\'es $A_i$ est (strictement) croissante, avec $A_0=A$, et
$w\not\in A_{n-1}$. Donc $w\not\in A$.  
\cqfd

\subsection{Exemples de complexes poly\'edraux pairs CAT$(-1)$}
 \label{sect:exemples_pairs}

(1)\qua 
Soient $k$ un entier pair avec $k\geq 4$, et $L$ le graphe
d'incidence d'un plan projectif sur un corps fini, ou plus
g\'en\'eralement n'importe quel immeuble \'epais fini de rang $2$
v\'erifiant la condition de Moufang (voir \cite{Ron}).  Cette
condition (plus le fait que $L$ soit \'epais) implique en particulier
que le fixateur de l'\'etoile d'un sommet de $L$ est non trivial. Donc
Aut$^+(A(k,L))$ est non trivial, d\`es que $W(k,L)$ est hyperbolique
(au sens de Gromov), c'est-\`a-dire si $k\geq 6$ ou $k=4$ et $L$ n'est
pas de type $A_1\times A_1$. Ceci concerne donc l'immeuble de Bourdon
$I_{p,q}$, avec $p$ pair, $p\geq 6$ et $q\geq 3$.

\medskip
\noindent (2)\qua 
\'Etant donn\'e un poly\`edre pair $C$, nous allons montrer 
comment construire un complexe poly\'edral pair CAT$(-1)$
ayant un gros groupe d'automorphismes, et dont toute cellule 
maximale est isomorphe (combinatoirement) \`a $C$.

\bprop 
Pour tout poly\`edre pair $C$, il existe un complexe poly\'edral pair
localement compact {\rm CAT}$(-1)$, dont les cellules maximales sont 
combinatoirement isomorphes \`a $C$, admettant un
groupe discret cocompact d'auto\-mor\-phis\-mes, et dont le groupe des
automorphismes engendr\'e par les fixateurs stricts de murs propres 
est non d\'enombrable.  Si $C$ n'est pas combinatoirement un produit, alors
on peut de plus supposer que tous les murs sont propres.
\eprop

\dem
Soit $(W,S)$ le syst\`eme de Coxeter fini associ\'e \`a $C$ par la
proposition \ref{prop:construction_polyedre_pair}. Consid\'erons une
fonction $\overline{n}$ de $S$ dans l'ensemble des entiers strictement
positifs, telle que, si $\overline{n} (s) >1$ et $\overline{n} (t)
>1$, on a $m_{s,t}>2$ (c'est-\`a-dire $s$ et $t$ sont li\'es par une
ar\^ete dans le graphe de Coxeter de $(W,S)$). Nous noterons
$K_{\overline{n} }$ le sous-graphe complet du graphe de Coxeter de
$(W,S)$ dont les sommets $s$ v\'erifient $\overline{n} (s)
>1$. Remarquons que par le th\'eor\`eme de classification des
syst\`emes de Coxeter fini (voir par exemple \cite[p.193]{Bourb}),
$K_{\overline{n}}$ est r\'eduit \`a un seul sommet ou \`a une seule
ar\^ete.

D\'efinissons $(\overline{W},\overline{S})$, l'unique syst\`eme de
Coxeter tel qu'il existe une application $\tau \co \overline{S}\ra S$
avec
\begin{itemize}
\item[i)]
 $\tau ^{-1}(\{s\} )$ poss\`ede $\overline{n} (s)$  \'el\'ements;
\item[ii)]
si $\overline{s}\neq\overline{t}$, ou bien $\tau (\overline{s})=\tau
(\overline{t})$, et dans ce cas $m_{\overline{s},\overline{t}}=
\infty$, ou bien $\tau (\overline{s})\neq \tau (\overline{t})$, et
dans ce cas $m_{\overline{s},\overline{t}}=m_{\tau (\overline{s}),\tau
(\overline{t})}$.
\end{itemize}

\begin{figure}[htbp]
\cl{\relabelbox\small
\epsfysize 3cm\epsfbox{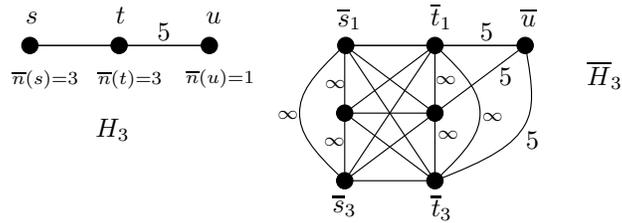}
\adjustrelabel <-1.5pt, 0pt> {oo}{$\scriptstyle\infty$}
\adjustrelabel <-1.5pt, 0pt> {oo1}{$\scriptstyle\infty$}
\adjustrelabel <-1.5pt, 0pt> {oo2}{$\scriptstyle\infty$}
\adjustrelabel <-1.5pt, 0pt> {oo3}{$\scriptstyle\infty$}
\adjustrelabel <-1.5pt, 0pt> {oo4}{$\scriptstyle\infty$}
\adjustrelabel <-1.5pt, 0pt> {oo5}{$\scriptstyle\infty$}
\relabel {5}{5}
\relabel {51}{5}
\relabel {52}{5}
\relabel {53}{5}
\relabel{H}{$H_3$}
\relabel{HH}{$\overline H_3$}
\adjustrelabel <-2pt, 0pt> {s}{$s$}
\adjustrelabel <-2pt, 0pt> {t}{$t$}
\adjustrelabel <-2pt, 0pt> {u}{$u$}
\adjustrelabel <-1pt, 0pt> {sb}{$\overline s_1$}
\adjustrelabel <-1pt, 0pt> {tb}{$\overline t_1$}
\adjustrelabel <-1pt, 0pt> {ub}{$\overline u$}
\adjustrelabel <4pt, 0pt> {s3}{$\overline s_3$}
\relabel {t3}{$\overline t_3$}
\adjustrelabel <-2pt, 1pt> {ns}{$\scriptstyle \overline n(s)=3$}
\relabel {nu}{$\scriptstyle \overline n(u)=1$}
\relabel {nt}{$\scriptstyle \overline n(t)=3$}
\endrelabelbox}
\caption{Exemple de syst\`eme de Coxeter 
hyperbolique non rigide}
\end{figure}

Il est imm\'ediat que $\tau$ s'\'etend en un homomorphisme de groupes
de $\overline{W}$ dans $W$, et est injective sur les parties
$\overline{T}$ de $\overline{S}$ telles que
$\overline{W}_{\overline{T}}$ est fini.  Donc les simplexes du nerf
fini de $(\overline{W},\overline{S})$ sont les parties de
$\overline{S}$ sur lesquelles $\tau$ est injective.  Les permutations
de $\overline{S}$ laissant $\tau$ invariante donnent des
automorphismes de $(\overline{W},\overline{S})$.  Si on suppose que
$\overline{n}$ atteint une valeur sup\'erieure ou \'egale \`a $3$, on
en d\'eduit que $(\overline{W},\overline{S})$ n'est pas rigide.

Montrons que $(\overline{W},\overline{S})$ est hyperbolique.  Si ce
n'est pas le cas, d'apr\`es Moussong, $(\overline{W},\overline{S})$
contient un sous-groupe sp\'ecial affine de rang au moins $3$, ou bien
deux sous-groupes sp\'eciaux infinis qui commutent. Dans le premier
cas, $\tau$ est n\'ecessairement injective sur le sous-groupe
sp\'ecial (car le graphe d'un tel syst\`eme de Coxeter ne contient pas
d'$\infty$), donc $(W,S)$ est infini, contradiction. Dans le deuxi\`eme
cas, un argument analogue au pr\'ec\'edent montre qu'il existe
$\overline{s_1}$, $\overline{s_2}$, $\overline{t_1}$ et
$\overline{t_2}$ tels que $m_{\overline{s_1},\overline{s_2}}=
m_{\overline{t_1},\overline{t_2}}=\infty$ et
$m_{\overline{s_i},\overline{t_j}}=2$ pour tous $i,j=1,2$. Mais alors
$\tau (\overline{s_1})=\tau (\overline{s_2})$ commute avec $\tau
(\overline{t_1})=\tau (\overline{t_2})$, en contradiction avec
l'hypoth\`ese de d\'epart sur $\overline{n}$.

Enfin, notons $P(C,\overline{n})$ la r\'ealisation g\'eom\'etrique de
Davis--Moussong de $(\overline{W},\overline{S})$. Alors les cellules
maximales de $P(C,\overline{n})$ correspondent aux sous-groupes
sp\'ec\-iaux finis maximaux de $(\overline{W},\overline{S})$, lesquels
sont tous isomorphes \`a $(W,S)$. Donc toutes les cellules maximales
de $P(C,\overline{n})$ sont isomorphes \`a $C$.

Supposons $(W,S)$ irr\'eductible. Par classification, son graphe de
Coxeter contient au plus une ar\^ete ayant un label pair.  Alors il
existe une application $\overline{n}$ telle que $\overline{n}(s)\geq
3$ si $s$ appartient \`a $K_{\overline{n}}$, et telle que s'il existe
une ar\^ete de label pair (diff\'erent de $2$), alors
$K_{\overline{n}}$ consiste en cette ar\^ete. Rappelons que si deux
sommets d'un graphe de Coxeter peuvent \^etre joints par un chemin
d'ar\^etes dont tous les labels sont impairs, alors les deux
r\'eflexions correspondantes sont conjugu\'ees dans le groupe de
Coxeter (voir \cite{Bourb}). Donc toute r\'eflexion de
$(\overline{W},\overline{S})$ est conjugu\'ee \`a un \'el\'ement de
$\tau^{-1}(K_{\overline{n}})$. Or le mur de toute r\'eflexion dans
$\tau^{-1}(K_{\overline{n}})$ est propre. Par cons\'equent, tout mur
est propre.  
\cqfd

Par exemple, lorsque $C$ est le polygone \`a $p=2k$ c\^ot\'es et 
$\overline{n}$ est constante \'egale \`a $q\geq 3$,
le poly\`edre $P(C,\overline{n})$ est l'immeuble de Bourdon $I_{p,q}$.

Lorsque $C$ est un cube de dimension $3$, $P(C,\overline{n})$ est le
produit d'un arbre r\'egulier par un carr\'e.

Lorsque $C$ est le poly\`edre pair du groupe $H_3$, d\'efinissons
$\overline{H_3}$ comme dans la figure pr\'ec\'edente.  Alors toutes
les 2--faces de $P(C,\overline{n})$ sont contenues dans $3$ copies de
$C$, sauf les d\'ecagones, qui ne sont contenus que dans une copie de
$C$.

\subsection{Automorphismes pr\'eservant le type de complexes
poly\-\'e\-draux pairs}
\label{sect:type}

Dans toute cette section, $P$ est un complexe poly\'edral pair
CAT$(0)$ dont toutes les cellules maximales (appel\'ees {\it chambres}
par la suite) sont isom\'etriques \`a une cellule $C$ fix\'ee (par
exemple, $P$ est un $(k,L)$--complexe, au sens de \cite{Hag,Ben}, voir
aussi \cite{BB}). La codimension des faces de $P$ est maintenant bien
d\'efinies.

\bdefi
Une {\it fonction type} de $P$ dans $C$ est une application
poly\-\'e\-drale $\tau \co P\ra C$ dont la restriction \`a chaque chambre de
$P$ est une isom\'etrie.
\edefi

\noindent {\bf Exemples}\qua 
(1)\qua Supposons que $C$ soit une cellule paire
de l'espace $E_{\chi}$ \`a courbure constante $\chi\leq 0$, dont les
faces de codimension $1$ font des angles di\`edres de la forme
$\frac{\pi}{n}$, avec $n\geq 2$.  Alors, par le th\'eor\`eme de
Poincar\'e (voir par exemple \cite{Har}), le sous-groupe $W(C)$ des
isom\'etries de $E_{\chi}$ engendr\'e par les r\'eflexions par rapport
aux faces de codimension 1 de $C$ est discret, et le quotient de
$E_{\chi}$ par $W(C)$ s'identifie naturellement \`a $C$. Cela signifie
que le pavage $P(C)$ de $E_{\chi}$ donn\'e par les $wC$, avec $w\in
W(C)$, admet une fonction type dans $C$.

(2)\qua Plus g\'en\'eralement, tout immeuble $P$ dont les appartements
sont iso\-m\'e\-tri\-ques \`a $P(C)$ admet une fonction type (on fixe une
certaine copie $A_0$ de $P(C)$ dans $P$, ainsi qu'une certaine chambre
$C_0$ de $A_0$, puis on consid\`ere la r\'etraction de $P$ sur $A_0$
bas\'ee en $C_0$, et on la compose par une quelconque fonction type
sur $A_0$).

(3)\qua Enfin, un arbre quelconque admet toujours une fonction type \`a
valeur dans l'une de ses ar\^etes.

Appelons {\it galerie de $P$} toute suite de chambres
$(C_0,C_1,\ldots,C_n)$ telles que $C_i\cap C_{i+1}$ contient une
cellule de codimension $1$. Nous laissons au lecteur le soin de
d\'emontrer la proposition suivante, qui ne servira pas dans ce texte.

\bprop 
Supposons que deux chambres de $P$ sont jointes par au moins une galerie.
Deux fonctions de type \'egales sur une chambre $C_0$ de $P$ sont
\'egales.  S'il est non vide, l'ensemble des fonctions types sur $P$
s'identifie avec l'ensemble (fini) des isom\'etries de $C_0$ sur
$C$. Dans ce cas, le link d'une face de codimension 2 de $P$ est
biparti.  

R\'eciproquement, si $P$ est de dimension $2$ avec $2$-cellules
r\'eguli\`eres, et si le link de chaque sommet de $P$ est un graphe
biparti connexe, alors l'ensemble des fonctions types sur $P$ est
non vide. \cqfd
\eprop

A partir de maintenant, nous supposerons que $P$ admet une fonction
type dans un poly\`edre pair $C$, et que deux chambres quelconques de
$P$ sont jointes par une galerie.

\bdefi
Nous noterons Aut$_0(P)$ le noyau de l'action par pr\'ecompos\-ition du
grou\-pe Aut$(P)$ sur l'ensemble des fonctions types de $P$ dans
$C$. Nous dirons que ses \'el\'ements {\it pr\'eservent le type}.
\edefi

\rem
Si $C'$ est isomorphe \`a $C$, un \'el\'ement de Aut$(P)$ pr\'eserve
le type dans $C$ si et seulement s'il pr\'eserve le type dans
$C'$. C'est ce qui justifie l'omission de $C$ dans la notation
Aut$_0(P)$. Remarquons que Aut$_0(P)$ est d'indice fini dans Aut$(P)$.

Notons Aut$_F (P)$ le sous-groupe caract\'eristique de Aut$(P)$
engendr\'e par les fixateurs de facettes (au sens de la d\'efinition
\ref{def:facette_bloc}).  Ses \'el\'ements seront appel\'es {\it
$F$-automorphismes}.  Notons $G_0$ et $G_1$ les sous-groupes de $G=$
Aut$(P)$ engendr\'es par les intersections avec Aut$_F(P)$ des
fixateurs de chambres d'une part, et des fixateurs de cellules de
codimension $1$ d'autre part.

Il est clair que $G_0\subset G_1\subset$ Aut$_F(P)$. Si $M$ est un mur 
propre de $P$,
son fixateur strict est dans $G_0$. Donc Aut$^+(P)$ est contenu dans
$G_0$. Si $f\in$ Aut$(P)$ fixe une face $F$
de codimension $1$ et envoie une chambre $C_2$ contenant $F$ sur
$C_1$, alors $f|_{C_2}$ commute avec la fonction type de $P$
restreinte \`a $C_1$ et \`a $C_2$.  Par connexit\'e par galeries de
$P$, $f$ pr\'eserve alors le type. Donc $G_1\subset \hbox{Aut}_0(P)$.
En r\'esum\'e, $\hbox{Aut}^+(P)\subset G_0 \subset G_1\subset
\hbox{Aut}_0(P)\cap$ Aut$_F(P)$.

Introduisons des propri\'et\'es de transitivit\'e, globales ou
locales:
\begin{itemize}\leftskip 10pt
\item[$(T_0)$]
 L'action de $G_0$ sur l'ensemble des chambres de $P$ est transitive.
\item[$(T_1)$]
 L'action de $G_1$ sur l'ensemble des chambres de $P$ est transitive.
\item[$(TL_1)$] 
Pour toute face $\sigma$ de codimension 1, le fixateur de $\sigma$ dans
Aut$_F (P)$ agit transitivement sur les chambres contenant $\sigma$.
\item[$(TL_0)$] 
Pour toute face $\sigma$ de codimension 1, le sous-groupe de 
Fix$(\sigma)\cap\hbox{Aut}_F (P)$ engendr\'e par les 
Fix$(\sigma)\cap\hbox{Aut}_F (P)\cap \hbox{Fix}(C)$, o\`u $C$ est une 
chambre de $P$, agit transitivement sur les chambres contenant $\sigma$.
\end{itemize}

Il est imm\'ediat que $(T_0)$ implique $(T_1)$ et $(TL_0)$ implique
$(TL_1)$. D'autre part:

\blemm\label{lem:tli_implique_ti}
Pour $i=0,1$, la condition $(TL_i)$ implique $(T_i)$, qui implique 
que $G_i=\hbox{\rm Aut}_0(P)\cap\hbox{\rm Aut}_F (P)$.
\elemm

\dem
La premi\`ere implication d\'ecoule de la connexit\'e par galerie de
l'en\-sem\-ble des chambres de $P$.  La deuxi\`eme de ce qu'un \'el\'ement
de Aut$_0(P)$ pr\'eservant une chambre la fixe n\'ecessairement.
\cqfd

Si $a$ et $b$ sont deux ar\^etes adjacentes \`a un sommet $x_0$ et
contenues dans un m\^eme polygone de $P$, nous noterons $m_{a,b}$ la
moiti\'e du nombre de c\^ot\'es de ce polygone. Nous obtenons ainsi
une fonction de l'ensemble des ar\^etes de $lk(x_0,P)$ dans l'ensemble
des entiers sup\'erieurs ou \'egaux \`a $2$. Il est clair qu'un
automorphisme $f$ de $P$ envoie la fonction $m$ du sommet $x_0$ sur la
fonction $m$ du sommet $f(x_0)$. Nous noterons $G_{x_0}$ le groupe des
automorphismes de $(lk(x_0,P),m)$ engendr\'e par les fixateurs de
facettes de $lk(x_0,P)$. Remarquons que si $P$ est de dimension $2$,
alors $m$ est constant (car tous les polygones ont le m\^eme nombre de
c\^ot\'es).

Voici maintenant deux propri\'et\'es de prolongement:
\begin{itemize}\leftskip 10pt
\item[$(P_0)$] 
  Pour tout sommet $x_0$ de $P$, tout \'el\'ement de
  $G_{x_0}$ s'\'etend \`a $P$.
\item[$(P^+)$] 
  Pour tout sommet $x_0$ de $P$, tout \'el\'ement de
  $G_{x_0}$ fixant une facette (voir d\'efinition \ref{def:facette_bloc})
 $\phi$ (transverse \`a une ar\^ete issue de $x_0$) s'\'etend \`a $P$ en un
  automorphisme fixant le mur $M$ passant par $\phi$ et fixant toute la
  moiti\'e de $P$ d\'efinie par $M$ et ne contenant pas $x_0$.
\end{itemize}

Nous avons maintenant des conditions permettant d'identifier
Aut$^+(P)$ et\break Aut$_0(P)\cap\hbox{Aut}_F (P)$, dans le cas o\`u $P$ est
de dimension 2.

\bprop\label{prop:aut_zero_egal_aut_plus} Soit $P$ un complexe
poly\'edral pair {\rm CAT}$(0)$ de dimension $2$ admettant un type et
dont deux chambres sont jointes par au moins une galerie.
\begin{enumerate}
\item
  Supposons que, pour tout sommet $x_0$ de $P$, le stabilisateur dans
  $G_{x_0}$ d'un sommet de $lk(x_0,P)$ agit transitivement sur les
  ar\^etes issues de ce sommet. Si $P$ v\'erifie $(P_0)$, alors
  {\rm Aut}$_0(P)\cap\hbox{\rm Aut}_F (P)= G_1$.
\item
  Supposons que, pour tout sommet $x_0$ de $P$, pour toute ar\^ete $a$
  issue de $x_0$, l'ensemble $E(a)$ des polygones de $P$ contenant $a$
  est de cardinal au moins trois, et que pour tout polygone $c$
  contenant $a$, le stabilisateur dans $G_{x_0}$ de $c$ agit
  transitivement sur $E(a)\setminus\{c\}$. Si $P$ v\'erifie $(P_0)$,
  alors $G_0 = G_1 =$ {\rm Aut}$_0(P)\cap\hbox{\rm Aut}_F (P)$.
\item
  Supposons que, pour tout sommet $x_0$ de $P$, toute ar\^ete $c$ de
  $lk(x_0,P)$, et tout $f\in G_{x_0}$ fixant $c$, on a une
  d\'ecomposition $f=f_1\circ f_2$, o\`u $f_1\in G_{x_0}$ fixe toute
  une facette de $lk(x_0,P)$ contenant une extr\'emit\'e de $c$, et
  $f_2\in G_{x_0}$ fixe toute la facette de $lk(x_0,P)$ contenant
  l'autre extr\'emit\'e de $c$. Supposons que la restriction d'un
  $F$-automorphisme fixant un sommet $x_0$ \`a $lk(x_0,P)$ est dans
  $G_{x_0}$. Si $P$ v\'erifie $(P^+)$ et si tous ses murs sont
  propres, alors {\rm Aut}$^+(P)=G_0$.
\end{enumerate}
\eprop

\dem
Pour la premi\`ere assertion, il suffit de remarquer que
l'hypoth\`ese, plus la propri\'et\'e $(P_0)$, entra\^{\i}nent la
propri\'et\'e $(TL_1)$. On applique alors le lemme
\ref{lem:tli_implique_ti} pr\'ec\'edent.

Pour la deuxi\`eme, par le lemme \ref{lem:tli_implique_ti}, il suffit
de v\'erifier que $P$ satisfait $(TL_0)$. Soit $a$ une ar\^ete de $P$
contenues dans deux polygones $c_1,c_2$. Fixons un sommet $x_0$ de $a$
et un troisi\`eme polygone $c$ contenant $a$ distinct de $c_1,
c_2$. Par hypoth\`ese, soit $f$ dans $G_{x_0}$ fixant $c$ et envoyant
$c_1$ sur $c_2$. La propri\'et\'e $(P_0)$ permet d'\'etendre $f$ en un
$F$-automorphisme de $P$, qui fixe $a$ et $c$, et envoie $c_1$ sur
$c_2$, ce qui montre $(TL_0)$.

Montrons la troisi\`eme assertion. Comme Aut$^+(P)$ est contenu dans
le groupe engendr\'e par les fixateurs de chambres, il suffit de
montrer que pour toute chambre $C$ de $P$, le groupe Fix$(C)$ est
contenu dans Aut$^+(P)\cap\hbox{Aut}_F (P)$. En fait, nous allons
montrer que si $\overline{f}$ est dans Fix$(C)$ et si $a$ et $b$ sont
deux ar\^etes du polygone $C$ adjacentes en un sommet $x_0$, alors il
existe $\overline{f}_a$ et $\overline{f}_b$ fixant strictement les
murs $M(a)$ et $M(b)$ tels que $\overline{f}=\overline{f}_b\circ
\overline{f}_a$.

\medskip
\noindent{\bf Affirmation 1}\qua Il existe un automorphisme
$\overline{f}_b$ de $P$ fixant $M(b)$ et toute la moiti\'e de $P$
d\'efinie par $M(b)$ ne contenant pas $x_0$, tel que $\overline{f}$
co\"{\i}ncide avec $\overline{f}_b$ sur l'ensemble des chambres de $P$
contenant l'ar\^ete $a$.

\dem 
L'automorphisme $\overline{f}$ fixe $x_0$ et $C$, donc induit un
automorphisme $f$ de $G_{x_0}$ fixant l'ar\^ete $c$ entre
les sommets du link correspondant aux ar\^etes $a$ et $b$ de $P$. Vu
l'hypoth\`ese sur $P$, il existe $f_a$ et $f_b$ dans $G_{x_0}$ tels
que $f=f_b\circ f_a$. D'apr\`es $(P^+)$, on peut prolonger ces deux
automorphismes locaux en \'el\'ements $\overline{f}_a$ et
$\overline{f}_b$ fixant strictement les murs $M(a)$ et $M(b)$ (ainsi
que les moiti\'es convenables). Maintenant l'\'egalit\'e
$\overline{f}=\overline{f}_b\circ \overline{f}_a$ sur l'\'etoile de
$x_0$ entra\^{\i}ne $\overline{f}=\overline{f}_b$ sur l'ensemble des
chambres de $P$ contenant l'ar\^ete $a$, puisque $\overline{f}_a$ agit
trivialement sur cet ensemble. \cqfd

Appelons {\it galerie g\'eod\'esique de $M(a)$ d'origine $(C,a)$}
toute galerie sans r\'ep\'e\-tition $(C_0,C_1,\ldots ,C_n)$ telle que
$C_0=C$, $C_0\cap C_1=a$, l'ar\^ete $C_i\cap C_{i+1}$ est parall\`ele
\`a $a$ et distincte de $a_{i-1}$. 
Comme $M(a)$ est un arbre, deux
galeries g\'eod\'esiques de $M(a)$ d'origine $(C,a)$ et de m\^emes
extr\'emit\'es sont \'egales.
Nous noterons $\delta(C_n)$ la longueur $n$ de cette galerie.
Soit alors ${\cal B}_n^+(C)$ l'ensemble des polygones $C'$ de $P$
qui sont extr\'emit\'es d'une galerie g\'eod\'esique de $M(a)$\break
\begin{center}
\mbox{\relabelbox\small\epsfysize 2cm\epsfbox{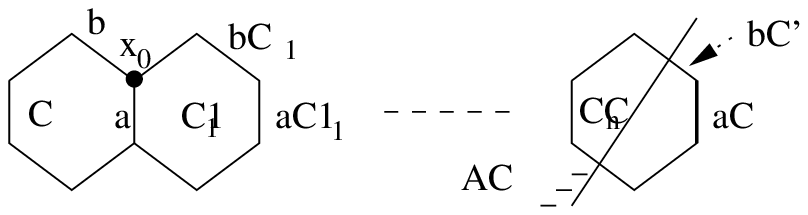}
\relabel {x}{$x_0$}
\relabel {b}{$b$}
\adjustrelabel <-2pt, 0pt> {a}{$a$}
\relabel {bC}{$b(C_1)$}
\adjustrelabel <0pt, 2pt> {bC'}{$b(C')$}
\relabel {C}{$C$}
\relabel {AC}{$A(C')$}
\adjustrelabel <-2pt, 0pt> {CC}{$C_n\!=\!C'$}
\relabel {aC}{$a(C')$}
\relabel {aC1}{$a(C_1)$}
\extralabel <-6.8cm, 0.9cm> $C_1$
\endrelabelbox}
\end{center}
d'origine $(C,a)$ de longueur au plus $n$.
  Nous noterons $a(C)$
l'ar\^ete $a$; pour $C'\in {\cal B}_n^+(C)$, avec $\delta (C')=n>0$,
soit $a(C')$ l'ar\^ete de $C'$ parall\`ele \`a $a$, non contenue dans
un polygone de ${\cal B}_{n-1}^+(C)$. Nous pouvons ensuite d\'efinir
$b(C')$ comme l'ar\^ete de $C'$ adjacente \`a $a(C')$, non s\'epar\'ee
de $b$ par $M(a)=M(a(C'))$.

Soient enfin $A(C')$ la moiti\'e ferm\'ee de $P$ d\'efinie par $M(b(C'))$
ne contenant pas $a(C')$, et ${\cal B}^+(C)$ l'union
des ${\cal B}_n^+(C)$.

\medskip
\noindent {\bf Affirmation 2}\qua 
Les murs $M(b(C'))$ sont deux \`a deux disjoints; la moiti\'e $A(C')$
contient strictement le mur $M(b(C''))$ (donc la chambre $C''$), d\`es
que la galerie g\'eod\'esique de $M(a)$ d'origine $(C,a)$ et
d'extr\'emit\'e $C''$ ne passe pas par
$C'$.

\dem
Prouvons d'abord que $M(b(C'))\cap M(b(C''))=\emptyset$, lorsque $C'$
et $C''$ sont deux chambres de ${\cal B}^+(C)$ telles que $C'\cap C''$
est une ar\^ete, et $\delta (C')\neq \delta (C'')$. Nous pouvons
supposer les notations telles que l'ar\^ete $a'$ commune \`a $C'$ et
$C''$ est l'ar\^ete $a(C')$ (autrement dit, $\delta (C')< \delta
(C'')$).  Pour all\'eger, nous notons alors $a''$, $b'$ et $b''$ les
ar\^etes $a(C'')$, $b(C')$ et $b(C'')$. Pour voir que $M(b')\cap
M(b'')=\emptyset$, il suffit de voir que les deux murs ont une
perpendiculaire commune (dans $C'\cup C''$): l'in\'egalit\'e CAT$(0)$
permet alors de conclure.

Soient $m$ le milieu de $a'$, et $p'$ (resp: $p''$) la projection
orthogonale de $m$ sur $M(b')$ (resp: $M(b'')$). Alors $p'$ s'obtient
comme l'intersection avec $M(b',C')$ de la g\'eod\'esique de $C'$
joignant $m$ \`a son image par $\sigma (b',C')$.  En particulier, $P'$
est \`a l'int\'erieur de $C'$, et $p''\neq p'$. Il reste \`a montrer
que la g\'eod\'esique de $P$ joignant $p'$ et $p''$ passe par $m$.

\begin{center}
\mbox{\relabelbox\small\epsfysize 2.5cm\epsfbox{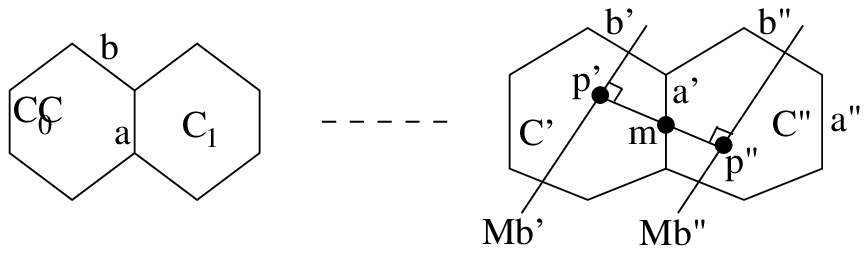}
\relabel {b}{$b$}
\relabel {C}{$C_1$}
\relabel {CC}{$C_0\!=\!C$}
\adjustrelabel <-1pt, 0pt> {a}{$a$}
\adjustrelabel <-1pt, 0pt> {a'}{$a'$}
\relabel {a"}{$a''$}
\relabel {b'}{$b'$}
\relabel {b"}{$b''$}
\relabel {C'}{$C'$}
\relabel {C"}{$C''$}
\adjustrelabel <-1pt, 0pt> {m}{$m$}
\relabel {p'}{$p'$}
\relabel {p"}{$p''$}
\relabel {Mb'}{$M(b')$}
\relabel {Mb"}{$M(b'')$}
\endrelabelbox}
\end{center}

La r\'eunion de $C'$ et $C''$ admet deux r\'eflexions orthogonales:
$\sigma _{a'}$ qui \'echange les extr\'emit\'es de $a'$, et $\rho
_{a'}$ qui fixe $a'$ en \'echangeant les deux chambres $C'$ et $C''$
(rappelons que $P$ admet un type).  Il est alors imm\'ediat que la
sym\'etrie centrale $\rho _{a'}\circ \sigma _{a'}$ envoie $M(b',C')$
sur $M(b'',C'')$ en fixant $m$, donc envoie le segment le $m$ \`a $p'$
sur le segment de $m$ \`a $p''$, de sorte que l'union de ces deux
segments est encore une g\'eod\'esique.

Il est maintenant clair que, si $(C_0,C_1,\ldots ,C_n)$ est une
galerie g\'eod\'esique de $M(a)$ d'origine $(C,a)$, la suite des
demi-espaces ferm\'es $A(C_i)$ est strictement croissante. En
particulier, $A(C_n)$ contient strictement les murs $M(b(C_i))$, pour
$0\leq i<n$.

\smallskip 
Montrons maintenant que $M(b(C'))\cap M(b(C''))=\emptyset$,
lorsque $C'$ et $C''$ sont deux chambres de ${\cal B}^+(C)$ telles que
$C'\cap C''$ est une ar\^ete, et $\delta (C')= \delta (C'')$. Dans ce
cas, l'ar\^ete $a_-$ formant $C'\cap C''$ est oppos\'ee \`a $a(C')$ et
$a(C'')$ dans $C'$ et $C''$ respectivement. Il existe alors deux
ar\^etes $b'_-$ et $b''_-$ de $C'$ et $C''$, oppos\'ees \`a $b(C')$ et
$b(C'')$ respectivement, donc adjacentes \`a $a_-$ en un sommet $x$,
avec $M(b'_-)=M(b(C'))$ et $M(b''_-)=M(b(C''))$.

\begin{center}
\mbox{\relabelbox\small
\epsfysize 3cm\epsfbox{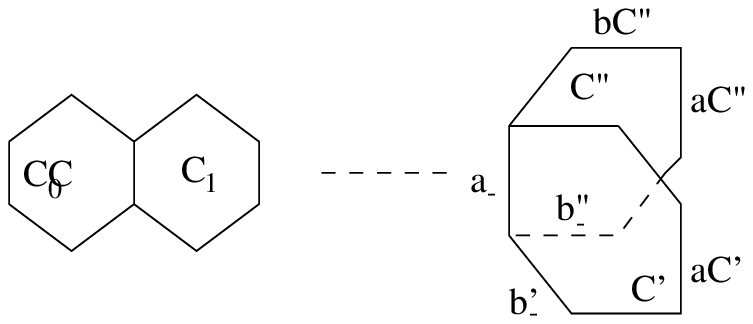}
\adjustrelabel <-2pt,0pt> {a}{$a_-$}
\relabel{b'}{$b'_-$}
\relabel{b"}{$b''_-$}
\relabel{C'}{$C'$}
\relabel{C"}{$C''$}
\relabel{C}{$C_1$}
\adjustrelabel <-2pt,0pt> {CC}{$C_0\!=\!C$}
\relabel{bC"}{$b(C'')$}
\relabel{aC"}{$a(C'')$}
\relabel{aC'}{$a(C')$}
\extralabel <-2.9cm, 0.7cm> $x$
\extralabel <-1.08cm, 1.2cm> $\longleftarrow b(C')$
\endrelabelbox}
\end{center}

Pour montrer que ces deux murs sont disjoints, on exhibe l\`a aussi
une perpendiculaire commune. Auparavant, on modifie la m\'etrique
CAT$(0)$ sur $P$, en rendant tous les polygones de $P$ r\'eguliers \`a
angle droit (donc hyperboliques, sauf si au d\'epart on avait des
carr\'es). La nouvelle m\'etrique est bien encore CAT$(0)$ (et m\^eme
souvent CAT$(-1)$), puisque tous les links de $P$ sont des graphes
bipartis ($P$ admet un type), donc ont des circuits de longueur au
moins $4$. Alors $b'_-\cup b''_-$ est g\'eod\'esique m\^eme en $x$, et
perpendiculaire aux deux murs.

Il est alors \'evident que $A(C')$ contient strictement la moiti\'e
ferm\'ee de $P$ d\'efinie par $M(b(C''))$ et contenant $a(C'')$
(i.e.~dont la r\'eunion avec $A(C'')$ est $P$ entier).

\begin{center}
\mbox{\relabelbox\small
\epsfysize 4cm\epsfbox{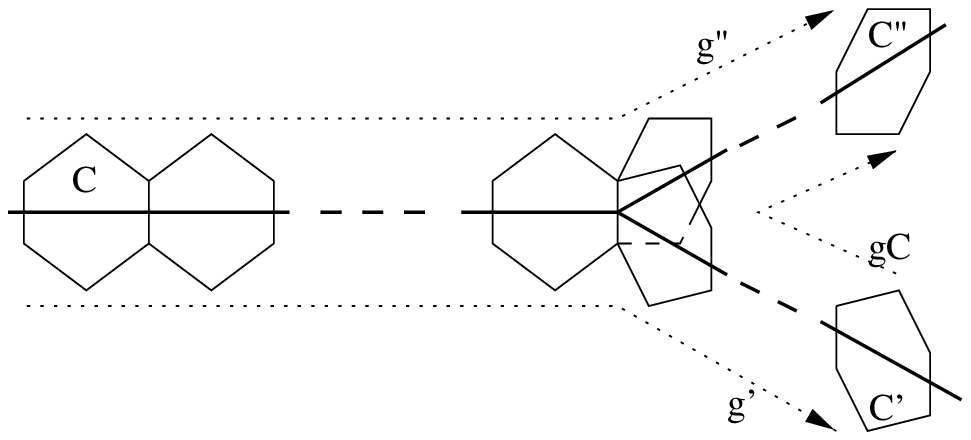}
\relabel{C'}{$C'$}
\relabel{C"}{$C''$}
\relabel{C}{$C$}
\relabel{g'}{$g'$}
\relabel{g"}{$g''$}
\adjustrelabel <0pt, 2pt> {gC}{$g(C',C'')$}
\endrelabelbox}
\end{center}

 Pour achever la preuve de l'affirmation, soient $C'$ et $C''$ deux
chambres distinctes de $B^+(C)$, telles que la galerie g\'eod\'esique
de $M(a)$ d'origine $(C,a)$ et d'extr\'emit\'e $C''$ ne passe pas par
$C'$.  Si $g'$ et $g''$ sont les galeries g\'eod\'esiques de $M(a)$
d'origine $(C,a)$ et d'extr\'emit\'e $C'$ et $C''$ respectivement, la
galerie $g(C',C'')$ obtenue \`a partir de $g'^{-1}.g''$ en \^otant les
r\'ep\'etitions permet, compte tenu des r\'esultats pr\'eliminaires
ci-dessus, de construire une suite strictement d\'ecroissante de
moiti\'es ferm\'ees dont la premi\`ere est $A(C')$ et la derni\`ere
est la moiti\'e compl\'ementaire de $A(C'')$.  Ceci conclut. 
\cqfd

Revenons \`a la preuve de la proposition. Pour $\overline{f}\in
\hbox{Fix}(C)\cap\hbox{Aut}_F(P)$, supposons avoir construit
$\overline{f}^n_b$ fixant strictement le mur $M(b)$, et co\"\i ncidant
avec $\overline{f}$ sur chaque chambre de ${\cal B}_n^+(C)$ (c'est
vrai pour $n=1$, d'apr\`es l'affirmation 1).  Alors
$(\overline{f}^n_b)^{-1}\circ \overline{f}$ agit trivialement sur
chaque chambre de ${\cal B}_n^+(C)$. Soient $C^1,\ldots ,C^k$ les
chambres de ${\cal B}^+(C)$ avec $\delta (C^i)=n$.  En appliquant
l'affirmation 1 \`a la chambre $C^1$, aux ar\^etes $a(C^1)$ et
$b(C^1)$, on trouve $\overline{f}^1$ fixant toute la moiti\'e $A(C^1)$
et co\"\i ncidant avec $(\overline{f}^n_b)^{-1}\circ\overline{f}$ sur
l'ensemble des chambres contenant $a(C^1)$. D'apr\`es l'affirmation 2,
$\overline{f}^1$ fixe strictement le mur $M(b)$, agit trivialement sur
toutes les chambres de ${\cal B}_n^+(C)$, et m\^eme sur les chambres
de ${\cal B}^+(C)$ adjacentes \`a $C^2,C^3,\ldots$ ou $C^k$.  Alors
$(\overline{f}^1)^{-1}\circ (\overline{f}^n_b)^{-1}\circ\overline{f} $
agit trivialement sur chaque chambre de ${\cal B}_n^+(C)$, et sur
chaque chambre contenant $a(C^1)$. En r\'eutilisant les affirmations 1
et 2, on trouve des automorphismes $\overline{f}^2,\overline{f}^3,
\ldots ,\overline{f}^k$ fixant tous strictement $M(b)$, tels que
$\overline{f}^{n+1}_b= \overline{f}^n_b\circ \overline{f}^1\circ
\overline{f}^2\circ\ldots\circ \overline{f}^k$ agit comme
$\overline{f}$ sur ${\cal B}_{n+1}^+(C)$. L'automorphisme
$\overline{f}^{n+1}_b$ fixe strictement $M(b)$ et co\"\i ncide avec
$\overline{f}$ sur chaque chambre de $B_{n+1}^+(C)$. En it\'erant ce
processus, et quitte \`a extraire une sous-suite convergente, on
trouve \`a la limite un $\overline{f}^+_b$ fixant strictement le mur
$M(b)$ et co\"\i ncidant avec $\overline{f}$ sur chaque chambre de
${\cal B}^+(C)$.  On peut imposer que $\overline{f}^+_b$ fixe point
par point la moiti\'e de $P$ d\'efinie par $M(b)$ et ne contenant pas
$a$.

En appliquant la construction pr\'ec\'edente sur l'autre moiti\'e de
$M(b)$, on trouve un $\overline{f}^-_b$ co\"\i ncidant avec
$\overline{f}$ sur chaque chambre de ${\cal B}^-(C)$, et fixant point
par point la moiti\'e de $P$ d\'efinie par $M(b)$ et contenant $a$. Si
on pose $\overline{f}_b=\overline{f}^+_b\circ \overline{f}^-_b$ et
$\overline{f}_a=(\overline{f}_b)^{-1}\circ \overline{f}$, on a
$\overline{f}=\overline{f}_b\circ \overline{f}_a$, avec
$\overline{f}_b$ fixant strictement le mur $M(b)$, et $\overline{f}_a$
fixant strictement le mur $M(a)$.  
\cqfd

Soit $k$ un entier pair au moins $4$ et $L$ un graphe fini de maille
au moins $5$ si $k=4$ et $4$ si $k\geq 6$.   Pour tout bloc $B$
de $A(k,L)$ (au sens de la d\'efinition \ref{def:facette_bloc}), notons
$F_B$ le sous-groupe caract\'eristique des automorphismes de $B$
engendr\'e par les fixateurs de facettes dans $B$. En fait, si $x$ est le 
centre du bloc $B$, alors $F_B=G_x$ avec les notations pr\'ec\'edant la
proposition \ref{prop:aut_zero_egal_aut_plus}. Remarquons que
$W(k,L)$ est un sous-groupe de Aut$_F A(k,L)$.

\blemm\label{lem:fixefacette}
Soit $B_0$ un bloc de $A(k,L)$. Si $\rho = \rho_{B_0}$ d\'esigne le 
morphisme de restriction de {\rm Stab}$\bigl(B_0,\hbox{\rm Aut}A(k,L)\bigr)$ 
dans {\rm Aut}$(B_0)$, alors 
$$\rho (\hbox{\rm Stab}\bigl(B_0,\hbox{\rm Aut}_F A(k,L)\bigr) ) = F_{B_0}.$$
\elemm

\dem 
Pour tout bloc $B$ de $A(k,L)$, notons ${\underline F} _{B}$
l'image r\'eciproque de $F_B$ par $\rho _B$. Par le prolongement $W(k,L)$-\'equivariant (voir les remarques avant la d\'efinition 
\ref{def:facette_bloc}), on a $\rho _B({\underline F} _{B}) = F_B$.

D'autre part, si ${\overline F} _{B}$ est le stabilisateur de $B$ dans 
Aut$_FA(k,L)$, alors  ${\underline F} _{B}\subset {\overline F} _{B}$. 
En effet, si $\hat \varphi\in {\underline F} _{B}$, par d\'efinition
$\rho _B (\hat \varphi)$ s'\'ecrit $\rho _B (\hat \varphi) = 
\varphi _1\cdots \varphi_n$, o\`u les
$\varphi_i$ sont des automorphismes de $B$ fixant une facette de $B$.
Comme $\rho _B({\underline F} _{B}) = F_B$, il existe 
$\hat \varphi_1, \cdots ,\hat \varphi_n$ \'el\'ements de
${\underline F} _{B}$ prolongeant les $\varphi_i$. On a donc
$\hat\varphi = \hat \varphi_1 \cdots \hat \varphi_n\varepsilon$, o\`u $\varepsilon$ vaut l'identit\'e sur
$B$. Chaque terme de la d\'ecomposition fixant une facette de $B$, on a $\hat\varphi\in $ Aut$_F A(k,L)$.

Pour montrer l'inclusion r\'eciproque ${\underline F} _{B_0}\supset 
{\overline F} _{B_0}$, introduisons le sous-groupe $H$ de $\hbox{Aut}A(k,L)$ 
engendr\'e par $W(k,L)$ et ${\underline F} _{B_0}$. Nous allons d'abord 
montrer que $H=$ Aut$_F A(k,L)$, puis que Stab$(B_0,H) = 
{\underline F} _{B_0}$, ce qui ach\`evera la preuve du lemme.

D'abord, comme ${\underline F} _{B_0}\subset {\overline F} _{B_0}$ 
et $W(k,L)\subset $ Aut$_F(A(k,L))$, on a bien $H\subset $ Aut$_F A(k,L)$. 
R\'eciproquement si $f$ est un automorphisme
de $A(k,L)$ fixant une facette $\phi$, il existe $w$ dans $W(k,L)$ tel que
$w(\phi)\subset B_0$. Alors
l'automorphisme $w f w^{-1}$ fixe une facette $\phi'$ de $B_0$. 
Si $s$ d\'esigne la r\'eflexion
de $W(k,L)$ par rapport au mur passant par cette facette, il existe un 
$k\in\lbrace 0,1\rbrace$
tel que $s^{k} w f w^{-1}$ pr\'eserve $B_0$ et fixe une 
facette $\phi'$ de $B_0$,
donc est dans ${\underline F} _{B_0}$.
Ainsi $f\in H$, et $H$ contient Aut$_F A(k,L)$.

Montrons maintenant que Stab$(B_0,H) = {\underline F} _{B_0}$, 
c'est-\`a-dire Stab$(B_0,H) \subset {\underline F} _{B_0}$.
Si nous v\'erifions que tout $h\in H$ peut s'\'ecrire $h=w f$, 
avec $w\in W(k,L)$ et $f\in {\underline F} _{B_0}$, alors on aura 
$h(B_0)=B_0$ implique $w=1$, donc $h\in {\underline F} _{B_0}$.
Pour \'etablir que $H$ co\"\i ncide avec l'ensemble $H'$ des 
automorphismes $f$ de $A(k,L)$ tels que $w_0^{-1} f\in
{\underline F} _{B_0}$, pour $w_0$ l'unique \'el\'ement de $W(k,L)$ 
tel que $w_0(B_0)=f(B_0)$,
introduisons l'ensemble $H''$ des automorphismes $f$ de $A(k,L)$ tels 
que, pour tout bloc $B$, $\ w^{-1} f\in
{\underline F} _{B}$, avec $w$ l'unique \'el\'ement de $W(k,L)$ tel 
que $w(B)=f(B)$. Il est clair que $H''$ est un sous-groupe de $H$ 
contenu dans $H'$. Pour conclure, montrons
que $H''=H$. Comme $W(k,L)\subset H''$, il suffit de montrer que 
${\underline F} _{B_0}\subset H''$, ce qui d\'ecoule
de l'affirmation suivante:
si $w_1^{-1} f$ est dans ${\underline F} _{B_1}$ et si $B_1\cap B_2$ 
est une facette $\phi$,
alors $w_2^{-1}\circ f$ est dans ${\underline F} _{B_2}$ (avec 
$w_i(B_i)=f(B_i)=B'_i$).
Pour voir ceci, soit $s$ (resp.~$s'$) la r\'eflexion de $W(k,L)$ 
\'echangeant $B_1$ et $B_2$
(resp.~$B'_1$ et $B'_2$), alors $w_2=s'w_1s$. Posons 
$\varepsilon _i=w_i^{-1} f$. Par hypoth\`ese
sur l'automorphisme $f$, on a $\varepsilon _1$ est dans 
${\underline F} _{B_1}$. Donc $s\varepsilon _1 s$ est dans 
${\underline F} _{B_2}$. Or $(s\varepsilon _1 s)^{-1} 
\varepsilon _2 = s  f^{-1}w_1 s s w_1^{-1} s' f =s  f^{-1}  s'  f $. 
Ce dernier automorphisme
fixe la facette $\phi$ et pr\'eserve le bloc $B_2$, donc est lui aussi 
dans ${\underline F} _{B_2}$, ce qui conclut. 
\cqfd

Appelons {\it facette de $L$} l'\'etoile d'un sommet de $L$ dans la
subdivision barycentrique $L'$.  Soit $F$ le sous-groupe
caract\'eristique de Aut$(L)$ engendr\'e par les fixateurs de facettes
de $L$. Si $L$ est le graphe biparti complet sur $p+q$ sommets avec 
$p,q\geq 3$, alors Aut$_0(L)(\simeq S_p\times S_q)= F$.

\bcoro\label{coro:fixefacette_indfini}
Le quotient de $\hbox{\rm Aut} A(k,L)$ par son sous-groupe distingu\'e
$\hbox{\rm Aut}_F A(k,L)$ est isomorphe au quotient de $\hbox{\rm Aut}(L)$ 
par son sous-groupe distingu\'e $F$.  
\ecoro

\dem 
Le groupe $\hbox{Aut} A(k,L)$ est transitif sur les sommets de
$A(k,L)$ car $W(k,L)$ l'est. Donc pour tout sommet $x_0$, centre du
bloc $B_0$, le quotient $\hbox{Aut} A(k,L)/\hbox{Aut}_F A(k,L)$ est
isomorphe \`a $\hbox{Fix~} x_0/\hbox{Fix~} x_0\cap\hbox{Aut}_F
A(k,L)$.  Par restriction, on a un morphisme $\hbox{Fix~} x_0 \ra
\hbox{Aut}(B_0)$, qui est surjectif par le paragraphe pr\'ec\'edant la
d\'efinition \ref{def:facette_bloc}. Son noyau est contenu dans
$\hbox{Fix~} x_0\cap\hbox{Aut}_F A(k,L)$. De plus, par le lemme pr\'ec\'edent,
l'image de $\hbox{Fix~} x_0\cap\hbox{Aut}_F A(k,L)$ est exactement $F_{B_0}$.
Donc $\hbox{Fix~} x_0/\hbox{Fix~} x_0\cap\hbox{Aut}_F
A(k,L)$ est isomorphe au quo\-tient $\hbox{Aut}(B_0)/F_{B_0}$. 
\cqfd

Soit $G$ un groupe de Chevalley fini de rang $2$, sur le corps fini $K$,
de syst\`eme de racines $\Phi$, de racines fondamentales 
$\alpha_1,\alpha_2$, de racines positives $\Phi^+$ et de groupes de racines 
$$X_\alpha=\{x_\alpha(t)\,/\, t\in K\}$$ pour $\alpha\in\Phi$. Nous
utiliserons les notations de \cite{Car}.  En particulier, $U$ est le
sous-groupe de $G$ engendr\'e par les racines positives. On a un
morphisme $h:\hbox{Hom}(\ZZ\alpha_1\oplus\ZZ\alpha_2,K^\times) \ra
\hbox{Aut} (G)$ qui, \`a un caract\`ere $\chi$ du r\'eseaux des
racines \`a valeurs dans le groupe multiplicatif de $K$, associe
l'automorphisme de $G$ induit par l'automorphisme
$$h(\chi):x_\alpha(t)\mapsto x_\alpha(\chi(\alpha)t)$$
sur chaque groupe de racine de $G$. On rappelle (voir \cite{Car}) que
$G$ est sans centre, est engendr\'e par les groupes de racines $X_\alpha$, 
et que chaque racine est combinaison lin\'eaire \`a c\oe fficients entiers 
(tous du m\^eme signe) de $\alpha_1,\alpha_2$.

On identifie $G$ \`a son image dans Aut$(G)$ par les automorphismes
int\'erieurs. On note $\widehat{H}$ l'image de $h$, $H=G\cap
\widehat{H}$ et $B=UH$. Il existe alors (voir \cite[page 101]{Car}) un
sous-groupe $N$ de $G$ tel que $(B,N)$ est une BN-paire de $G$. Soit
$L$ le $m$-gone g\'en\'eralis\'e associ\'e \`a cette BN-paire, muni de
son action de $G$, de sa chambre fondamentale $c$ de fixateur $B$, et
de son appartement fondamental $\Sigma$ de fixateur $H$ \cite[page
102]{Car}.  On identifie $\Phi$ avec l'ensemble des demi-appartements de
$\Sigma$, de sorte que $\Phi^+$ corresponde \`a ceux contenant $c$, et
que $X_\alpha$ soit le fixateur de la r\'eunion de $\alpha$ et des
ar\^etes de $L$ rencontrant $\alpha$ en un sommet int\'erieur de
$\alpha$.

Notons que $\widehat{H}$ pr\'eserve chaque groupe de racine
$X_\alpha$. Par cons\'equent, il agit sur $L$ en fixant $\Sigma$ (et
en particulier en pr\'eservant le type). Pour $i=1,2$, notons $x_i$ le
sommet de $c$ appartenant au bord de $\alpha_i$, $c_i$ la chambre de
$\Sigma$ adjacente \`a $c$ en $x_i$, et $\phi_i$ la facette de $L$ de
centre $x_i$. Les ar\^etes de $\phi_i$ sont les moiti\'es contenant
$x_i$ des chambres disjointes $\{c\}\cup\{x_{\alpha_i}(t)c_i\,/\, t\in
K\}$, car $X_{\alpha_i}$ agit simplement transitivement sur l'ensemble
des chambres contenant $x_i$ diff\'erentes de $c$. Tout caract\`ere
$\chi:\ZZ\alpha_1\oplus\ZZ\alpha_2\ra K^\times$ s'\'ecrit comme un produit
de caract\`eres $\chi_1\chi_2$ avec $\chi_i$ valant $1$ sur
$\alpha_i$. Comme $h$ est un morphisme, on a donc
$h(\chi)=h(\chi_1)h(\chi_2)$. De plus $h(\chi_i)$ fixe $\phi_i$
par la description pr\'ec\'edente.

\bprop\label{prop:f1f2_algebrique} 
Si $F$ est le sous-groupe caract\'eristique de $\hbox{\rm Aut}(L)$ 
engendr\'e par les fixateurs de facettes de $L$, alors
\begin{enumerate}
\item
$F=G\widehat{H}$,
\item
Le fixateur {\rm Fix}$_F(c)$ de $c$ dans $F$ est $U\widehat{H}$,
\item
$\hbox{\rm Aut}_0(L)/F$ est isomorphe au groupe $\hbox{\rm Aut}(K)$ des 
automorphismes du corps $K$.
\end{enumerate}
\eprop

\dem (1) D'apr\`es les rappels pr\'ec\'edents, l'inclusion de
$G\widehat{H}$ dans $F$ est claire. Soit $\phi$ une facette de $L$ et
$f$ un automorphisme de $L$ fixant $\phi$. On veut montrer que $f$
appartient \`a $G\widehat{H}$. Quitte \`a composer \`a gauche par un
\'el\'ement de $G$, on peut supposer que $f$ fixe $\Sigma$ et l'une
des facettes $\phi_1$ ou $\phi_2$, disons $\phi_1$.
Comme $f$ pr\'eserve $\phi_2$ en fixant $c$ et $c_2$, on peut \'ecrire
$$f(x_{\alpha_2}(1)c_2) = x_{\alpha_2}(\xi)c_2$$ pour un certain $\xi$
dans $K-\{0\}$. Soit $\chi$ le caract\`ere qui \`a $\alpha_1$ associe
$1$ et \`a $\alpha_2$ associe $\xi$. Montrons alors que $f=h(\chi)$.
Posons $\theta = h(\chi)^{-1} f$. C'est un automorphisme de $L$ fixant
$\Sigma, \phi_1$ et $x_{\alpha_2}(1)c_2$. Notons que $\theta$
normalise $G$, et notons encore $\theta $ l'automorphisme de $G$
induit. Alors $\theta$ pr\'eserve chaque $X_\alpha$,
$\theta(x_{\alpha_1}(t))=x_{\alpha_1}(t)$ pour tout $t\in K$ et
$\theta(x_{\alpha_2}(1))=x_{\alpha_2}(1)$.  Il d\'ecoule de la preuve
du th\'eor\`eme 12.5.1 de \cite[page 211]{Car} que l'ensemble des
automorphismes $\theta$ de $G$ qui pr\'eservent chaque $X_\alpha$, et fixent
$x_{\alpha_i}(1)$ pour $i=1,2$, est un sous-groupe de $G$
isomorphe \`a Aut$(K)$, et que si de plus
$\theta(x_{\alpha_1}(t))=x_{\alpha_1}(t)$ pour tout $t\in K$, alors
$\theta$ vaut l'identit\'e.

(2) L'inclusion de $U\widehat{H}$ dans Fix$_F(c)$ est claire.
R\'eciproquement, soit $f$ dans $F$ fixant $c$. Alors $f=g\widehat{f}$
avec $g\in G$ et $\widehat{f}\in\widehat{H}$ par (1). Comme
$\widehat{f}$ fixe $c$, on en d\'eduit que $g$ fixe $c$. Or le
fixateur de $c$ dans $G$ est $B=UH$, et comme $H\subset \widehat{H}$,
le r\'esultat en d\'ecoule.

(3) Il est clair que $F$ est contenu et distingu\'e dans Aut$_0(L)$.
Soit $\theta\in \hbox{Aut}_0(L)$. Quitte \`a le multiplier par un
\'el\'ement de $G$, on peut supposer que $\theta$ fixe $\Sigma$.
Quitte \`a le multiplier par un \'el\'ement de $\widehat{H}$, on peut
supposer que $\theta$ fixe $x_{\alpha_i}(1)$ pour $i=1,2$. Soit $Z$ le
fixateur dans $\hbox{Aut}_0(L)$ de $\Sigma\cup
\{x_{\alpha_1}(1),x_{\alpha_2}(1)\}$.  On a donc un isomorphisme entre
$\hbox{\rm Aut}_0(L)/F$ et $Z/Z\cap F$.  Or si $f=g\widehat{f}$ fixe
$\Sigma$, avec $g\in G$ et $\widehat{f}\in\widehat{H}$, alors $g$ fixe
$\Sigma$. Donc $g$ appartient au fixateur de $\Sigma$ dans $G$, qui
est $H$. Par cons\'equent $f\in\widehat{H}$. Or un \'el\'ement
$h(\chi)$ de $\widehat{H}$ fixant $x_{\alpha_1}(1)$ et
$x_{\alpha_2}(1)$ vaut l'identit\'e, car on aurait $\chi(\alpha_1)=1$
et $\chi(\alpha_2)=1$. D'o\`u $Z\cap F=\{1\}$, ce qui montre le
r\'esultat, $Z$ \'etant isomorphe \`a Aut$(K)$, d'apr\`es le dernier 
argument de (1).  \cqfd

\bcoro\label{coro:F1F2quoi} Pour $i=1,2$, soit $F_i$ le fixateur de la
facette $\phi_i$ dans {\rm Aut}$(L)$.  Alors $$\hbox{\rm
Fix}_F(c)=F_1F_2=F_2F_1.$$ \ecoro

\dem L'\'egalit\'e $F_1F_2=F_2F_1$ vient du fait que $F_1$ et $F_2$
fixent $c=\phi_1\cap\phi_2$ donc $F_i$ pr\'eserve $\phi_{3-i}$.
L'inclusion de $F_1F_2$ dans $\hbox{\rm Fix}_F(c)$ est claire.  Pour
montrer l'inclusion inverse, comme $\hbox{\rm Fix}_F(c)=U\widehat{H}$,
il suffit de le faire pour $U$ et pour $\widehat{H}$. Or $U$ est
engendr\'e par les $X_\alpha$ pour $\alpha$ racine positive, et un tel
$X_\alpha$ est contenu soit dans $F_1$, soit dans $F_2$.  De plus, on
a vu avant la proposition \ref{prop:f1f2_algebrique} que pour tout
caract\`ere $\chi$, $h(\chi)=h(\chi_1)h(\chi_2)$ avec $h(\chi_i)$
fixant $\phi_i$.  \cqfd

\bcoro\label{coro:plus_implique_type} 
Si $L$ est un $m$--gone g\'en\'eralis\'e \'epais fini classique, alors
le grou\-pe {\rm Aut}$_0 A(k,L)\cap\hbox{\rm Aut}_F (P)$ des $F$-automorphismes
pr\'eservant le type de $A(k,L)$, co\"{\i}n\-ci\-de avec le groupe
{\rm Aut}$^+A(k,L)$ des automorphismes de $A(k,L)$ engendr\'e par les
fixateurs stricts de murs propres, et est distingu\'e dans\break {\rm Aut}$_0
A(k,L)$, de quotient trivial si $m=2$, et sinon isomorphe au groupe fini 
des automorphismes de corps du corps fini de d\'efinition de $L$.
 \ecoro

\dem Nous allons v\'erifier les hypoth\`eses de la proposition
\ref{prop:aut_zero_egal_aut_plus} (2) et (3) pour montrer que
Aut$^+A(k,L)=G_0 =$ Aut$_0(P)\cap\hbox{Aut}_F (P)$.  Puisque $A(k,L)$
est un immeuble, il admet un type. Puisque $A(k,L)$ est la
r\'ealisation g\'eom\'et\-rique de Davis--Moussong d'un syst\`eme de
Coxeter, il v\'erifie la propri\'et\'e $(P^+)$. Tous ses murs sont
propres par le lemme \ref{lemm:tout_mur_propre}.  Par hypoth\`ese, $L$
est \'epais et de Moufang, et pour tout sommet $x_0$ de $A(k,L)$, le
bord du bloc de centre $x_0$ s'identifie avec $L$, donc l'hypoth\`ese
de \ref{prop:aut_zero_egal_aut_plus} (2) est v\'erifi\'ee.

Par le lemme \ref{lem:fixefacette}, si $\rho$ est le morphisme de
restriction \`a un bloc $B_0$ de centre $x_0$ des automorphismes de
$A(k,L)$ fixant $x_0$, alors l'image $I$ par $\rho$ du fixateur de
$x_0$ dans Aut$_F A(k,L)$ est exactement $G_{x_0}$. L'inclusion de
$G_{x_0}$ dans $I$ montre la propri\'et\'e $(P_0)$ et l'inclusion
r\'eciproque montre la deuxi\`eme hypoth\`ese de
\ref{prop:aut_zero_egal_aut_plus} (3).

La premi\`ere hypoth\`ese de
\ref{prop:aut_zero_egal_aut_plus} (3) d\'ecoule du corollaire
\ref{coro:F1F2quoi} si $m\geq 3$, et est claire si $m=2$.

Enfin, par le corollaire \ref{coro:fixefacette_indfini} et la
proposition \ref{prop:f1f2_algebrique} (3) si $m\geq 3$, le quotient
Aut$_0A(k,L)/\hbox{Aut}^+A(k,L)$ est isomorphe \`a Aut$(K)$.
\cqfd

\section{Simplicit\'e de groupes d'automorphismes d'espaces \`a murs}
\label{sect:simple}

\btheo\label{theo:simplicite}
Soient $(X,\M)$ un espace \`a murs hyperbolique, de graphe associ\'e
$\G$, et $G$ un groupe d'au\-to\-mor\-phis\-mes de $(X,\M)$, dont
l'action sur $\G$ est non \'el\'ementaire, d'ensemble limite \'egal
\`a $\partial \G$.  Supposons que $G$ v\'erifie la condition {\rm(P)}.
Soit $G^+$ le sous-groupe de $G$ engendr\'e par les fixateurs stricts
de murs propres et $H$ un sous-groupe distingu\'e de $G^+$.  Alors ou
bien $H$ est contenu dans le noyau de l'action de $G^+$ sur $\partial
\G$, ou bien $H$ est \'egal \`a $G^+$.  
\etheo

\bcoro
Si l'action de $G^+$ sur $\partial \G$ est fid\`ele, alors $G^+$ 
est simple. \hfill$\Box$
\ecoro

Remarquons que le sous-groupe $G^+$ est distingu\'e dans $G$, 
et qu'il peut \^etre trivial.

\dem
Soit $H$ un sous-groupe distingu\'e non trivial de $G^+$. Supposons
que $H$ n'est pas contenu dans le noyau de l'action de $G^+$ sur
$\partial \G$. Rappelons que $X$ est le sous-ensemble des sommets de
$\G$.

\blemm\label{lem:beaucoup_d_hyperbolique}
Pour tout demi-espace $A$ avec $\partial X\setminus \partial A$ non
vide, il existe une cha\^{\i}ne propre $(A_i)_{i\in\ZZ}$ et un \'el\'ement
$h$ dans $H$ tels que $A\subset A_0\setminus A_1$, $h(A_i)=A_{i+1}$
pour tout $i$.
\elemm

\begin{figure}[htbp]
\cl{\relabelbox\small
\epsfysize 6cm\epsfbox{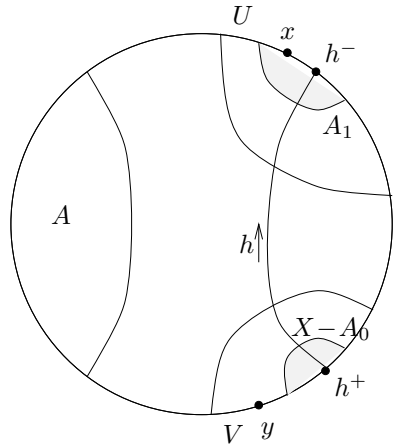}
\relabel {A}{$A$}
\adjustrelabel <-2pt, -2pt> {A1}{$A_1$}
\relabel {U}{$U$}
\relabel {V}{$V$}
\adjustrelabel <-3pt, 1pt> {X-A}{$X\!-\!A_0$}
\adjustrelabel <0pt, -2pt> {x}{$x$}
\relabel {y}{$y$}
\adjustrelabel <-1pt, 0pt> {h}{$h$}
\adjustrelabel <-1pt, 0pt> {hp}{$h^+$}
\adjustrelabel <-1pt, 0pt> {hm}{$h^-$}
\endrelabelbox}
\caption{Construction de cha\^{\i}ne invariante par un 
\'el\'ement hyperbolique}\label{fig3_2}
\end{figure}

\dem 
Puisque $\Lambda G=\partial \G$ n'a pas de point isol\'e ($G$ est
non \'el\'ementaire), et par la condition (H), il existe (voir figure
\ref{fig3_2}):
\begin{itemize}
\item 
$x,y$ deux points distincts dans l'ouvert 
$\partial X\setminus \partial A
=(\overline{X}\setminus \overline A)\cap \partial X$ de $\partial X$,
\item 
$U,V$ deux ouverts disjoints de $\overline X$, contenus dans
$\overline{X}\setminus \overline{A}$ et contenant respectivement
$x,y$,
\item 
$A_1$ un demi-espace contenu dans $U$, avec $\overline{A_1}$ 
un voisinage de $x$, et dont le mur est propre.
\end{itemize}

Par une application double du lemme \ref{prop:non_elementaire} (\`a
$G^+\subset G$ et \`a $H\subset G^+$), les couples des points fixes
d'\'el\'ements hyperboliques de $H$ sont denses dans $\partial^2
\G$. Soit donc $h$ un \'el\'ement hyperbolique de $H$ dont un point
fixe au bord est contenu dans l'int\'erieur de $\partial A_1$, et
l'autre dans $V$. Quitte \`a remplacer $h$ par une puissance
suffisamment grande (en valeur absolue), pour que $h(A_1)$ soit
strictement contenu dans $A_1$ et que $h^{-1}(X\setminus A_1)$ soit
contenu dans $V$, la suite de demi-espaces $(h^{i-1}(A_1))_{i\in\ZZ}$
est une cha\^{\i}ne. Le lemme est alors facile \`a v\'erifier.
\endproof

\blemm \label{lem:commutateur_diabolique} 
Soient $h\in H$ et $C=(A_i)_{i\in\ZZ}$ une cha\^{\i}ne propre tels que 
$h(A_i)=A_{i+1}$ pour tout $i$.  Pour tout $g\in G$ fixant 
strictement $C$, il existe $f\in G^+$ tel que $g=[h,f]$.  
\elemm

\dem
On note $[u,v]=uvu^{-1}v^{-1}$. Soient $h,g$ comme dans l'\'enonc\'e.
Si $u\in G$ fixe strictement les $M_i=\{A_i,X\setminus A_i\}$, notons
$u_i$ la restriction de $u$ \`a $A_i\setminus A_{i+1}$.  Alors
$g=[h,f]$ si et seulement si, pour tout $i\in \ZZ$,
$$g_i=h f_{i-1}h^{-1} f_i^{-1}$$
ou encore
$$f_i=g_i^{-1} h f_{i-1}h^{-1}.$$

Posons $f_0$ la restriction \`a $A_0\setminus A_1$ de l'identit\'e de
$G$.  Alors la relation de r\'ecurrence ci-dessus (ou
$f_{i-1}=h^{-1}g_i f_i h$ pour les $i$ strictement n\'egatifs) permet
de d\'efinir une application $f_i$ sur $A_i\setminus A_{i+1}$, qui est
par r\'ecurrence restriction \`a $A_i\setminus A_{i+1}$ d'un
\'el\'ement $\tilde{f}_i$ de $G$ fixant strictement $C$.  En effet, le
fixateur strict de $C$, qui contient $\tilde{f}_{i-1}$, est
distingu\'e dans le stabilisateur de $C$ (qui contient $h$).  Par la
propri\'et\'e (P), il existe un \'el\'ement $f$ dans $G$ fixant
strictement $C$, dont les restrictions sont les $f_i$, et la remarque
pr\'eliminaire montre que $g=[h,f]$.

Par d\'efinition, le fixateur strict dans $G$ d'une cha\^{\i}ne propre
est contenu dans $G^+$. Ceci conclut la preuve.  
\endproof

\bcoro \label{coro:contient_fixateur_demiespace}
Le groupe $H$ contient le fixateur strict dans $G$ de tout mur propre.
\ecoro

\dem 
Soit $M=\{A,X\setminus A\}$ un mur propre, donc tel que $\partial
X\setminus\partial A$ est non vide. Soit $g\in G$ fixant strictement
$M$. Par le lemme \ref{fixateur_mur_stabilisateur_demiespace}, pour
montrer que $g$ appartient \`a $H$, il suffit de le montrer en
supposant de plus que $g$ fixe (point par point) $X\setminus A$.  Par
le lemme \ref{lem:beaucoup_d_hyperbolique}, il existe $h\in H$ et une
cha\^{\i}ne propre $C=(A_i)_{i\in\ZZ}$ tels que $h(A_i)=A_{i+1}$ et $A
\subset A_0\setminus A_1$. En particulier $g$ fixe strictement
$C$. Par le lemme \ref{lem:commutateur_diabolique}, il existe $f$ dans
$G^+$ tel que $g=[h,f]= h(fh^{-1}f^{-1})$. Comme $H$ est distingu\'e
dans $G^+$, il contient $g$, d'o\`u le r\'esultat.
\endproof

Le corollaire \ref{coro:contient_fixateur_demiespace}  d\'emontre le 
th\'eor\`eme. \endproof

\section{Applications}
\label{sect:applications}

\btheo Soit $P$ un complexe poly\'edral pair {\rm CAT}$(0)$, dont la
m\'etrique est hyperbolique au sens de Gromov, dont le groupe des
automorphismes est non \'el\'ementaire et d'ensemble limite \'egal \`a
$\partial P$.  Soit {\rm Aut}$^+(P)$ le sous-groupe de {\rm Aut}$(P)$
engendr\'e par les fixateurs stricts de murs propres et $H$ un
sous-groupe distingu\'e de {\rm Aut}$^+(P)$.  Alors ou bien $H$ est contenu
dans le noyau de l'action de $G$ sur $\partial \G$, ou bien $H$ est
\'egal \`a {\rm Aut}$^+(P)$.  \etheo

\dem 
D'apr\`es le th\'eor\`eme \ref{theo:resume_prop_complexe_pair},
l'espace \`a murs $(X_P,\M_P)$ associ\'e \`a $P$ est un espace \`a
murs hyperbolique, et le bord de $P$ s'identifie au bord du graphe
associ\'e \`a $(X_P,\M_P)$. D'apr\`es le th\'eor\`eme
\ref{theo:meme_groupe_automorphisme}, le groupe des automorphismes de
$P$ (resp.~le groupe engendr\'e par les fixateurs stricts de murs
propres de $P$) co\"{\i}ncide avec le groupe $G$ des automorphismes de
l'espace \`a murs $(X_P,\M_P)$ (resp.~le groupe engendr\'e par les
fixateurs stricts de murs propres de $(X_P,\M_P)$).  Par le lemme
\ref{lemm:espace_amur_associe_verifie_Mprime}, l'espace \`a murs
$(X_P,\M_P)$ v\'erifie la condition (${\rm M}'$). Donc $G$ v\'erifie la
condition (P) par le lemme \ref{Mprime_entraine_P}. Le r\'esultat
d\'ecoule alors du th\'eor\`eme \ref{theo:simplicite}.  
\endproof

\bcoro \label{coro:final}
Sous les hypoth\`eses du th\'eor\`eme pr\'ec\'edent:
\begin{enumerate}
\item
Si $P$ est localement compact alors $H$ est 
relativement compact, ou \'egal \`a {\rm Aut}$^+ P$.
\item
Si le seul \'el\'ement de {\rm Aut}$^+\,P$ agissant trivialement sur
le bord de $P$ est l'identit\'e, alors {\rm Aut}$^+ P$ est simple.
\item
Si $P$ est {\rm CAT}$(-1)$ et tout point de $P$ est
contenu dans une droite g\'eo\-d\'e\-si\-que, alors {\rm Aut}$^+ P$ est simple.
\end{enumerate}
\ecoro

\dem 
Si $H$ est contenu dans le noyau de l'action sur le bord, et si $P$ est 
localement compact, alors par le lemme
\ref{lemm:compact_fidele}, $H$ est relativement compact. Sinon, par 
le th\'eor\`eme pr\'ec\'edent 
on a $H=$Aut$^+ P$, ce qui montre (1) et (2). L'assertion (3) d\'ecoule 
de (2) par le lemme \ref{lemm:compact_fidele}, car puisque Aut$(P)$
est non \'el\'ementaire, d'ensemble limite \'egal \`a tout $\partial \G$,
il n'y a pas de point isol\'e dans $\partial \G$.
\endproof

Le th\'eor\`eme \ref{theo:intro_pair} de l'introduction d\'ecoule de 
ce corollaire et de la remarque pr\'ec\'edant le lemme 
\ref{lemm:compact_fidele}.
 
\bcoro \label{coro:simple_Coxeter} 
Soit $(W,S)$ un syst\`eme de Coxeter,
avec $W$ hyperbolique au sens de Gromov. Alors le quotient, par 
son sous-groupe distingu\'e localement compact form\'e
des \'el\'ements fixant l'infini, du sous-groupe $G^+$ des
automorphismes du graphe de Cayley de $(W,S)$ engendr\'e par les
fixateurs stricts de murs propres, est simple.  Il est non trivial (et
donc non d\'enombrable) si et seulement si $(W,S)$ est non rigide.
\ecoro

\dem 
D'apr\`es la remarque (1) de la section
\ref{sect:espace_amur_associe_pair}, le complexe poly\'edral $|W|_0$ 
est pair. Il est localement compact, et $W$ agit discr\`etement avec 
quotient compact sur lui. Le r\'esultat de simplicit\'e d\'ecoule 
du th\'eor\`eme 
pr\'ec\'edent. La derni\`ere assertion d\'ecoule du 
th\'eor\`eme \ref{theo:coxeter_non_rigide}, la non trivialit\'e 
de Aut$^+|W|_{0}$ impliquant sa non d\'enombrabilit\'e par le lemme 
\ref{lem:nontrivial_nondenombrable}. 
\endproof

Le th\'eor\`eme \ref{theo:intro_coxeter} de l'introduction d\'ecoule de ce
corollaire, car le groupe des automorphismes du graphe de Cayley de $(W,S)$
s'identifie avec le groupe des automorphismes poly\'edraux de la r\'ealisation
g\'eom\'etrique au sens de Davis--Moussong de $(W,S)$ (voir section
\ref{sect:groupe_auto_polyedre_pair}).

Pour terminer, d\'emontrons les th\'eor\`emes \ref{theo:intro_bourdon}
 et \ref{theo:intro_Akl} de l'introduction. Par le lemme
 \ref{coro:plus_implique_type}, le groupe  des $F$-automorphismes
 pr\'eservant le type des immeubles hyperboliques $A(k,L)$
 co\"{\i}ncide avec le groupe engendr\'e par les fixateurs strict de
 murs propres, est d'indice fini dans Aut$_0A(k,L)$ et est simple par
 le corollaire \ref{coro:final} (3).  Comme $L$ est non rigide (par
 exemple si $m\geq 3$, un groupe de racine est non trivial et fixe
 l'\'etoile d'un sommet), il est non d\'enombrable, par le lemme
 \ref{lem:nontrivial_nondenombrable}.  Il est \'evidemment ferm\'e
 dans le groupe de tous les automorphismes, donc est localement
 compact.

Enfin, pour montrer que Aut$^+A(k,L)$ est non lin\'eaire, il suffit, par le
th\'eor\`eme de Schur--Kaplansky (voir par exemple \cite[page
154]{Dix}), de montrer qu'il contient un sous-groupe de type fini, de
torsion et infini. Supposons que $k$ est multiple de $4$ et que $L$ est
ou bien un graphe biparti complet $K_{p,p'}$, ou l'immeuble
sph\'erique d'un groupe de Chevalley fini sur un corps $F_q$ de
caract\'eristique $p$ diff\'erente de $2$.  En utilisant les m\'ethodes
de l'affirmation 2 de la proposition
\ref{prop:aut_zero_egal_aut_plus}, il est alors possible de montrer
que $G$ contient une copie du $p$--groupe infini \`a deux
g\'en\'erateurs $\tau,\alpha$ de Grigorchuk--Gupta--Sidki (voir
\cite[page 19]{Bau}).

\noindent{\small\parskip 0pt  
Laboratoire de Topologie et Dynamique URA 1169 CNRS\\Universit\'e Paris-Sud\\
B\^at. 425 (Math\'ematiques)\\91405 ORSAY Cedex\\FRANCE\\\smallskip\\
Email:\qua{\tt  haglund@math.u-psud.fr, Frederic.Paulin@math.u-psud.fr}}

\recd


\begin{thebibliography}


\bibitem{And}  
{\bf E\,M Andreev}, {\it On convex polyhedra in Lobachevskii space}, 
Math. USSR Sb. {12} (1970) 413--440

\bibitem{Bau}
{\bf G~Baumslag}, {\it Topics in Combinatorial Group Theory}, 
Lectures in Math. Birk\-hauser (1993)

\bibitem{BB}
{\bf W~Ballmann}, {\bf M~Brin}, 
{\it Polygonal complexes and combinatorial group theory}, 
Geom. Dedicata  50 S\'erie I (1994) 165-191

\bibitem{Ben}
{\bf N~Benakli}, {\it Polygonal complexes I: Combinatorial and geometric 
properties}, 
J. Pure Appl. Alg. {97} (1994) 247--263
 
\bibitem{Bourb}
{\bf N~Bourbaki}, {\it Groupes et alg\`ebres de Lie}, 
chap.~4,5,6, Hermann, Paris (1968)
  
\bibitem{Bourd1}
{\bf M~Bourdon}, {\it Structure conforme au bord et flot g\'eod\'esique
d'un CAT$(-1)$ espace},  L'Ens. Math. {41} (1995) 63--102

\bibitem{Bourd2}
{\bf M~Bourdon}, {\it Immeubles hyperboliques, dimension conforme et 
rigidit\'e de Mostow}, GAFA {7} (1997) 245--268

\bibitem{Bri}
{\bf M\,R~Bridson}, {\it Geodesics and curvature in metric 
sim\-pli\-cial com\-ple\-xes}, dans: ``Group 
theory from a geometrical viewpoint'' (E~Ghys, A~Haefliger, A~Verjovsky,
editeurs) World Scientific
(1991) 373--463

\bibitem{BH}
{\bf M\,R~Bridson}, {\bf A~Haefliger}, {\it Metric spaces of 
non-positive curvature},
Grund.~math.~Wiss. 319, Springer--Verlag (1998)

\bibitem{Bro} 
{\bf K~Brown}, {\it Buildings}, Springer--Verlag (1989)
 
\bibitem{BM} 
{\bf M~Burger}, {\bf S~Mozes}, {\it Finitely presented simple groups and 
products of trees}, 
C. R. Acad. Sci. Paris, Ser. I, {324} (1997) 747--752

\bibitem{Car}
{\bf R\,W Carter}, {\it Simple groups of Lie type}, 
Pure Appl, Math. {28}, Wiley (1972)
 
\bibitem{Cha}
{\bf C~Champetier},  {\it L'espace des groupes de type fini}, 
\`a para\^{\i}tre dans Topology

\bibitem{CD} 
{\bf R~Charney}, {\bf M~Davis}, {\it Singular metrics of
nonpositive curvature on branched covers of Riemannian manifolds},
Amer.~J.~Math. {115} (1993) 929--1009

\bibitem{Dav}
{\bf  M~Davis},  {\it Groups generated by reflections and aspherical manifolds 
not covered by Euclidean space}, Annals of Math. {117} (1983) 293--324

\bibitem{Dav3} 
{\bf M~Davis}, {\it Negative curvature and reflection groups}, to appear in
Handbook of Geometric Topology

\bibitem{DJS} 
{\bf M~Davis}, {\bf T~Januszkiewicz}, {\bf R~Scott}, 
{\it Nonpositive curvature of blow-ups},  to appear in Selecta Math.


\bibitem{Dix}
{\bf J~Dixon}, {\it The structure of linear groups}, 
Van Nostrand (1971)

\bibitem{GP}
{\bf D~Gaboriau}, {\bf F~Paulin}, {\it Sur les immeubles hyperboliques}, 
pr\'epublication  18, Univ.~Orsay (Janv 1998)

\bibitem{GH} {\bf E~Ghys}, {\bf P~de\,La\,Harpe}, editeurs, {\it Sur les
groupes hyperboliques d'apr\`es Mikhael Gromov}, Prog. in Math. {83},
Birkh\"auser (1990)

\bibitem{Gro}
{\bf M~Gromov}, {\it Hyperbolic groups}, dans: ``Essays in group theory'', 
(S~Gersten, editeur)  MSRI Pub. {8}, Springer--Verlag (1987) 75--263

\bibitem{Hae} 
{\bf A~Haefliger}, {\it Complexes of groups and orbihedra}, dans: ``Group 
theory from a geometrical viewpoint'' 
(E~Ghys, A~Haefliger, A~Verjovsky, editeurs) 
World Scientific (1991) 504--540

\bibitem{Hag}
{\bf F~Haglund}, {\it R\'eseaux de Coxeter--Davis et commensura\-teurs}, 
Ann. Inst. Fourier  48 (1998) 649-666
 
\bibitem{Har}
{\bf P~de\,la\,Harpe}, {\it An invitation to Coxeter groups}, dans: ``Group
theory from a geometrical viewpoint'',
(E~Ghys, A~Haefliger, A~Verjovsky, editeurs)
World Scientific (1991) 193--253

\bibitem{Mou} 
{\bf G~Moussong}, {\it Hyperbolic Coxeter group},
Doctoral Dissertation, Ohio State University (1988)

\bibitem{NR} 
{\bf G~Niblo}, {\bf L~Reeves}, {\it Groups acting on CAT$(0)$ cube complexes}, 
Geometry and Topology {1} (1997) 1--7

\bibitem{Ron} 
{\bf M\,A~Ronan}, {\it Lectures on buildings}, 
Persp. Math. {7}, Academic Press (1989)

\bibitem{Sag}
{\bf M~Sageev}, {\it Ends of group pairs and non-positively curved
cube complex},
Proc. London Math. Soc. {71} (1995) 585--617

\bibitem{Ser}
{\bf J-P~Serre}, {\it Arbres, amalgames, SL$_2$}, 
Ast\'erisque  {46}, Soc. Math. France (1983)

\bibitem{Whi} 
{\bf E\,H~Spanier}, {\it Algebraic topology}, Tata--McGraw--Hill (1966)

\bibitem{Tit}
{\bf J~Tits}, {\it Sur le groupe des automorphismes d'un arbre},
dans: ``Essays on Topology (M\'emoires d\'edi\'ees \`a G~de\,Rham)'', 
 Springer--Verlag (1970) 188--211


\end{thebibliography}
\end{document}